\input amstex
\input amsppt.sty
\NoBlackBoxes
\nologo
\topmatter
\title Lectures on Factorization of Birational Maps
\endtitle
\dedicatory in memory of Prof. Kunihiko Kodaira
\enddedicatory
\author Kenji Matsuki
\endauthor
\endtopmatter
\leftheadtext{}
\rightheadtext{}
\document

This is an expanded version of the notes for the lectures given by the author at RIMS in the
summer of 1999 to provide a detailed account of the proof for the (weak) factorization theorem of
birational maps by Abramovich-Karu-Matsuki-W{\l}odarczyk [1].  All the main ideas of these
lecture notes, therefore, are derived from the collaboration and joint work with D. Abramovich,
K. Karu and J. W{\l}odarczyk.  The author would like to emphasize that the theory of ``Birational
Cobordism", which plays the central role in our proof, was developed in W{\l}odarczyk [2] and
that its origin may be traced back to the earlier work of Morelli [1][2] in his brilliant
solution to the (strong) factorization conjecture for toric birational maps.  The author would
also like to emphasize that there is another independent proof for the (weak) factorization
theorem by W{\l}odarczyk [3], which is NOT discussed in these notes.  It is noteworthy that not
only the ideas of but also some essential steps of the two proofs rely on the work of Morelli:
our proof in Abramovich-Karu-Matsuki-W{\l}odarczyk [1] reduces the factorization of general
birational maps to that of toroidal birational maps and then uses the combinatorial algorithm by
Morelli to factor toroidal ones.  The proof in W{\l}odarczyk [2] uses the process of
``$\pi$-desingularization", which is the most difficult and subtle part of the algorithm by
Morelli.  In fact, all that we cover in these notes is the above mentioned reduction step
``from general to toroidal", taking the combinatorial algorithm for the factorization of toroidal
birational maps as a black box.  We refer the reader to Abramovich-Matsuki-Rashid [1] for what
is inside of the black box, including the details of the process of $\pi$-desingularization and
clarifications of some discrepancies contained in the original arguments of Morelli [1][2].

The deepest gratitude of the author goes to Prof. S. Mori for his generosity to provide an
opportunity to deliver these lectures at RIMS and to Profs. Y. Miyaoka,  N. Nakayama and the
participants of the seminar, who listened to and suggested innumerable improvements to the
presentation and contents of the lectures.  In the talks subsequently given at Tokyo University,
another set of invaluable comments was given to the author by Profs. S. Iitaka, Y. Kawamata and
others.  The present notes are made more comprehensible due to their constructive criticism than the
original ones.  The author hopes that our sweat in the midst of the humidity and heat of ``TSUYU"
(the rainy season in Japan) served some purpose.

\newpage

$$\bold{CONTENTS}$$

Chapter 0. Introduction

\vskip.1in

\ \ \S 0-1. Outline of the strategy

\ \ \S 0-2. Some historical remarks

\vskip.1in

Chapter 1. Preliminaries

\vskip.1in

\ \ \S 1-1. Brief review of Geometric Invariant Theory and Toric Geometry

\ \ \S 1-2. ``Locally Toric" structures vs. ``Toroidal" structures

\ \ \S 1-3. $K^*$-action on locally toric and toroidal structures

\ \ \S 1-4. Elimination of points of indeterminacy

\vskip.1in

Chapter 2. Birational Cobordism

\vskip.1in

\ \ \S 2-1. Definition of a birational cobordism and the toric main example

\ \ \S 2-2. Construction of a (collapsible) birational cobordism

\ \ \S 2-3. Interpretation by Geometric Invariant Theory

\ \ \S 2-4. Factorization into locally toric birational maps

\vskip.1in

Chapter 3. Torification

\vskip.1in

\ \ \S 3-1. Construction of the torific ideal

\ \ \S 3-2. The torifying property of the torific ideal

\vskip.1in

Chapter 4. Recovery of Nonsingularity

\vskip.1in

\ \ \S 4-1. Application of the canonical resolution of singularities

\ \ \S 4-2. Conclusion of the proof for the weak factorization theorem

\vskip.1in

Chapter 5. Generalizations

\vskip.1in

\ \ \S 5-1. Factorization of bimeromorphic maps

\ \ \S 5-2. Equivariant factorization under group action

\ \ \S 5-3. Factorization over a non-algebraically closed field

\ \ \S 5-4. Factorization in the logarithmic category

\vskip.1in

Chapter 6. Problems

\vskip.1in

\ \ \S 6-1. Effectiveness of the construction

\ \ \S 6-2. Toroidalization Conjecture and Strong Factorization Conjecture

\vskip.2in

$$\bold{CHAPTER\ 0.\ INTRODUCTION}$$

\vskip.1in

The purpose of these lecture notes is to provide a detailed account of the proof given by
Abramovich-Karu-Matsuki-W{\l}odarczyk [1] for the following (weak) factorization conjecture for
birational maps.  (We note that there is another proof independently given by W{\l}odarczyk
[1], which we will NOT discuss in these notes.)

\proclaim{Weak Factorizatioin Theorem of Birational Maps} Let $\phi:X_1
\dashrightarrow X_2$ be a birational map between complete nonsingular varieties over an
algebraically closed field $K$ of characteristic zero.  Let $X_1 \supset U \subset X_2$ be a
common open subset over which $\phi$ is an isomorphism.  Then $\phi$ can be factored into a
sequence of blowups and blowdowns with smooth irreducible centers disjoint from $U$.  That is to
say, there exists a sequence of birational maps between complete nonsingular varieties
$$X_1 = V_1 \overset{\psi_1}\to{\dashrightarrow} V_2 \overset{\psi_2}\to{\dashrightarrow}
\cdot\cdot\cdot \overset{\psi_{i-1}}\to{\rightarrow} V_i \overset{\psi_i}\to{\rightarrow}
V_{i+1} \overset{\psi_{i+1}}\to{\rightarrow} \cdot\cdot\cdot
\overset{\psi_{l-2}}\to{\dashrightarrow} V_{l-1}
\overset{\psi_{l-1}}\to{\rightarrow} V_l = X_2$$ where

\ \ (i) $\phi = \psi_{l-1} \circ \psi_{l-2} \circ \cdot\cdot\cdot \circ \psi_2 \circ \psi_1$,

\ \ (ii) $\psi_i$ are isomorphisms over $U$, and

\ \ (iii) either $\psi_i:V_i \dashrightarrow V_{i+1}$ or $\psi_i^{-1}:V_{i+1} \dashrightarrow V_i$
is a morphism obtained by blowing up a smooth irreducible center disjoint from $U$.

Furthermore, if both $X_1$ and $X_2$ are projective, then we can choose a factorization so that all
the intermediate varieties $V_i$ are projective.  More precisely, in the factorization we
construct, there exists an index $i_o$ such that the birational map $V_i \dashrightarrow X_1$ is a
projective birational morphism for $i \leq i_o$ and that the birational map $V_i \dashrightarrow
X_2$ is a projective birational morphism for $i_o \leq i$.

\endproclaim

It is called the WEAK factorization theorem as we allow the sequence to consist of blowups and
blowdowns in any order.  If we insist on the sequence in the factorization to consist only of
blowups immediately followed by blowdowns, then we obtain the following STRONG factorization
conjecture, which remains open.

\proclaim{Strong Factorization Conjecture} Let $\phi:X_1 \dashrightarrow X_2$
be as above.  Then a sequence as above can be chosen with the extra condition that there
exists an index $i_o$ so that $\psi_i^{-1}$ for
$1 \leq i \leq i_o - 1$ are blowups of $V_i$ while $\psi_i$ for $i_o \leq i \leq l - 1$ are
blowdowns of $V_i$. 
\endproclaim

We note that, if both $X_1$ and $X_2$ are projective, in a sequence asserted in the strong
factorization conjecture all the intermediate varieties are automatically projective, while in an
arbitrary sequence in the weak factorization it may have some complete but nonprojective intermediate
varieties and that the ``Moreover" part of the weak factorization theorem asserts the
existence of a sequence preserving the projectivity. 

Several generalizations of the theorem will be discussed in Chapter 5.

\vskip.2in

\S 0-1. Outline of the strategy

\vskip.1in

First, we simply state the strategy of our proof in a symbolic way with a brief description of the
following five main steps:

\vskip.1in

$$\boxed{\text{Strategy\ of\ the\ Proof}}$$

\vskip.1in
 
Step 1: Reduction to the case where $\phi$ is a projective birational morphism

\vskip.1in

\hskip.4in $\bold{Elimination\ of\ points\ of\ indeterminacy}$

\newpage

Step 2: Factorization into locally toric birational maps

\vskip.1in

\hskip.4in $\bold{Construction\ of\ a\ birational\ cobordism}$

\vskip.1in

Step 3: Factorization into toroidal birational maps

\vskip.1in

\hskip.4in $\bold{Torification:Blowing\ up\ the\ torific\ ideal}$

\vskip.1in

Step 4: Recovery of Nonsingularity

\vskip.1in

\hskip.4in $\bold{Canonical\ resolution\ of\ singularities}$

\vskip.1in

Step 5. Factorization of toroidal birational maps among nonsingular toroidal embeddings

\vskip.1in

\hskip.4in $\bold{Morelli's\ combinatorial\ algorithm}$.

\vskip.2in

Now we will explain each step more in detail and describe the contents of the chapters of these
notes accordingly.

\vskip.1in

Let $\phi:X_1 \dashrightarrow X_2$ be a birational map between complete nonsingular varieties over
an algebraically closed field $K$ of characteristic zero.

\vskip.1in

Step 1: Reduction to the case where $\phi$ is a projective birational morphism

\vskip.1in

Hironaka's elimination of points of indeterminacy asserts that there exist sequences of
blowups with smooth centers $X_1' \rightarrow X_1$ and $X_2' \rightarrow X_2$ so that the induced map
$\phi':X_1'
\rightarrow X_2'$ is a projective birational morphism.  By replacing the original $\phi$ with
$\phi'$ we may assume in the factorization problem that $\phi$ is a projective birational
morphism.  This process is discussed in \S 1-4 of Chapter 1.

\vskip.1in

Step 2: Factorization into locally toric birational maps

\vskip.1in

This is the key step of our proof, discussed in Chapter 2, which allows us to bring in a ``Morse
Theoretic" view point to approach the factorization problem via the theory of ``Birational
Cobordism".

Let us briefly recall what the Morse theory tells us: We have a usual cobordism $C = C(M_1,M_2)$
between differentiable manifolds $M_1$ and $M_2$ equipped with a Morse function $f:C(M_1,M_2)
\rightarrow {\Bbb R}$.  (By extending the cobordism appropriately, we may assume that $f$ is
surjective.)

\newpage

\proclaim{Picture 0-1-1}\endproclaim

\vskip3in

By looking at the vector field $grad(f)$ (with respect to a Riemannian metric on $C(X_1,X_2)$), we
can introduce the action of the ``time" $t \in {\Bbb R}$ on the cobordism, namely, we have a map
from the time ${\Bbb R}$ to the group $Diff(C)$ of diffeomorphisms of $C$ denoted by
$$t \in {\Bbb R} \mapsto \varphi_t:C \rightarrow C \in Diff(C)$$
where $\varphi_t(p)$ is the position of a point $p \in C$ after time $t$ along the flow given by
integrating the vector field $grad(f)$.  Taking the exponential $exp(t) \in {\Bbb R}_{> 0}$ of the
time $t \in {\Bbb R}$, we see that the multiplicative group $({\Bbb R}_{> 0},\times) = exp({\Bbb
R},+)$ acts as well.  The crucial information in the Morse Theory can be read off in terms of
the action of the multiplicative group:

\ \ (C-o) the multiplicative group ${\Bbb R}_{> 0}$ acts on a cobordism $C$,

\ \ (C-i) $M_2$ and $M_1$ being on the top and at the bottom of the cobordism $C$, respectively, we
can describe them as the quotients
$$\align
M_2 &\cong C_+/{\Bbb R}_{> 0} \\
M_1 &\cong C_-/{\Bbb R}_{> 0} \\
\endalign$$
where the sets $C_+$ and $C_-$ are defined to be
$$\align
C_+ &:= \{x \in C;\lim_{t \rightarrow \infty}t(x) = \lim_{exp(t) \rightarrow \infty}exp(t)(x) \text{\
does\ NOT\ exist\ in\ }C\}
\\
C_- &:= \{x \in C;\lim_{t \rightarrow - \infty}t(x) = \lim_{exp(t) \rightarrow 0}exp(t)(x) \text{\
does\ NOT\ exist\ in\ }C\},
\\ 
\endalign$$

\ \ (C-ii) the critical points of $f$ are precisely the fixed points of the action and the homotopy
types of fibers of $f$ change as we pass through these fixed points.

\vskip.1in

This interpretation of the Morse Theory in terms of the action of the multiplicative group led
W{\l}odarczyk [2] naturally to the notion of a birational cobordism $B = B_{\phi}(X_1,X_2)$ for a
birational map
$\phi$:

\ \ (B-o) the multiplicative group $t \in K^*$ acts on a (nonsingular) variety $B$ of $\dim B = \dim
X_1 + 1 =
\dim X_2 + 1$,

\ \ (B-i) we can
describe $X_2$ and $X_1$ as the quotients
$$\align
X_2 &\cong B_+/K^* \\
X_1 &\cong B_-/K^* \\
\endalign$$
where the sets $B_+$ and $B_-$ are defined to be
$$\align
B_+ &:= \{x \in B;\lim_{t \rightarrow \infty}t(x) \text{\ does\ NOT\ exist\ in\ }B\}
\\
B_- &:= \{x \in B;\lim_{t \rightarrow 0}t(x) \text{\ does\ NOT\ exist\ in\ }B\}
\\ 
\endalign$$
and the birational map induced by the natural inclusions
$$X_1 = B_+/K^* \supset (B_+ \cap B_-)/K^* \subset B_-/K^* = X_2$$
coincides with $\phi$,

\ \ (B-ii) ``passing through" the fixed points of the action induces the birational transformations.

\vskip.1in

The precise definition of the birational cobordism is found in \S 2-1.  We give the construction of
a birational cobordism (only for a projective birational morphism, after the reduction Step 1, and
leaving the construction for a general birational map to the original W{\l}odarczyk [2]) in \S
2-2, together with the discussion of the ``collapsibility" in order to line up the fixed points
nicely in the above algebraic setting where we lack a Morse function
$f$, which, in the usual differentiable setting, would assign the levels to the critical points
and hence line them up nicely.  We also provide an interpretation of the notion of the birational
cobordism in
\S 2-3 from a view point of Geometric Invariant Theory (reviewed briefly in \S 1-1
together with Toric Geometry), in connection with the work of Thaddeus and others.  This view
point allows us to preserve the projectivity of the intermediate varieties in our factorization if
both
$X_1$ and $X_2$ are projective.

\vskip.1in

We would like to take a closer look at the birational transformations as we pass through the fixed
points.

\vskip.1in

Let $p \in B^{K^*}$ be a fixed point of the birational cobordism $B$.  Analytically locally, the
action of $K^*$ on $B$ at $p$ is equivalent to the action of $K^*$ on the tangent space $T_{B,p}$,
which is a vector space $V_p$ of $\dim V_p = n$.  As the multiplicative group $K^*$ is reductive,
the vector space $V_p$ splits into the eigenspaces according to the characters of the action.  By
choosing a basis $\{z_1, \cdot\cdot\cdot, z_n\}$ from the eigenspaces, we can describe the action
of $t \in K^*$ on $T_{B,p} = V_p$ as
$$t(z_1, \cdot\cdot\cdot, z_n) = (t^{\alpha_1} \cdot z_1, \cdot\cdot\cdot, t^{\alpha_n} \cdot
z_n)$$
where the $\alpha_i \in {\Bbb Z}$ are the characters for the $z_i$.  The crucial sets $(V_p)_+$ and
$(V_p)_-$ for the local behavior of the cobordism, defined similarly as $B_+$ and $B_-$, have simple
descriptions: 

$$\align
(V_p)_+ &= \{z = (z_1, \cdot\cdot\cdot, z_n) \in V_p; z_i \neq 0 \text{\ for\ some\ }i \text{\
with\ }\alpha_i > 0\}\\
(V_p)_- &= \{z = (z_1, \cdot\cdot\cdot, z_n) \in V_p; z_i \neq 0 \text{\ for\ some\ }i \text{\
with\ }\alpha_i < 0\}.\\
\endalign$$

If we identify the vector space $V_p$ with the affine $n$-space ${\Bbb A}^n = X(N,\sigma)$, regarded
as the toric variety associated to the standard cone $\sigma = \langle v_1, \cdot\cdot\cdot,
v_n \rangle \subset N_{\Bbb R} = N
\otimes {\Bbb R}$ generated by the standard ${\Bbb Z}$-basis of the lattice $N \cong {\Bbb Z}^n$, and
identify the action of
$K^*$ with the one-parameter subgroup corresponding to the point $a = (\alpha_1, \cdot\cdot\cdot,
\alpha_n) \in N$, we have a simple local description of the birational transformation, namely the
birational transformation between $(V_p)_+/K^*$ and $(V_p)_-/K^*$ in terms of toric geometry
associatde to the projection
$\pi:N_{\Bbb R} \rightarrow N_{\Bbb R}/{\Bbb R} \cdot a$;

$$\align
(V_p)_+/K^* &= X(\pi(N),\pi(\partial_+\sigma)) \\
(V_p)_-/K^* &= X(\pi(N),\pi(\partial_-\sigma)) \\
\endalign$$
where the upper boundary face and lower boundary face of $\sigma$ are defined to be
$$\align
\partial_+\sigma &:= \{v \in \sigma;v + \epsilon \cdot (-a) \not\in \sigma \text{\ for\ }\epsilon >
0\} \\
\partial_-\sigma &:= \{v \in \sigma;v + \epsilon \cdot a \not\in \sigma \text{\ for\ }\epsilon >
0\}. \\
\endalign$$
Therefore, analytically locally, the birational transformation induced by passing through the fixed
point $p \in B$ is equivalent to a toric birational transformation between toric varieties
corresponding to (the projections of) the upper boundary face and lower boundary face of the standard
cone
$\sigma$ (with respect to the one-parameter subgroup $a \in N$), which is a typical example of the
``polyhedral cobordism" used by Morelli to solve the factorization probelm for toric birational
maps.  Thus via the birational cobordism we succeed in factoring $\phi$ into locally toric birational
transformations as we pass through the fixed points (at the several levels ordered nicely by the
collapsibility).  This is the content of
\S 2-4.

\proclaim{Picture 0-1-2}\endproclaim

\newpage

One might naively expect at this point that since $\phi$ is now factored into ``locally toric"
birational maps and since ``toroidal" birational maps can be factored by Morelli's cmbinatorial
algorithm, we had already completed a proof for the weak factorization theorem, especially when the
notion of ``locally toric" and that of ``toroidal" are mixed and sometimes confused in the existing
literatures.  But there is a huge GAP between locally toric birational maps and toroidal birational
maps.  Actually the distinction between the two notions is one of the key issues in our argument
and the main bodies \S 1-2 and \S 1-3 of Chapter 1 are devoted to the detailed discussion and
clarification of them.

\vskip.1in

Step 3: Factorization into toroidal birational maps

\vskip.1in

The birational cobordism $B$ has a locally toric structure, as we observed above, each point
$p \in B$ having an analytic neghborhood isomorphic to a toric variety $V_p$, and the action of
$K^*$ is translated into that of a one-parameter subgroup in this locally toric chart.  But the
choice of coordinates (or equivalently the coordinate divisors) is not canonical and hence the
embedded tori vary from chart to chart as we move from point to point.  Thus these locally toric
coordinates do NOT patch together.  On the other hand, we require a toroidal structure $(U_X,X)$ to
have a fixed open subset
$U_X \subset X$ so that
$X - U_X$ provides the global coordinate divisors.  Thus in order to give $B$ (or its modification) a
toroidal structure we have to create such global divisors.  ``Torification" achieves exactly this
goal by blowing up the ``torific" ideal and hence creating these desired global coordinate divisors
defined by the pull-back of the torific ideal.  The torific ideal is canonically defined in terms of
the action of
$K^*$ only.

Here arises one technical but essential difficulty.  The torific ideal is only well-defined
over smaller pieces of the birational cobordism, called the quasi-elementary pieces of $B$, but not
on the entire birational cobordism $B$.  (The quasi-elementary pieces correspond to the semi-stable
loci of different linearizations of the action of $K^*$, when the birational cobordism is
interpreted via Geometric Invariant Theory.)   Thus we can only introduce toroidal structures,
after blowing up the torific ideals, to the quasi-elementary pieces and these toroidal structures may
not be compatible with each other.  As a consequence, though we succeed in factoring $\phi$ into
toroidal birational maps by torification, we have to pay the price that these toroidal birational
maps are only enough to provide weak factorization but not strong factorization, due to the
incompatibility of the toroidal structures.

The precise definition of the torific ideal, together with the subtlety why it is only well-defined
over the quasi-elementary pieces, is given in \S 3-1.  The core of ``torification", i.e., the
torifying property of the torific ideal asserting that blowing it up induces a toroidal structure
so that the action of $K^*$ is that of a one-parameter subgroup in a toroidal chart of each point,
is discussed in \S 3-2.

\vskip.1in

Step 4: Recovery from singular to nonsingular

\vskip.1in

Now that we succeeded in factoring $\phi$ into toroidal birational maps after Step 3, are we done
with the proof ?  Not yet !!

\vskip.1in

Though the birational maps are toroidal, the toroidal embeddings we have constructed through the
steps may be, and in most of the cases they actually are, SINGULAR.  (In fact, even at the stage of
factoring
$\phi$ into locally toric birational maps the intermediate varieties may be singular.)  Since we
seek to factor
$\phi$ into blowups and blowdowns with smooth centers on nonsingular varieties, we have to bring
back the whole situation from SINGULAR to NONSINGULAR.  This is done via the application of the
canonical resolution of singularities (and the canonical principalization of some ideals) in \S
4-1.  The canonicity of the resolution allows us to preserve the toroidal structure through the
process of desingularization.

\vskip.1in

Step 5. Factorization of toroidal birational maps among nonsingular toroidal embeddings

\vskip.1in

Finally, applying Morelli's combinatorial algorithm to toroidal birational maps between nonsingular
toroidal embeddings, we complete the proof of the weak factorization in \S 4-2.

\vskip.1in

In Chapter 5, we discuss several generalizations of the weak factorization theorem, modifying the
arguments discussed above according to the different purposes.  \S 5-1 accomplishes the
factorization of bimeromorphic maps in the analytic category.  In \S 5-2, we consider the case where
there is a group acting and establish the factorization equivariantly.  (We remark that we may have
to blow up several smooth centers simultaneously to maintain the equivariance.)  As an application
in \S 5-3, we show that we can remove the assumption of the base field $K$ being algebraically
closed.  It is noteworthy that we can NOT remove the assumption of the characteristic of $K$ being
zero for the time being, as the canonical resolution of singularities sits in the center of our
method, and that it is the only place where we use the characteristic assumption and hence any
future development of the method of the canonical resolution of singularities in positive
characteristics would allow us to remove the assumption.  \S 5-4 discusses the factorization in
the logarithmic category, which has an application to the study of the behavior of the Hodge
structures under birational transformations.  

\vskip.1in

In Chapter 6, we mention some problems related to our proof of the weak factorization conjecture. 
Though the proof, providing a specific procedure to factor a birational map, is constructive, it
falls short of being effective at several places, of which \S 6-1 make a list.  In \S 6-2, we discuss
briefly the toroidalization conjecture, which lies at the heart of the whole circle of ideas,
inspired by the work of De Jong [1] and others, involving the factorization problem, resolution of
singularities, semi-stable reductions, the log categories of Kato [1] and beyond.

\vskip.2in

\S 0-2. Some historical remarks

\vskip.1in

The following is the auhtor's personal and prejudiced view of how the factorization problem evolved
over years.  It is by NO means meant to be a complete account of the subject and the author
apologizes in advance for the omission of any references or topics which should be included.

\vskip.1in

\noindent $\bold{Early\ origins\ of\ the\ problem}$

\vskip.1in

A birational map $\phi:X_1 \dashrightarrow X_2$ between algebraic varieties is by definition a
rational map which induces an isomorphism of the function fields $\phi^*:K(X_2)
\overset{\sim}\to{\rightarrow} K(X_1)$ or equivalently a rational map which induces an isomorphism
over a common dense open subset $X_1 \supset U \subset X_2$.  

If we resrict our attention to complete
nonsingular varieties, a birational map is nothing but an isomorphism for curves, i.e., varieties
of dimension 1.  In dimension 2 or higher, however, the picture changes drastically.  There are many
examples of birational maps which are not isomorphisms, typical examples of which are BLOWUPS
with smooth centers, their inverses BLOWDOWNS and their composites.  A natural and fundamental
question arises then if these exhaust all the possible birational maps, i.e., if a given birational
map is a composite of blowups and blowdowns with smooth centers, which came to be known as the
factorization problem of birational maps.

The history of the factorization problem of birational maps could be traced back to the Italian
school of algebraic geometers, who already knew that the operation of blowing up points on surfaces
is a fundamental source of richness for surface geometry.  The importance of the strong
factorization theorem in dimension 2 by Zariski cannot be overestimated in the analysis of the
birational geometry of algebraic surfaces.  The strong factorization problem was stated in the
form of a question as ``Question $(F')$" in Hironaka [2], who recognized the connection with the
problem of resolution of singularities and gave a partial answer in the form of Elimination of
Points of Indeterminacy.  Though the question of the weak factorization was also raised in Oda
[1], the problem remained largely open in higher dimensions despite the efforts and interesting
results of many, e.g., Crauder [1] Kulikov [1] Schaps [1] and others.  These were summarized in
Pinkham [1], where the weak factorization conjecture is explicitly stated in the form
presented here.

\vskip.1in

\noindent $\bold{Toric\ Case}$

\vskip.1in

For toric birational maps, the equivariant version of the factorization problem under the torus
action was posed in Oda [1] and came to be known as Oda's weak and strong conjectures.  Though the
toric geometry can be interpreted in terms of the geometry of convex cones and hence
sometimes is considered to be a ``baby" or easier version of the real case, the factorization problem
presented a substantial challenge and difficulty in combinatorics.  (Actually the fact that all the
known examples (cf. Hironaka [1] Sally
[1] Shannon [1]) demonstrating the difficulties in higher dimensions were toric made some of us
even suspect that it is not a coincidence but an indication that all the combinatorial
complications in general are somehow concentrated in the toric case.)  In dimension 3, Danilov [2]'s
proof for the weak factorization was later supplemented by Ewald [1].  Oda's weak conjecture was
solved in arbitrary dimension by W{\l}odarczyk [1].  A big breakthrough was brought by Morelli
[1][2], who not only solved Oda's strong conjecture in arbitrary dimension but also introduced the
notion of a ``cobordism" between fans into the analysis.  This theory of the polyhedral cobordism
made the combinatorial approach to the problem conceptually transparent, so much so that it made us
eager to hope for the extension of the theory to be applied to general birational maps.  But because
of the apparent disguise as a purely combinatorial object with little hope of extension, we had to
wait for W{\l}odarczyk [2] to reveal the true meaning as the algebraic version of the Morse Theory.

\vskip.1in

\noindent $\bold{Minimal\ Model\ Program}$

\vskip.1in

It is worthwhile to note the relation of the factorization problem to the development of the Mori
program.  Hironaka [1] used the cone of effective curves to study the properties of birational
morphisms.  This direction was further developed and given a decisive impact by Mori [2], who,
motivated by Kleiman's criterion for ampleness and inspired by the solution to a
conjecture of Hartshorne (cf. Mori [1]), introduced the notion of extremal rays and systematically
used it in an attempt to construct minimal models in higher dimension, called the minimal model
program.  Danilov [2] introduced the notion of canonical and terminal singularities in conjunction
with the factorization problem.  This was developed by Reid [1][2][4] into a general theory of these
singularities, which appear in an essential way in the minimal model program.  The minimal model
program is so far proven up to dimension 3 (See Mori [3] Kawamata [1][2][3] Koll\'ar [1] Shokurov
[1].) and for toric varieties in arbitrary dimension (See Reid [3]).  In the steps of the minimal
model program one is only allowed to contract a divisor into a variety with terminal singularities
or to perform a flip, modifying some codimension 2 loci.  This allows a factorization of a given
projective birational morphism into such ``elementary operations".  An algorithm to factor
birational maps among uniruled varieties, known as ``Sarkisov's program", has also been developed
and established in dimension 3 in general (See Sarkisov [1] Reid [5] Corti [1].) and in arbitrary
dimension for toric varieties (See Matsuki [1]).  It was hoped that the detailed understanding of
such ``elementary operations" would yield a solution, at least in dimension 3, to the classical
factorization problem.  Still, we do not know of a way to carry out this approach.

\vskip.1in

\noindent $\bold{Resurgence\ of\ Toroidal\ Geometry}$

\vskip.1in

It was a shock when De Jong [1] showed that the resolution of singularities in positive
characteristics, though up to finite alterations, can be achieved and reduced to the toroidal case
by an elegant application of the theory of the moduli space of (marked) stable curves.  Resolution
of singularities (easy but weak) in characteristic zero, without finite alteration, was also achieved
shortly after by Abramovich-DeJong [1] and Bogomolov-Pantev [1] (See also Paranjape [1]), again
reducing the general case to the toroidal case.  When we saw that Abramovich-Karu [1] solved the
problem of semi-stable reduction over higher dimensional base as an application of the method above,
extending the original theorem by Kempf-Knudsen-Mumford-SaintDonat [1], the inspiration was
developing into a belief that toric (toroidal) geometry is not just a baby version testing the real
case but reflects the geometry in general and that the factorization problem can be approached in a
similar manner by some reduction step yet to be found.

\vskip.1in

\noindent $\bold{Local\ Version}$

\vskip.1in

There is a local version of the factorization problem, formulated and proved in dimension 2 by
Abhyankar [1].  Chritensen [1] posed the problem in general and solved it for some special cases in
dimension 3.  Here the varieties $X_1$ and $X_2$ for a birational map $\phi:X_1 \dashrightarrow X_2$
are replaced by appropriate birational local rings dominated by a fixed valuation, and blowups are
replaced by monoidal transforms subordinate to the valuation.  This
local conjecture was recently solved by S.D. Cutkosky in its strongest form in a series of papers
(cf. Cutkosky [1][2][3]).  The main ingredient of his proof is the ``monomialization of a given
birational map", i.e., the local version of the reduction step to the toroidal case.  Thus the
results of Cutkosky, though we do not use any of his method in our proof, not only gave us an
psycological support that there is no local obstruction to solving the global factorization
conjecture but also presented another evidence supporting the belief mentioned above.

\vskip.1in

\noindent $\bold{Birational\ Cobordism\ and\ Connection\ with\ G.I.T.}$

\vskip.1in

W{\l}odarczyk [2] provides the long-sought-after and key ingredient to reduce the factorization
problem of GENERAL birational maps to that of TOROIDAL birational maps by revealing the true nature
of Morelli's cobordism and hence opening the way to applying his combinatorial analysis to the
general case via the theory of ``birational cobordism", which can be regarded as an
algebraic version of the Morse Theory.  His theory can also be interpreted in the frame work of
Geometric Invariant Theory and has a direct connection with the subject of the change of G.I.T.
quotients associated to the change of linearizations, which has been intensively studied by
Thaddeus [1][2] and others.  

\vskip.1in

\noindent $\bold{Canonical\ Resolution\ of\ Singularities}$

\vskip.1in

There has been much progress after Hironaka's original work toward the understanding and
simplification of the algorithms of resolution of singularities, notably by Bierstone-Milman [1]
Villamayor [1] Encinas-Villamayor [1] among others.  Since the elimination of points of
indeterminacy, it has been a common consensus that the problem of resolution of singularities is in
close connection with the factorization problem.  In fact, the canonicity of
these algorithms for resolution of singularities plays an indispensible role in our proof for
the weak factorization.

\vskip.1in

\noindent $\bold{Solutions}$

\vskip.1in

The theory of birational cobordism by W{\l}odarczyk [2] almost immediately factors a given
birational map into locally toric birational transformations, the factorization which is a
consequence of a local analysis in nature.  But in order to apply Morelli's combinatorial
algorithm, which is a global in nature, we still have to bridge the gap between ``local" and
``global".

W{\l}odarczyk [3] provides this bridge by constructing a global combinatorial object, much like the
conical complexes associated to toroidal embeddings in Kempf-Knudsen-Mumford-SaintDonat [1], and
hence a solution for the weak factroization problem.  Abramovich-Karu-Matsuki-W{\l}odarczyk [1],
meanwhile, provides this bridge by transforming the locally toric structure into a (global)
toroidal structure by the process of torification.  It is obtained by blowing up the ``torific"
ideal.  The idea originates from the above-mentioned work of Abramovich-DeJong [1] for the weak
resolution of singularities.

Both proofs share essentially the same difficulty toward the strong factorization conjecture, which
remains as an open question.

\vskip.1in

\noindent $\bold{Today\ and\ Future}$

\vskip.1in

Though based upon the same belief, our proof presented here is different from what we had originally
conceived as a possible approach to the factorization problem.  King, in a joint paper with
Abkulut, essentially presented what is now known as the toroidalization conjecture,
our original approach, stated explicitly in Abramovich-Karu [1] Abramovich-Matsuki-Rashid [1].

\proclaim{Toroidalization Conjecture} Let $\phi:X_1 \dashrightarrow X_2$ be as above.  Then there
exist sequences of blowups and blowdowns with smooth centers $X_1' \rightarrow X_1$ and $X_2'
\rightarrow X_2$ such that the induced morphism $\phi':X_1' \rightarrow X_2'$ becomes toroidal (with
respect to suitably chosen toroidal structures on $X_1'$ and $X_2'$).
\endproclaim

The conjecture has been known to hold in diemsnion 2 by the original works of Abkulut-King [1],
Abramovich-Karu [1] Karu [1] and that of Abramovich presented in
Abramovich-Karu-Matsuki-W{\l}odarczyk [1], as well as by a recent local algorithmic proof by
Cutkosky-Piltant [1].  The toroidalization conjecture makes sense not only for a birational morphism
but also for a morphism $f:X \rightarrow Y$ between varieties of different dimensions, the
flexibility which makes us hopeful for the inductional structure in a possible proof.  (Cutkosky
announced a proof for the case
$\dim X = 3$ and $\dim Y = 2$ in the summer of 1999.)  

An affirmative solution for the toroidalization conjecture
will automatically provide one for the strong factorization conjecture.  

Also by resolution of
singularities, it is easily reduced to the following case of morphisms between toroidal
embeddings: Let $f:(U_X,X) \rightarrow (U_Y,Y)$ be a morphism between complete nonsingular toroidal
embeddings with $f:U_X = f^{-1}(U_Y) \rightarrow U_Y$ being smooth.  Then there should exist
sequences of blowups and blowdowns $\sigma_X:(U_{X'},X') \rightarrow
(U_X,X)$ and $\sigma_Y:(U_{Y'},Y') \rightarrow (U_Y,Y)$ with smooth centers sitting in the
boundaries so that the induced morphism
$f':(U_{X'},X') \rightarrow (U_{Y'},Y')$ is ``log-smooth", or equivalently, toroidal.  The notion of
a morphism being log-smooth was introduced by Kato [1] together with that of the log category.  (See
also the log category and logarithmic ramification formula of Iitaka [1].)  It gives us an intrinsic
justification of the toroidalization conjecture, which otherwise might stand as a mere technical
speculation.  Interpreted in this way, the toroidalization conjecture may be regarded as the
problem of resolution of singularities of morphisms (maps) in the log-category, similar to the
problem of resolution of singularities of varieties in the usual category.  The pursuit of a
canonical algorhithm for toroidalization is the future topic which invokes intensive research of
today.

\newpage

$$\bold{CHAPTER\ 1.\ PRELIMINARIES}$$

\vskip.2in

In the first part of this chapter, we recall some basic definitions and facts from Geometric
Invariant Theory (G.I.T. for short, cf. Kirwan-Fogarty-Mumford [1]) and Toric
Geometry (cf. Danilov [1] Fulton [1] Kempf-Knudsen-Mumford-SaintDonat [1] Oda [1]),
mainly to fix the notation for these lecture notes.  

The second part is devoted to the discussion of
``$\bold{locally\ toric}$" and ``$\bold{toroidal}$" structures.  Since one of the critical
issues in our argument is the distinction between the two notions and since in the existing
literatures they appear in some mixed terminologies, the understanding of this part is crucial
as well as basic.  

\vskip.1in

In the third part of this chapter, we discuss the actions of the multiplicative group $K^*$
of the base field on the locally and/or toroidal structures.  In the course of studying
birational transformations via birational cobordisms, we would like to introduce locally
toric and/or toroidal structures not only on the birational cobordisms but also on their
quotients by the actions of $K^*$.  This naturally leads us to the (well-known) notion of
a ``strongly \'etale" $K^*$-equivariant morphism, utilized in the form of $\bold{Luna's\
Fundamental\ Lemma}$ and $\bold{Luna's\ \acute{E}tale\ Slice\ Theorem}$.

\vskip.1in

In these lecture notes, all the varieties are assumed to be irreducible and reduced schemes of
finite type over an algebraically closed field
$K$ of characteristic zero, EXCEPT for Chapter 5 where we specifically try to loosen
some of the restrictions. 

\vskip.1in

\S 1-1. Brief Review of G.I.T. and Toric Geometry

\vskip.1in

\proclaim{Review 1-1-1 (Geometric Invariant Theory)}\endproclaim

Suppose a reductive group $G$
acts on an algebraic variety
$X$.  We denote by $X/G$ what we call the geometric quotient, i.e., the space of orbits, and
by $X//G$ the space of equivalence classes of orbits, where the equivalence 
relation is generated by the condition that two orbits are equivalent if their closures intersect.

These spaces are endowed with the structures (of algebraic varieties) which satisfy the usual
universal property (if they exist in the category of algebraic varieties).
 
\vskip.1in

In most of our situations, 

\ \ $\circ$ $G$ is taken to be the multiplicative group $K^*$ or a finite
group,
 
\ \ $\circ$ when we take
the quotient $X//G$, a variety $X$ is such that the closure of 

\ \ any orbit contains
a unique closed orbit, i.e., contains a unique fixed point (See the definition of a
quasi-elementary cobordism in Chapter 2),

\ \ $\circ$ if a variety $X$ is normal, then the quotient $X/G$ or $X//G$ is
also a normal variety (if it exists in the category of algebraic varieties) .  

\vskip.1in

If $X$ is endowed with an ample line bundle ${\Cal L}$ with a linear $G$-action compatible with
the
$G$-action on $X$ (called a linearization of the action), then we denote by $X_{\Cal L}^{ss}$ the
set of semi-stable points with respect to this linearization 
$$X_{\Cal L}^{ss} = \{x \in X; \exists s \in H^0(X,{\Cal L}^{\otimes n})^G \text{\ for\ some\
}n \in {\Bbb N}\text{\ s.t.\ }s(x) \neq 0\}$$ and by $X^{ss}_{\Cal L}//G$ the so-called
G.I.T. quotient, which has the structure of a projective algebraic variety
$$X^{ss}_{\Cal L}//G = \roman{Proj} \oplus_{n \geq 0}H^0(X,{\Cal L}^{\otimes n})^G.$$

\vskip.1in

We discuss more on G.I.T. related to the work of Thaddeus and others, i.e., in the subject of
the change of the G.I.T. quotients associated to the change of linearizations in Chapter 2.

\vskip.1in

\proclaim{Review 1-1-2 (Toric Geometry)}\endproclaim

Let $N \cong {\Bbb Z}^n$ be a lattice and $\sigma \subset N_{\Bbb R} = N \otimes {\Bbb R}$ be
a strictly convex rational polyhedral cone.  We denote the dual lattice by $M = Hom_{\Bbb
Z}(N,{\Bbb Z})$ and the dual cone by $\sigma^{\vee} \subset M_{\Bbb R} = M \otimes {\Bbb R}$. 
The affine toric variety $X$ associated to $\sigma$ is defined to be
$$X = X(N,\sigma) = \roman{Spec}\ K[M \cap \sigma^{\vee}],$$
containing the torus
$$X \supset T_X = \roman{Spec}\ K[M].$$
The lattice $N$ is identified with the group of one-parameter subgroups of $T_X$ and the dual
lattice with their characters.  As we only consider toric varieties associated to the saturated 
lattices $M \cap \sigma^{\vee}$ as above, they are all normal by definition. 

More generally the toric variety associated to a fan $\Sigma$ in $N_{\Bbb R}$ is denoted by
$X(N,\Sigma)$.

If $X_1 = X(N,\Sigma_1)$ and $X_2 = X(N,\Sigma_2)$ are toric varieties having the same lattice
$N$, the embeddings of the torus $T_{X_1} = T_{X_2} = \roman{Spec}\ K[M]$ define an equivariant
birational map $\phi:X_1
\dashrightarrow X_2$, called a toric birational map.  This map $\phi$ is a morphism if and
only if every cone in $\Sigma_1$ is contained in a cone in $\Sigma_2$ and proper if and only
if the support of $\Sigma_1$ coincides with that of $\Sigma_2$.

\vskip.1in

Suppose that $K^*$ acts on an affine toric variety $X = X(N,\sigma)$ as a
one-parameter subgroup of the torus $T_X$, corresponding to a lattice point $a \in N$.  If $t
\in K^*$ and $m \in M$, the action on the monomial $z^m$ is given by
$$t^*(z^m) = t^{(a,m)} \cdot z^m,$$
where $(\cdot,\cdot)$ is the natural pairing on $N \times M$.  (For $m \in M$, we denote its
image in the semi-group algebra $K[M \cap \sigma^{\vee}]$ by $z^m$.)

The $K^*$-invariant monomials correspond to the lattice points $M \cap a^{\perp}$, hence
$$X//K^* = \roman{Spec}\ K[M \cap \sigma^{\vee} \cap a^{\perp}].$$
Let $\pi:N_{\Bbb R} \rightarrow N_{\Bbb R}/a \cdot {\Bbb R}$ be the projection.  Then the
lattice $M \cap a^{\perp}$ is dual to the lattice $\pi(N)$ and 
$$M \cap
\sigma^{\vee} \cap a^{\perp} = Hom_{\Bbb Z}(\pi(N),{\Bbb Z}) \cap \pi(\sigma)^{\vee}.$$
Therefore, $X//K^* = X(\pi(N),\pi(\sigma))$ is again an affine toric variety defined by the
lattice $\pi(N)$ and the cone $\pi(\sigma)$ (though $\pi(\sigma)$ may no longer be strictly
convex).  This quotient is a geometric quotient precisely when $\pi:\sigma \rightarrow
\pi(\sigma)$ is a linear isomorphism (while the map between the lattices $\pi:N \cap \sigma
\rightarrow \pi(N) \cap \pi(\sigma)$ may not be an isomorphism).

\vskip.1in

\S 1-2. ``Locally Toric" Structures vs. ``Toroidal" Structures 

\vskip.1in

There is some confusion in the existing literatures among the terminologies expressing the
notion of toroidal embeddings and toroidal maps (cf. Kempf-Knudsen-Mumford-SaintDonat [1]
Abramovich-Karu [2]) and that of toroidal varieties (cf. Danilov [1] W{\l}odarczyk [2]), the
latter of which we prefer to call locally toric varieties (and locally toric maps).  Since
one of the critical issues in our argument is the distinction of the two notions ``locally
toric" and ``toroidal", we clarify our usage of these two terminologies below.

\proclaim{Definition 1-2-1 (Locally Toric $\&$ Toroidal Structures)}

(i) A variety $W$ is locally toric if for every closed point $p \in W$ there exists a closed
point $x \in X$ in an affine toric variety such that we have an isomorphism of complete local
rings over $K$
$$\widehat{{\Cal O}_{W,p}} \hskip.1in \overset{\sim}\to{\rightarrow} \hskip.1in \widehat{{\Cal
O}_{X,x}}.$$

(ii) A pair $(U_W,W)$ consisting of a dense open subset $U_W \subset W$ contained in a
variety $W$ is called a toroidal embedding (or alternatively we say $(U_W,W)$ has a toroidal
structure) if for every closed point $p \in W$ there exists a closed point $x \in X$ in an
affine toric variety such that we have an isomorphism of complete local rings over $K$
$$\widehat{{\Cal O}_{W,p}} \hskip.1in \overset{\overset{\eta^*}\to{\sim}}\to{\rightarrow}
\hskip.1in
\widehat{{\Cal O}_{X,x}}$$ 
which induces an isomorphism of ideals defining the boundaries
$$\widehat{{\Cal I}_{W - U_W}} = {\Cal I}_{W - U_W} \otimes \widehat{{\Cal O}_{W,p}} \hskip.1in \overset{\overset{\eta^*}\to{\sim}}\to{\rightarrow}
\hskip.1in \widehat{{\Cal I}_{X - T_X}} = {\Cal I}_{X - T_X} \otimes \widehat{{\Cal
O}_{X,x}}$$
where ${\Cal I}_{W - U_W}$ (resp. ${\Cal I}_{X - U_X}$) is the ideal defining the reduced
structure on the boundary $W - U_W$ (resp. $X - T_X$).

If all the irreducible components of the boundary $W - U_W$ is normal, then $(U_W,W)$ is
called a toroidal embedding $\bold{without\ self-intersection}$.

\endproclaim

\proclaim{Definition 1-2-2 (Locally Toric $\&$ Toroidal Birational Morphisms)} 

\ \ (i) A proper birational morphism $f:W_1 \rightarrow W_2$ between locally toric varieties is
locally toric if for every closed point $p_1 \in W_1$ mapping to $p_2 = f(p_1) \in W_2$, there
exists a toric morphism $\varphi:X_1 \rightarrow X_2$ mapping a closed point
$x_1
\in X_1$ to $x_2 = \varphi(p_1) \in X_2$ such that we have a commutative diagram of
homomorphisms between complete local rings
$$\CD
\widehat{{\Cal O}_{W_1,p_1}} \hskip.1in @.\overset{\overset{\eta_1^*}\to{\sim}}\to{\rightarrow}
@.\hskip.1in \widehat{{\Cal O}_{X_1,x_1}} \\
@A{f^*}AA @. @AA{\varphi^*}A \\
\widehat{{\Cal O}_{W_2,p_2}} \hskip.1in @.
\overset{\overset{\eta_2^*}\to{\sim}}\to{\rightarrow} @.
\hskip.1in \widehat{{\Cal O}_{X_2,x_2}}. \\
\endCD$$

\ \ (ii) A proper birational morphism $f:(U_{W_1},W_1) \rightarrow (U_{W_2},W_2)$ between
toroidal embeddings is toroidal if for every closed point $p_1 \in W_1$ mapping to $p_2 =
f(p_1) \in W_2$, there exists a toric morphism $\varphi:X_1 \rightarrow X_2$ mapping a
closed point
$x_1
\in X_1$ to $x_2 = \varphi(p_1) \in X_2$ such that we have a commutative diagram of
homomorphisms between complete local rings
$$\CD
\widehat{{\Cal O}_{W_1,p_1}} \hskip.1in @.
\overset{\overset{\eta_1^*}\to{\sim}}\to{\rightarrow} @.\hskip.1in 
\widehat{{\Cal O}_{X_1,x_1}} \\
@A{f^*}AA @.@AA{\varphi^*}A \\
\widehat{{\Cal O}_{W_2,p_2}} \hskip.1in @.
\overset{\overset{\eta_2^*}\to{\sim}}\to{\rightarrow} @.
\hskip.1in \widehat{{\Cal O}_{X_2,x_2}} \\
\endCD$$
which induces a commutative diagram of ideals defining the boundaries
$$\CD
\widehat{({\Cal I}_{W_1 - U_{W_1}})_{p_1}} \hskip.1in @.
\overset{\overset{\eta_1^*}\to{\sim}}\to{\rightarrow} @.\hskip.1in 
\widehat{({\Cal I}_{X_1 - U_{X_1}})_{x_1}} \\
@A{f^*}AA @.@AA{\varphi^*}A \\
\widehat{({\Cal I}_{W_1 - U_{W_1}})_{p_2}} \hskip.1in @.
\overset{\overset{\eta_2^*}\to{\sim}}\to{\rightarrow} @.\hskip.1in 
\widehat{({\Cal I}_{X_1 - U_{X_1}})_{x_2}}. \\
\endCD$$
\endproclaim

Recall that a birational map $f:W_1 \dashrightarrow W_2$ is proper (cf. Iitaka [1]) if both
projections $p_1:\Gamma_f \rightarrow W_1$ and $p_2:\Gamma_f \rightarrow W_2$ are proper
morphisms in the usual sense, where $\Gamma_f \subset W_1 \times W_2$ is the graph of $f$.  (If
$f$ is well-defined as a morphism over a dense open subset $U \subset W_1$, then $\Gamma_f$ is the
closure in $W_1
\times W_2$ of the graph of the morphism $f|_U:U \rightarrow W_2$.  This closure is
independent of the choice of the open subset $U$.)

\proclaim{Definition 1-2-3 (Locally Toric $\&$ Toroidal Birational Maps)}

\ \ (i) A proper birational map $f:W_1 \dashrightarrow W_2$ between locally toric
varieties is locally toric if there exists another locally toric variety $Z$ which dominates
both $W_1$ and $W_2$ by proper birational morphisms
$$W_1 \overset{f_1}\to{\leftarrow} Z \overset{f_2}\to{\rightarrow} W_2$$
where we require $f_1$ and $f_2$ to be locally toric as in Definition 1-2-2 (i).

\ \ (ii) A proper birational map $f:(U_{W_1},W_1) \dashrightarrow (U_{W_2},W_2)$ between
toroidal embeddings is toroidal if there exists another toroidal embedding $(U_Z,Z)$
which dominates both $(U_{W_1},W_1)$ and $(U_{W_2},W_2)$ by proper birational morphisms
$$(U_{W_1},W_1) \overset{f_1}\to{\leftarrow} (U_Z,Z) \overset{f_2}\to{\rightarrow}
(U_{W_2},W_2)$$ 
where we require $f_1$ and $f_2$ to be toroidal as in Definition 1-2-2
(ii).
\endproclaim

\proclaim{Remark 1-2-4}\endproclaim

\ \ (i) Though we presented the notions of ``locally toric" and ``toroidal" as if they were
parallel and of equal importance, the collection of locally toric varieties with locally toric
birational morphisms and/or birational maps does NOT form a category, in clear contrast to the
fact that the collection of toroidal embeddings with toroidal birational morphisms and/or
toroidal birational maps does form a category as shown in Proposition 1-2-5. 

A composition of locally toric birational morphisms or maps is NOT locally toric in
general (cf. Proposition 1-2-5 (i) (ii)): We take $g = f_2 \circ f_1$ to be a composite of the
blowup
$f_2:W_2
\rightarrow W_3$ of a point on a smooth threefold $W_3$ and the blowup $f_1:W_1 \rightarrow
W_2$ of a smooth curve on
$W_2$ which is not transversal to the exceptional divisor for $f_2$.  Then while $f_1$ and
$f_2$ are locally toric, the composition $g$ is not.

It is NOT clear and probably not true (cf. Proposition 1-2-5 (iii)) that the definition of
being locally toric for a birational morphism in Definition 1-2-2 (i) is compatible with the
definition of being locally toric for a birational map in Definition 1-2-3 (i), though the
author does not have a specific example showing incompatibility.

\ \ (ii) There are examples where a locally toric birational morphism $f:W_1 \rightarrow W_2$
cannot be toroidal no matter how we introduce toroidal structures on $W_1$ and $W_2$: Let
$W_2$ be a smooth surface.  First we blow up a point on $W_2$ and then blow up 3 or more
distinct points on the exceptional divisor to obtain $f:W_1 \rightarrow W_2$. 

\vskip.1in

\proclaim{Proposition 1-2-5}

\ \ (i) A composition of two (proper and) toroidal birational morphisms is again (proper and)
toroidal.

\ \ (ii) A composition of two (proper and) toroidal birational maps is again (proper and)
toroidal.

\ \ (iii) A proper birational morphism between toroidal embeddings is toroidal in the sense
of Definition 1-2-3 (ii) if and only if it is toroidal in the sense of Definition 1-2-2
(ii).

In particular, the collection of toroidal embeddings with toroidal birational morphisms and/or
toroidal birational maps forms a category.
\endproclaim

\demo{Proof}\enddemo The assertion (i) is proved in Abramovich-Karu [1] Karu [1] for a
composition of general (not necessarily birational) toroidal morphisms.

The assertions (ii) and (iii) can be proved as an application of the following lemmas, after
reducing the case of general toroidal embeddings to that of toroidal embeddings without
self-intersection via some appropriate toroidal blowups and base changes.

\proclaim{Lemma 1-2-6} Let $(U_W,W)$ be a toroidal embedding.  Then there exists a sequence
$(U_{\tilde W},{\tilde W}) \rightarrow (U_W,W)$ of toroidal blowups with smooth centers
(i.e., the centers analytically locally correspond to the smooth closures of orbits) such
that the resulting toroidal embedding $(U_{\tilde W},{\tilde W})$ is nonsingular and without
self-intersection.
\endproclaim

\demo{Proof}\enddemo First take the canonical resolution of singularities $r:W^{res}
\rightarrow W$ and set $U_{W^{res}} = r^{-1}(U_W)$.  Then by the property
$(\spadesuit^{res}-1)$ of the canonical resolution (We refer the reader to Remark 4-1-1 for
the details of the properties of the canonical resolution of singularities.) we conclude that
$r:(U_{W^{res}},W^{res}) \rightarrow (U_W,W)$ is a toroidal morphism (between toroidal
embeddings) obtained as a sequence of toroidal blowups with smooth centers.  Now starting
from blowing up the 0-dimensional strata of the boundary components of $(U_{W^{res}},W^{res})$, we
consecutively blow up (the closures of) the strict transforms of the original $k$-dimensional
strata of the boundary of the original toroidal structure $(U_{W^{res}},W^{res})$ for
$k = 1, \cdot\cdot\cdot, \dim W - 1$ to obtain $s:(U_{\tilde W} = s^{-1}(U_{W^{res}}),{\tilde W})
\rightarrow (U_{W^{res}},W^{res})$.  (Note that the closures of the stata of the boundary
components are the closures of the generic points which, analytically locally, correspond to
the intersections of the boundary divisors of the toroidal charts.)  This is a sequence of
toroidal blowups with smooth centers so that the resulting toroidal embedding
$(U_{\tilde W},{\tilde W})$ is nonsingular and without self-intersection.

\proclaim{Lemma 1-2-7} Let $f:(U_{W_1},W_1) \rightarrow (U_{W_2},W_2)$ be a proper and
toroidal birational morphism between toroidal embeddings without self-intersection.  Then $f$
is allowable (in the sense of Kempf-Knudsen-Mumford-SaintDonat [1]).  In particular, such
proper and toroidal birational morphisms over $(U_{W_2},W_2)$ are in one-to-one correspondence
with the subdivisions of the conical complex associated to the toroidal embedding
$(U_{W_2},W_2)$ without self-intersection.
\endproclaim

\demo{Proof}\enddemo For a proof, we refer the reader to Abramovich-Karu [1] Karu [1] and
Kempf-Knudsen-Mumford-SaintDonat [1].

\vskip.1in

In order to see the assertion (ii), suppose
$$f:(U_{W_1},W_1) \dashrightarrow (U_{W_2},W_2) \text{\ and\ } g:(U_{W_2},W_2) \dashrightarrow
(U_{W_3},W_3)$$
are proper and toroidal birational maps among toroidal embeddings.  Let
$$(U_{W_1},W_1) \overset{f_1}\to{\leftarrow} (U_Z,Z) \overset{f_2}\to{\rightarrow}
(U_{W_2},W_2) \text{\ and\ } (U_{W_2},W_2) \overset{g_2}\to{\leftarrow} (U_{Z'},{Z'})
\overset{g_3}\to{\rightarrow} (U_{W_3},W_3)$$
be some toroidal embeddings dominating $(U_{W_1},W_1), (U_{W_2},W_2)$ and $(U_{W_3},W_3)$ by
proper and toroidal birational morphisms as required in Definition 1-2-3 (ii).  By taking a
sequence of toroidal blowups with smooth centers of $(U_{W_2},W_2)$ using Lemma 1-2-6, we can
make it without self-intersection.  Since $f_2$ and $g_2$ are toroidal, we can ``pull-back" these
toroidal blowups to the toroidal blowups on $(U_Z,Z)$ and on $(U_{Z'},Z')$.  Again using Lemma
1-2-6 and the assertion (i), we may assume that $(U_Z,Z), (U_{W_2},W_2)$ and $(U_{Z'},Z')$ are
all without self-intersection.  By Lemma 1-2-7, the morphisms $f_2$ and $g_3$ are associated with
the subdivisions $\Delta_Z$ and $\Delta_{Z'}$ of the conical complex $\Delta_{W_2}$.  Let
$(U_{\hat Z},{\hat Z}) \rightarrow (U_{W_2},W_2)$ be the toroidal embedding associated to
the common refinement $\Delta_{\hat Z}$ of $\Delta_Z$ and $\Delta_{Z'}$.  Then by
construction $(U_{\hat Z},{\hat Z})$ dominates $(U_Z,Z)$ and $(U_{Z'},Z')$ by proper and
toroidal birational morphisms
$$(U_{W_1},W_1) \overset{f_1}\to{\leftarrow} (U_Z,Z) \overset{{\hat f}_1}\to{\leftarrow}
(U_{\hat Z},{\hat Z}) \overset{{\hat g}_3}\to{\rightarrow} (U_{Z'},Z')
\overset{g_3}\to{\rightarrow} (U_{W_3},W_3).$$
Since the compositions $f_1 \circ {\hat f}_1$ and $g_3 \circ {\hat g}_3$ are both toroidal by
Lemma 1-2-6, we conclude that the composition $g \circ f$ is also a toroidal birational map by
Definition 1-2-3 (ii).

\vskip.1in

In order to see the assertion (iii), let $f:(U_{W_1},W_1) \rightarrow (U_{W_2},W_2)$ be
a proper birational morphism between toroidal embeddings.  If $f$ is toroidal as a birational
morphism according to Definition 1-2-2 (ii), then it is obviously toroidal as a birational
map according to Definition 1-2-3 (ii) by taking the diagram
$$(U_{W_1},W_1) \overset{identity}\to{\leftarrow} (U_{W_1},W_1) \overset{f}\to{\rightarrow}
(U_{W_2},W_2).$$
Thus we only have to show the converse, i.e., if $f$ is toroidal as a birational map as in
Definition 1-2-3 (ii), then it is toroidal as a birational morphism as in Definition 1-2-2 (ii). 

According to Definition 1-2-3 (ii), there exists another toroidal embedding $(U_Z,Z)$
which dominates both $(U_{W_1},W_1)$ and $(U_{W_2},W_2)$ by proper birational morphisms
$$(U_{W_1},W_1) \overset{f_1}\to{\leftarrow} (U_Z,Z) \overset{f_2}\to{\rightarrow}
(U_{W_2},W_2)$$ 
where $f_1$ and $f_2$ are proper and toroidal birational morphisms according to Definition
1-2-2 (ii).  By Lemma 1-2-6 and Proposition 1-2-5 (i), we may assume that $(U_Z,Z)$ is a
nonsingular toroidal embedding without self-intersection.  We denote by $E_k$ the irreducible
components of the boundary $Z - U_Z$. 

Let $p_1 \in W_1$ be a closed point mapping to $p_2 = f(x_1) \in W_2$.

Via the analytic isomorphisms as in the definition of the toroidal structures
$$\align
\widehat{{\Cal O}_{{W_1},p_1}} &\overset{\sim}\to{\rightarrow} \widehat{{\Cal O}_{{X_1},x_1}}
\\
\widehat{{\Cal O}_{{W_2},p_2}} &\overset{\sim}\to{\rightarrow} \widehat{{\Cal O}_{{X_2},x_2}}
\\
\endalign$$
we pull back the standard coordinates for the tori
$$\align
\{z_1, \cdot\cdot\cdot, z_n\} &\text{\ for\ the\ torus\ }T_{X_1} = \roman{Spec}\ K[z_1^{\pm
1}, \cdot\cdot\cdot, z_n^{\pm 1}] \\
\{w_1, \cdot\cdot\cdot, w_n\} &\text{\ for\ the\ torus\ }T_{X_2} = \roman{Spec}\ K[w_1^{\pm
1}, \cdot\cdot\cdot, w_n^{\pm 1}], \\
\endalign$$
which we denote by the same letters by abuse of notation.

We only have to show the claim that there exists a unimodular $n \times n$ matrix $A =
\{a_{ij}\}$ and units $u_{i} \in \widehat{{\Cal O}_{{W_1},p_1}}$ for $i, j = 1,
\cdot\cdot\cdot, n$ such that
$$f^*w_i = u_i \prod z_j^{a_{ij}} \text{\ in\ the\ field\ of\ fractions\ of\ }\widehat{{\Cal
O}_{{W_1},p_1}}.$$

Remark that
$$\align
\widehat{f_1}:& (U_Z,Z) \times \roman{Spec}\ \widehat{{\Cal O}_{{W_1},p_1}} \rightarrow
(U_{W_1},W_1) \times \roman{Spec}\ \widehat{{\Cal O}_{{W_1},p_1}} \\
\widehat{f_2}:& (U_Z,Z) \times \roman{Spec}\ \widehat{{\Cal O}_{{W_2},p_2}} \rightarrow
(U_{W_2},W_2) \times \roman{Spec}\ \widehat{{\Cal O}_{{W_2},p_2}} \\ 
\endalign$$
are both proper and birational toroidal morphisms between toroidal embeddings WITHOUT
self-intersection and hence allowable and that they correspond to the subdivisions of the conical
complexes by Lemma 1-2-7, which we denote by
$$\align
\Delta_{\widehat{f_1}}:&\Delta_{Z,p_1} \rightarrow \Delta_{W_1,p_1} \\
\Delta_{\widehat{f_2}}:&\Delta_{Z,p_2} \rightarrow \Delta_{W_2,p_2}.\\
\endalign$$ 

Remark also that $\Delta_{Z,p_1}$ can be considered as a subcomplex of $\Delta_{Z,p_2}$,
since $f_1^{-1}(p_1) \subset f_2^{-1}(p_2)$.

Suppose $\dim \Delta_{Z,p_1} = \dim \Delta_{Z,p_2}$ (while in general we have $\dim
\Delta_{Z,p_1} \leq \dim \Delta_{Z,p_2}$).  Let $d_1$ be the minimum among the dimensions of
the strata whose closures are the connected components of the intersections of the
$E_k$ and which have nonempty intersection with $f_1^{-1}(p_1)$.  Let $d_2$ be the minimum
among the dimensions of the strata whose closures are the connected components of the
intersections of the $E_k$ and which have nonempty intersection with $f_2^{-1}(p_2)$.  Then the
assumption is equivalent to $d_1 = d_2$ (while in general we have $d_1 \geq d_2$).  In this case,
$\Delta_{Z,p_1}$ is a conical complex of full dimension embedded in another conical complex
$\Delta_{Z,p_2}$ with the compatible integral structures.  More precisely, since $f_2$ factors
as $f_2 = f_1 \circ f$, we can take a ${\Bbb Z}$-basis of the integral structure for
$\Delta_{W_1,p_1}$ (which hence provides a common ${\Bbb Z}$-basis for each of the cones of the
maximal dimension in $\Delta_{Z,p_1}$) which is compatible with a ${\Bbb Z}$-basis of
the integral structure of $\Delta_{W_2,p_2}$ (which hence provides a common ${\Bbb Z}$-basis for each of the
cones of the maximal dimension in $\Delta_{Z,p_2}$).  It follows that there exists a
unimodular matrix $A = \{a_{ij}\}$ such that for all $z \in f_1^{-1}(p_1) \subset Z$ we have
$$f_2^*w_i = (u_i)_z \prod z_j^{a_{ij}} \text{\ in\ the\ field\ of\ fractions\ of\
}\widehat{{\Cal O}_{Z,z}}$$
for some units $(u_i)_z \in \widehat{{\Cal O}_{Z,z}}$ for $i, j = 1, \cdot\cdot\cdot, n$.
Since $f_1:Z \rightarrow W_1$ is proper and $z \in f_1^{-1}(p_1)$ is arbitrary, we have the
claim.

Suppose $\dim \Delta_{Z,p_1} < \dim \Delta_{Z,p_2}$, i.e., $d_1 > d_2$.  Let $S$ be a
stratum of dimension $\dim S = d_1$ with nonempty intersection with $f_1^{-1}(p_1)$.  Then
there exists another stratum $S' \subset S$ of dimension $\dim S' = d_2$ with nonempty
intersection with $f_2^{-1}(p_2)$.  Take a point $z' \in S' \cap f_2^{-1}(p_2)$ and let
$p_1' = f_1(z')$.  By the argument for the previous case, for the pull-backs via the
isomorphism as in the definition of the toroidal structure
$$\widehat{{\Cal O}_{{W_1},p_1'}} \overset{\sim}\to{\rightarrow} \widehat{{\Cal
O}_{{X_1'},x_1'}}$$
of the standard coordinates
$$\{z_1', \cdot\cdot\cdot, z_n'\} \text{\ for\ the\ torus\ }T_{X_1'} =
\roman{Spec}\ K[{z_1'}^{\pm 1}, \cdot\cdot\cdot, {z_n'}^{\pm 1}]$$
there exists a unimodular matrix $A = \{a_{ij}\}$ such that
$$f_2^*w_i = (u_i)' \prod (z_j')^{a_{ij}}$$
for some units $(u_i)' \in \widehat{{\Cal O}_{{W_1},p_1'}}$.

Note that every divisor $E_k$ with $E_k \cap f_1^{-1}(p_1) \neq \emptyset$ we have $E_k \cap
f_1^{-1}(p_1')$ by the minimality of the dimension $d_1 = \dim S$ of the stratum $S$ and by the
allowability of the morphism 
$$\widehat{f_1}: (U_Z,Z) \times \roman{Spec}\ \widehat{{\Cal O}_{{W_1},p_1}} \rightarrow
(U_{W_1},W_1) \times \roman{Spec}\ \widehat{{\Cal O}_{{W_1},p_1}}.$$
These also imply that there exists an isomorphism as in the definition of the toroidal structure
$$\widehat{{\Cal O}_{{W_1},p_1}} \overset{\sim}\to{\rightarrow} \widehat{{\Cal
O}_{{X_1},x_1}}$$
which pulls back the standard coordinates
$$\{z_1, \cdot\cdot\cdot, z_n\} \text{\ for\ the\ torus\ }T_{X_1} = \roman{Spec}\ K[z_1^{\pm
1}, \cdot\cdot\cdot, z_n^{\pm 1}]$$
with the property
$$\roman{ord}_{E_k}(z_j) = \roman{ord}_{E_k}(z_j')$$
for all $E_k$ with $E_k \cap f_1^{-1}(p_1)$ and $j = 1, \cdot\cdot\cdot, n$.  (Roughly put, we
can find a point $p_1''$ which sits in the same stratum as $p_1$ with respect to the natural
stratification of the boundary $W_1 - U_{W_1}$ and hence the toroidal structures at $p_1$ and
$p_1''$ are isomorphic, and which is close to $p_1'$ and hence admits the same locally toric
charts as $p_1'$.)  Therefore, we
have again for an arbitrary point $z \in f_1^{-1}(p_1)$
$$f_2^*w_i = (u_i)_z \prod z_j^{a_{ij}} \text{\ in\ the\ field\ of\ fractions\ of\
}\widehat{{\Cal O}_{Z,z}}$$
for some units $(u_i)_z \in \widehat{{\Cal O}_{Z,z}}$ for $i, j = 1, \cdot\cdot\cdot, n$.
Since $f_1:Z \rightarrow W_1$ is proper, we have the
claim.

This completes the proof of Proposition 1-2-5.

\vskip.1in

\proclaim{Remark 1-2-8}\endproclaim

\ \ (i) As we have seen in Remark 1-2-4 (i), a composition of two locally toric birational
morphisms is not necessarily locally toric, since the locally toric structures for the two
locally toric morphisms may not be compatible in the middle.  Even worse, the locally toric
structures for the morphisms
$W_1
\overset{f_1}\to{\leftarrow} Z$ and $Z \overset{f_2}\to{\rightarrow} W_2$ in Definition 1-2-3
(i) of a locally toric birational map may not be compatible on $Z$.  Thus it is nontrivial to
prove the assertions (i), (ii) and (iii) in Proposition 1-2-5, settleing a rather subtle issue
of the compatibility of toroidal structures for several different toroidal morphisms.

\vskip.1in

\ \ (ii) We define a stronger version of locally toric
and toroidal birational maps, namely V-locally toric and V-toroidal birational maps
below.  These turn out to be the only maps we have to consider in our note, the fact which
enables us to avoid the subtlety as above.  

\vskip.1in

Therefore, except for Lemma 1-2-6 and Lemma 1-2-7, we will NOT use the results of Proposition 1-2-5
in the later argument of our proof.

\vskip.1in

Actually the only place where Lemma 1-2-7 appears in our note is when we use the allowablity in
the proof of Theorem 4-2-1.  In the case of $V$-toroidal birational maps the allowablity (after
reducing the factorization to that of a birationa map between toroidal embeddings without
self-intersection) in Theorem 4-2-1 is automatic and hence our argument does not use Lemma
1-2-7, though it appears in Theorem 4-2-1 dealing with the factorization of general toroidal
birational maps.  

\vskip.1in

\proclaim{Remark 1-2-9 (Comparison with the definitions in Abramovich-Karu-Matsuki-W{\l}odarczyk
[1])}\endproclaim

In Abramovich-Karu-Matsuki-W{\l}odarczyk [1], they give more
restrictive definitions for locally toric and toroidal structures and locally toric and
toroidal birational morphisms/maps via the use of (allowable) Zariski locally toric and toroidal
charts.  

\vskip.1in

\noindent Definition of locally toric and toroidal structures in
Abramovich-Karu-Matsuki-W{\l}odarczyk [1]: 

\ \ (i) A variety $W$ is locally toric if for every closed
point $p \in W$ there exists a Zariski open neighborhood $V_p \subset W$ of $p$ and an \'etale
morphism $\eta_p:V_p \rightarrow X_p$ to a toric variety $X_p$.

\ \ (ii) A pair $(U_W,W)$ consisting of a dense open subset $U_W \subset W$ in a variety $W$ is
called a toroidal embedding if for every closed point $p \in W$ there exists a Zariski open
neighborhood $V_p$ and an
\'etale morphism
$\eta_p:V_p \rightarrow X_p$ to a toric variety such that $\eta_p^{-1}(T_{X_p}) = U_W$ where
$T_{X_p}$ is the torus for $X_p$. 

\vskip.1in

\noindent Definition of locally toric and toroidal morphisms in
Abramovich-Karu-Matsuki-W{\l}odarczyk [1]:

\ \ (i) A proper birational morphism $f:W_1 \rightarrow W_2$ between locally toric varieties is
called locally toric if for every closed point $p_2 \in W_2$ there exists a Zariski open
neighborhood $V_{p_2}$ and an
\'etale morphism $\eta_{p_2}:V_{p_2} \rightarrow X_{p_2}$ to a toric variety, together with a
toric proper birational morphism $Y_{p_2} \rightarrow X_{p_2}$ such that we have a fiber square
$$\CD
f^{-1}(V_{p_2}) = V_{p_2} \times_{X_{p_2}}Y_{p_2} \hskip.1in @. \rightarrow @. \hskip.2in Y_{p_2}
\\ @VVV \square @. \hskip.2in @VVV \\
V_{p_2} \hskip.1in @. \rightarrow @. \hskip.2in X_{p_2} \\
\endCD$$ 

\ \ (ii) A proper birational morphism $f:(U_{W_1},W_1) \rightarrow (U_{W_2},W_2)$ between
toroidal embeddings is called toroidal if for every closed point $p_2 \in W_2$ there exists a Zariski open
neighborhood $V_{p_2}$ and an
\'etale morphism $\eta_{p_2}:V_{p_2} \rightarrow X_{p_2}$ to a toric variety, together with a
toric proper birational morphism $Y_{p_2} \rightarrow X_{p_2}$ such that we have a fiber square
of toroidal embeddings
$$\CD
(U_{W_1} \cap f^{-1}(V_{p_2}),f^{-1}(V_{p_2})) = V_{p_2} \times_{X_{p_2}}(T_{Y_{p_2}},Y_{p_2})
\hskip.1in @.
\rightarrow @.
\hskip.2in (T_{Y_{p_2}},Y_{p_2}) \\ 
@VVV \square @. \hskip.2in @VVV \\
(U_{W_2} \cap V_{p_2},V_{p_2}) = V_{p_2} \times_{X_{p_2}}(T_{X_{p_2}},X_{p_2}) \hskip.1in @.
\rightarrow @.
\hskip.2in (T_{X_{p_2}},X_{p_2}) \\
\endCD$$ 

\noindent Definition of locally toric and toroidal birational maps in
Abramovich-Karu-Matsuki-W{\l}odarczyk [1]:

\ \ (i) A proper birational map $f:W_1 \dashrightarrow W_2$ between locally toric
varieties is locally toric if there exists another locally toric variety $Z$ which dominates
both $W_1$ and $W_2$ by proper birational morphisms
$$W_1 \overset{f_1}\to{\leftarrow} Z \overset{f_2}\to{\rightarrow} W_2$$
where we require $f_1$ and $f_2$ to be locally toric as in (i) above.

\ \ (ii) A proper birational map $f:(U_{W_1},W_1) \dashrightarrow (U_{W_2},W_2)$ between
toroidal embeddings is toroidal if there exists another toroidal embedding $(U_Z,Z)$
which dominates both $(U_{W_1},W_1)$ and $(U_{W_2},W_2)$ by proper birational morphisms
$$(U_{W_1},W_1) \overset{f_1}\to{\leftarrow} (U_Z,Z) \overset{f_2}\to{\rightarrow}
(U_{W_2},W_2)$$ 
where we require $f_1$ and $f_2$ to be toroidal as in (ii) above.

\vskip.2in

The above definitions differ from ours in several subtle ways.  Here are a couple of
differences:

\ \ $\circ$ The toroidal embeddings in their sense are automatically without self-intersection,
while those in our sense may not.  For example, a pair $(U_S,S)$ consisting of a nonsingular
surface
$S$ with an irreducible nodal curve as a boundary component $S - U_S$ is a toroidal embedding in
our sense but not in their sense.

\ \ $\circ$ Locally toric birational morphisms in our sense may not be locally toric in their
sense.  For example, as in Remark 1-2-4, if we blowup a point on a smooth surface $W_2$ and
further blowup three or more points on the exceptional divisor to obtain $W_1$, then the
morphism $f:W_1 \rightarrow W_2$ is locally toric in our sense but not in their sense. 

\vskip.1in

On the other hand, we would like to emphasize
that for the category of toroidal embeddings WITHOUT self-intersection and proper birational
morphisms and maps, theur definitions and our definitions coincide.

\vskip.1in

In fact, we have the following:

\ \ $\circ$ The toroidal embedding without self-intersection in our sense always has a Zariski
toroidal chart as required in their definition.

Let $p \in (U_W,W)$ be a closed point in a toroidal embedding without self-intersection in our
sense and let $x \in X$ be a closed point in an affine toric variety as described in Definition
1-2-1 (ii).  The Cartier divisors in $\roman{Spec}\ \widehat{{\Cal O}_{W,p}} =
\roman{Spec}\ \widehat{{\Cal O}_{X,x}}$ supported on the support of $\widehat{{\Cal
O}_{W,p}}/\widehat{I_{W - U_W}} = \widehat{{\Cal O}_{X,x}}/\widehat{I_{X - T_X}}$ give rise to
a lattice in a cone of effective Cartier divisors (the dual cone sits in the conical complex
associated to
$(U_W,W)$, which is well-defined as
$(U_W,W)$ is without self-intersection).  We only provide a proof in case the dual cone has the
full dimension (being equal to $n = \dim W$), leaving the general case as an exercise
to the reader.  Let $z_1, \cdot\cdot\cdot, z_n$ be the standard coordinates for the torus $T_X
= \roman{Spec}\ K[z_1^{\pm}, \cdot\cdot\cdot, z_n^{\pm}]$.  Since $z_1, \cdot\cdot\cdot, z_n$
define Cartier divisors which form a ${\Bbb Z}$-basis of the lattice and since $(U_W,W)$ is
without self-intersection, we can find rational functions $r_1, \cdot\cdot\cdot, r_n \in K(W)$
for $W$ such that the Cartier divisors associated to them coincide, analytically locally at $p$,
with those associated to
$z_1, \cdot\cdot\cdot, z_n$.  Now let $m_1, \cdot\cdot\cdot, m_l$ be the monomials in $z_1,
\cdot\cdot\cdot, z_n$ generating the affine coordinate ring for $X$ and hence the maximal ideal
corresponding to the point
$x
\in X$, and let $m_1', \cdot\cdot\cdot, m_l'$ be the same monomials but in $r_1, \cdot\cdot\cdot,
r_l$.  In a suitable affine Zariski open neighborhood $V_p$, we may assume all the monomials
$m_1',
\cdot\cdot\cdot, m_l'$ are regular over $V_p$.  It is easy to see $r_1, \cdot\cdot\cdot, r_n$
are transcendental over $K$, as well as $z_1, \cdot\cdot\cdot, z_n$ are over $K$.  Therefore,
we have a well-defined morphism
$$\eta_p:V_p \rightarrow X_p = X$$
defined by sending $m_1, \cdot\cdot\cdot, m_l$ to $m_1', \cdot\cdot\cdot, m_l'$,
respectively.  Now by construction we see that the homomorphism
$${\eta_p}^*:\widehat{{\Cal O}_{X,x}} \rightarrow \widehat{{\Cal O}_{V,p}}$$
is surjective and that it is also injective as so is ${\eta_p}^*:A(X) \rightarrow A(V_p)$,
hence an isomorphism.  (Note that we are in the situation where we DO have the preservation of
injectivity at the completion level, unlike the one in Gabrielov [1][2]'s counter-example to a
conjecture of Grothendieck [1].  See, e.g., H\"ubl [1].)  Therefore, by shrinking $V_p$, we may
assume
$\eta_p$ is
\'etale.  By shrinking
$V_p$ further, we may also assume that all the irerducible divisors in $V_p \cap (W - U_W)$
contain $p$ and hence $\eta_p^{-1}(T_X) = U_W \cap V_p$.   Thus $\eta_p:V_p \rightarrow X_p$ is
a desired Zariski toroidal chart. 

\vskip.1in

\ \ $\circ$ The toroidal birational morphisms in our sense are toroidal in their sense.

Let $f:(U_{W_1},W_1) \rightarrow (U_{W_2},W_2)$ be a proper birational morphism which is toridal
in our sense according to Definition 1-2-2 (ii).  For a point $p_2 \in W_2$ we take a Zariski
toroidal chart $\eta_p:V_{p_2} \rightarrow X_{p_2}$ constructed as above.  By shrinking $X_{p_2}$
we may assume all the orbits in $X_{p_2}$ have nonempty intersection with the image
$\eta_p(V_{p_2})$.  Then it is easy to see that the Cartier divisors in $V_{p_2}$ supported on
$V_{p_2}
\cap (W - U_W)$ are in one-to-one correspondence with cartier divisors in $X_{p_2}$ supported on
$X_{p_2} - T_{X_{p_2}}$ and that the conical complexes $\Delta_{V_{p_2}}$ associated to $V_{p_2}
\cap (W - U_W)$ and
$\Delta_{X_{p_2}}$ associated to $X_{p_2}$ coincide.  Now by the results of Abramovich-Karu [1]
Karu [1],
$f$ is allowable and hence corresponds to a subdivision of $\Delta_{V_{p_2}}$ over $V_{p_2}$. 
By taking $Y_{p_2} \rightarrow X_{p_2}$ to be the toric birational morphism corresponding to the
same subdivision of $\Delta_{X_{p_2}}$, we obtain the commutative diagram of fiber squares
$$\CD
(U_{W_1} \cap f^{-1}(V_{p_2}),f^{-1}(V_{p_2})) = V_{p_2} \times_{X_{p_2}}(T_{Y_{p_2}},Y_{p_2})
\hskip.1in @.
\rightarrow @.
\hskip.2in (T_{Y_{p_2}},Y_{p_2}) \\ 
@VVV \square @. \hskip.2in @VVV \\
(U_{W_2} \cap V_{p_2},V_{p_2}) = V_{p_2} \times_{X_{p_2}}(T_{X_{p_2}},X_{p_2}) \hskip.1in @.
\rightarrow @.
\hskip.2in (T_{X_{p_2}},X_{p_2}). \\
\endCD$$ 
Therefore, $f$ is also toridal in their sense.

\vskip.1in

Now that we see that the their definitions and our definitons coincide for the category
of toroidal embeddings without self-intersection, the statements of Proposition 1-2-5 (i), (ii)
and (iii) also hold using their definitions in the category
of toroidal embeddings without self-intersection.

\vskip.1in

\proclaim{Definition 1-2-10}

\ \ (i) A proper birational map $f:W_1 \dashrightarrow W_2$ between locally toric varieties
is called V-locally toric if there exists another locally toric variety $Y$ which both $W_1$
and $W_2$ dominate by proper birational morphisms
$$W_1 \rightarrow Y \leftarrow W_2$$
such that for every closed point $y \in Y$ we can find \'etale morphisms $i:V \rightarrow
Y$ with $y \in i(V)$ and $\eta:V \rightarrow X_Y$ to an affine toric variety $X_Y$ which some
toric varieties $X_1$ and $X_2$ dominate by proper toric birational morphisms
$$X_1 \rightarrow X_Y \leftarrow X_2$$
with the property that when we take the fiber-products with $V$ they coincide, i.e., we have
a commutative diagram of the form 

$$\CD
W_1 @. \rightarrow @. \hskip.25in Y @. \hskip.25in \leftarrow @. W_2 \\
@AAA \square @. \hskip.25in @AAA \hskip.25in \square @. @AAA \\
W_1 \times_Y V @. \rightarrow @. \hskip.25in V @. \hskip.25in \leftarrow @. W_2 \times_Y V \\
@| @. \hskip.25in @| \hskip.25in @.@| \\
X_1 \times_{X_Y} V @. \rightarrow @. \hskip.25in  V @. \hskip.25in \leftarrow @. X_2 \times_{X_Y} V \\
@VVV \square @. \hskip.25in @VVV \hskip.25in \square @. @VVV \\
X_1 @. \rightarrow @. \hskip.25in X_Y @. \hskip.25in \leftarrow @. X_2.\\
\endCD$$

(Note that $X_Y$ and hence $X_1\ \&\ X_2$ are allowed to vary depending on the points $y \in Y$.)

\ \ (ii) A proper birational map $f:(U_{W_1},W_1) \dashrightarrow (U_{W_2},W_2)$ between
toroidal embeddings is called V-toroidal if there exists another toroidal embedding $(U_Y,Y)$
which both
$(U_{W_1},W_1)$ and
$(U_{W_2},W_2)$ dominate by proper birational morphisms
$$(U_{W_1},W_1) \rightarrow (U_Y,Y) \leftarrow (U_{W_2},W_2)$$
such that for every closed point $y \in Y$ we can find \'etale morphisms $i:V \rightarrow
Y$ with $y \in i(V)$ and $\eta:V \rightarrow X_Y$ to an affine toric variety $X_Y$ which some
toric varieties $X_1$ and $X_2$ dominate by proper toric birational morphisms
$$X_1 \rightarrow X_Y \leftarrow X_2$$
with the property that when we take the fiber-products with $V$ they coincide, i.e., we have
a commutative diagram of the form 

$$\CD
(U_{W_1},W_1) @. \rightarrow @. \hskip.25in (U_Y,Y) @. \hskip.25in \leftarrow @. (U_{W_2},W_2) \\
@AAA \square @.\hskip.25in @AAA \hskip.25in \square @. @AAA \\
(U_{W_1},W_1) \times_Y V @. \rightarrow @. \hskip.25in (U_Y,Y) \times_Y V @. \hskip.25in \leftarrow @.
(U_{W_2},W_2)
\times_Y V
\\ @| @. \hskip.25in @|@. \hskip.25in  @| \\
(T_{X_1},X_1) \times_{X_Y} V @. \rightarrow @. \hskip.25in (T_{X_Y},X_Y) \times_Y V @. \hskip.25in
\leftarrow @. (U_{X_2},X_2)
\times_{X_Y} V
\\ @VVV \square @. \hskip.25in @VVV \hskip.25in \square @.@VVV \\
(T_{X_1},X_1) @. \rightarrow @. \hskip.25in (T_{X_Y},X_Y) @. \hskip.25in \leftarrow @. (T_{X_2},X_2).\\
\endCD$$

(Note again that $X_Y$ and hence $X_1\ \&\ X_2$ are allowed to vary depending on the points $y \in
Y$.)
\endproclaim

\vskip.2in

\S 1-3. $K^*$-action on locally toric and toroidal structures

\vskip.1in

\proclaim{Definition 1-3-1 (Strongly \'Etale)} Let $V$ and $X$ be affine varieties with
$K^*$-actions and let
$$\eta:V \rightarrow X$$
be a $K^*$-equivariant \'etale morphism.  Then $\eta$ is said to be $\bold{strongly\ {\acute
e}tale}$ if

\ \ (i) the quotient map
$$\eta//K^*:V//K^* \rightarrow X//K^*$$
is \'etale, and

\ \ (ii) the natural map
$$V \rightarrow X \times_{X//K^*}V//K^*$$
is an isomorphism.

\endproclaim 

\proclaim{Remark 1-3-2}\endproclaim

\ \ (i) The statement below does NOT hold in general: 

\noindent For a $K^*$-equivariant \'etale morphism $\eta:V \rightarrow X$ between
affine varieties with $K^*$-actions, the quotient map $\eta//K^*:V//K^* \rightarrow X//K^*$ is
also \'etale.

For example, take 
$$V = \roman{Spec}\ K[u,u^{-1},y] \text{\ with\ the\ $K^*$-action\ given\ by\ }t \cdot (u,
y) = (tu, t^{-1}y)$$ 
and 
$$X = \roman{Spec}\ K[w,w^{-1},y] \text{\ with\ the\ $K^*$-action\ given\ by\ }t \cdot
(w,y) = (t^2w,t^{-1}y),$$ 
while a $K^*$-equivariant \'etale morphism $\eta:V \rightarrow X$ is
associated to the ring homomorphism $\eta^*:K[w,y] \rightarrow K[u,y]$ defined by $\eta^*(w) =
u^2$.  Then the quotient map is given by
$$\eta//K^*:V//K^* = \roman{Spec}\ K[uy] \rightarrow X//K^* = \roman{Spec}\ K[wy^2]$$
which ramifies over the origin as $wy^2 = (uy)^2$.

The ``reason" why $\eta//K^*$ fails to be \'etale is that in general the stabilizer $Stab(v)
\subset K^*$ of a point $v \in V$ is ``smaller" than the stabilizer $Stab(x) \subset K^*$ of the
image
$x =
\eta(v)$.  (The inclusion $Stab(v) \subset Stab(x) \subset K^*$ is obvious.)  Thus, roughly
speaking, the quotient $X//K^*$ is divided more than $V//K^*$ by the difference
$Stab(x)/Stab(v)$.  In the above example, if a point $p \in V$ has the coordinate $y = 0$ the
stabilizer is $Stab(v) = \{1\}
\subset K^*$, while the image $x = \eta(v)$ has the stabilizer $Stab(x) = \{\pm 1\} \subset K^*$.

It is not difficult to come up with an example where $V//K^*$ is nonsingular while $X//K^*$ is
singular, and hence $\eta//K^*:V//K^* \rightarrow X//K^*$ obviously fails to be \'etale.

\vskip.1in

\ \ (ii) On the other hand, if $\eta:V \rightarrow X$ is strongly \'etale between affine
varieties with $K^*$-actions, it follows from the condition (ii) that
$$Stab(v) = Stab(x) \hskip.1in\forall v \in V \text{\ where\ }x = \eta(v).$$

\ \ (iii) Luna's Fundamental Lemma stated below can be considered as the statement in (i) with
the extra condition that $Stab(v) = Stab(x)$ for all $v \in V$ where $x = \eta(v)$ is the image,
when properly interpreted.  Therefore, in a very rough sense, a $K^*$-equivariant \'etale
morphism $\eta:V \rightarrow X$ is strongly \'etale if and only if the stabilizers are preserved
by $\eta$.

\vskip.1in

\proclaim{Definition 1-3-3 (Locally Toric and Toroidal $K^*$-actions)}

\ \ (i) Let $W$ be a locally toric variety with a $K^*$-action.  We say that the action is
locally toric if for every closed point $p \in W$ we can find a $K^*$-invariant affine
neighborhood
$p
\in U_p \subset W$, an affine variety $V_p$ with a $K^*$-action and an affine toric variety
$X_p$ with $K^*$ acting as a one-parameter subgroup such that we have $K^*$-equivariant and
strongly \'etale morphisms $\eta_p$ and $i_p$
$$X_p \overset{\eta_p}\to{\leftarrow} V_p \overset{i_p}\to{\rightarrow} U_p.$$
We call such a set of stronly \'etale morphisms $\bold{Luna's\ locally\ toric\ chart\ at\ p}$
for the locally toric $K^*$-action. 

\ \ (ii) Let $(U_W,W)$ be a toroidal embedding with a $K^*$-action.  We say that the action is
toroidal if for every closed point $p \in W$ we can find a $K^*$-invariant affine
neighborhood
$p \in U_p \subset W$, an affine toroidal embedding $(U_{V_p},V_p)$ with a $K^*$-action and an
affine toric variety
$(T_{X_p},X_p)$ with $K^*$ acting as a one-parameter subgroup such that we have $K^*$-equivariant
and strongly \'etale morphisms $\eta_p$ and $i_p$
$$(T_{X_p},X_p) \overset{\eta_p}\to{\leftarrow} (U_{V_p},V_p) \overset{i_p}\to{\rightarrow}
(U_W \cap U_p,U_p),$$
where in the condition (ii) of $\eta_p$ and $i_p$ being strongly \'etale as stated in
Definition 1-3-1 we require that the natural morphisms
$$\align
(U_{V_p},V_p) &\rightarrow (T_{X_p},X_p)\times_{X_p//K^*}V_p//K^* \\
(U_{V_p},V_p) &\rightarrow (U_W \cap U_p,U_p)\times_{U_p//K^*}V_p//K^* \\
\endalign$$
to be the isomorphisms of the toroidal embeddings and hence $\eta_p$ and $i_p$ coincide with the
obvious \'etale toroidal morphisms
$$\CD
(T_{X_p},X_p) @<{\eta_p}<< (U_{V_p},V_p) \\ 
@| @| \\
(T_{X_p},X_p)\times_{X_p//K^*}X_p//K^*
@<<< (T_{X_p},X_p)\times_{X_p//K^*}V_p//K^*\\
(U_{V_p},V_p) @>{i_p}>> (U_W \cap U_p,U_p) \\ 
@| @| \\
(U_W \cap U_p,U_p)\times_{U_p//K^*}V_p//K^* @>>> (U_W \cap
U_p,U_p)\times_{U_p//K^*}U_p//K^*.\\
\endCD$$
We call such a set of strongly \'etale (toroidal) morphisms $\bold{Luna's\ toroidal\ chart\ at\
p}$ for the toroidal $K^*$-action.

\endproclaim

\proclaim{Proposition 1-3-4} Any $K^*$-action on a nonsingular variety $W$ is locally toric,
i.e., for every closed point $p \in W$ we can find Luna's locally toric chart at $p$.
\endproclaim

\proclaim{Remark 1-3-5}\endproclaim

We give two alternative proofs for Proposition 1-3-4:

\ \ A) One which uses the canonical resolution of
singularities and Luna's Fundamental Lemma.  In this proof, we show that Luna's locally toric chart
$$X_p \overset{\eta_p}\to{\leftarrow} V_p \overset{i_p}\to{\rightarrow} U_p$$
can be taken so that $i_p$ is an isomorphism, i.e., $V_p = U_p$ is a $K^*$-invariant Zariski
open neighborhood of $p$ and that $X_p$ is nonsingular.  Moreover, we show that in case $p \in W$ is
not a fixed point Luna's locally toric chart can be chosen so that none of $U_p = V_p$ or $X_p$ have
any fixed points. 

\ \ B) One which uses Luna's \'Etale Slice Theorem.  In this proof, we show that Luna's locally toric chart can be taken so that $i_p$ is
surjective in general and $X_p$ is nonsingular and that $i_p$ is an isomorphism if $p \in W$ is a
fixed point of the
$K^*$-action.  Moreover, we show that in case $p \in W$ is not a fixed point Luna's
locally toric chart can be chosen so that none of $U_p$, $V_p$ or $X_p$ has any fixed points. 

While Luna's Fundamental Lemma is characteristic free, the canonical resolution of singularities
(for the moment) and Luna's \'Etale Slice Theorem are only valid in characteristic 0.

\vskip.1in

\demo{Proof}\enddemo A) By Sumihiro's equivariant completion theorem (cf. Sumihiro [1][2]), we
can embed $W$ $K^*$-equivariantly into a complete variety $\overline{W}$ with a $K^*$-action. 
By taking the canonical resolution of singularities, which is a sequence of blowups with centers
outside of $W$ (cf. the condition $(\spadesuit^{res}-0)$ of the canonical resolution in Remark
4-1-1), we may assume $\overline{W}$ is nonsingular as well.  Note that the $K^*$-action lifts to
the canonical resolution by the condition $(\spadesuit^{res}-1)$.

Let $p \in W$ be a closed point.  Since $\overline{W}$ is complete, there exists a fixed point
$q \in \overline{O(p)}$ in the closure of the orbit of $p$.  By another theorem of Sumihiro,
there exists a $K^*$-invariant affine open neighborhood $q \in V_q = \roman{Spec}\ A(V_q)$.  Since
the maximal ideal $m_q \subset A(V_q)$ associated to the fixed point $q$ is $K^*$-invariant, it
splits into a direct sum of eigenspaces 
$$m_q = \oplus_{\alpha \in {\Bbb Z}}m_{q,\alpha}$$
where $m_{q,\alpha}$ consists of all the eigenfunctions (in $m_q$) of $K^*$-character $\alpha$. 
Therefore, we can choose eigenfunctions 
$$f_1, \cdot\cdot\cdot, f_n \in m_q\ (n = \dim W)$$
with $K^*$-characters $\alpha_1, \cdot\cdot\cdot, \alpha_n$, i.e.,
$$t^*(f_j) = t^{\alpha_j} \cdot f_j \text{\ for\ }t \in K^*,$$
such that they form a basis of
$m_q/m_q^2$ and hence that they generate the maximal ideal over ${\Cal O}_{W,q}$.  Consider the
morphism
$$\eta:V_q = \roman{Spec}\ A(V_q) \rightarrow X = {\Bbb A}^n = \roman{Spec}\ [z_1, \cdot\cdot\cdot,
z_n]$$ defined by
$$\eta^*(z_1) = f_1, \cdot\cdot\cdot, \eta^*(z_n) = f_n.$$
Letting $t \in K^*$ act on ${\Bbb A}^n$ by $t^*(z_j) = t^{\alpha_j} \cdot z_j$, i.e., letting
$K^*$ act as the one parameter subgroup $a = (\alpha_1, \cdot\cdot\cdot, \alpha_n) \in N$ on
${\Bbb A}^n$ regarded as a standard toric variety, and shrinking
$V_q$ if necessary, we see that $\eta$ is a $K^*$-equivariant \'etale morphism.

\proclaim{Lemma 1-3-6 (Luna's Fundamental Lemma)} Let $G$ be a reductive group acting on affine
varieties $V$ and $X$.  Let $\eta:V \rightarrow X$ be a $G$-equivariant morphism.  Let $T$ be a
closed orbit of $G$ in $V$ such that

\ \ (i) $\eta$ is \'etale at some point of $T$,

\ \ (ii) $\eta(T)$ is closed in $X$,

\ \ (iii) $\eta$ is injective on $T$, and

\ \ (iv) $V$ is normal along $T$.

Then there are $G$-stable open subsets $V' \subset V$ and $X' \subset X$, with $T \subset V'$, such
that $\eta|_{V'}:V' \rightarrow X'$ is a strongly \'etale $G$-equivariant morphism of $V'$ onto
$X'$.
\endproclaim

First we apply Luna's Fundamental Lemma to $\eta:V \rightarrow X$ with $T = \{q\}$ and $G =
K^*$ to observe that for the choice of $K^*$-invariant open subsets $V' \subset V$ and $X' \subset
X$ as above
$$V' \rightarrow X' \times_{X'//K^*}V'//K^*$$
is an isomorphism.  

Since $O(p) \subset V'$, wthis implies that $\eta$ is injective on $O(p)$. 
The morphism $\eta$ is \'etale at all points of $O(p)$.  Let $J = \{j; f_j(p) \neq 0\} \subset
\{1, \cdot\cdot\cdot, n\}$.   By replacing $V$ and $X$ with
$$\align
&\{v \in V; \prod_{j \in J}f_j(v) \neq 0\} \\
&\{x \in X; \prod_{j \in J}z_j(x) \neq 0\}, \\
\endalign$$
respectively, we may assume $\eta(O(p))$ is closed in $X$, which is still an affine
nonsingular toric variety.  By assumption $V
\subset W$ is nonsingular.  With these conditions satisfied, we apply Luna's Fundamental Lemma
secondly to
$\eta:V
\rightarrow X$ with $T = O(p)$ and $G = K^*$ to have
$$\eta|_{V'}:V' \rightarrow X' (\subset X)$$
being strongly \'etale and hence
$$V' \rightarrow X' \times_{X'//K^*}V'//K^*$$
being an isomorphism.

\vskip.1in

Case: $p \in W$ is not a fixed point. 

\vskip.1in

In this case, by the first application of Luna's Fundamental
Lemma we conclude that there exists $j \in J$ with $f_j(p) \neq 0$ and $\alpha_j \neq 0$.  By the
choice of $X$, we see that $X$ has no fixed points.  Thus we see that $X'//K^* \subset
X//K^*$ is an open subset.  

Therefore, we conclude

\ \ (i) $\eta_p = \eta|_{V'}:V_p = V' \rightarrow X_p = X$ is \'etale, and

\ \ (ii) $V_p = V' \rightarrow X' \times_{X'//K^*}V'//K^* = X_p \times_{X_p//K^*}V_p//K^*$ is an
isomorphism.

Thus
$$X_p \overset{\eta_p}\to{\leftarrow} V_p \overset{i_p}\to{=} U_p$$
is Luna's locally toric chart at $p$.

\vskip.1in

Case: $p$ is a fixed point.

\vskip.1in

In this case, $\eta(p)$ is also a fixed point.  Let $F_X$ be the set of fixed points in $X$ with
respect to the $K^*$-action.  Let
$$\pi_X:X \rightarrow X//K^*$$
be the quotient map.  Since $F_X - X'$ is a $K^*$-invariant closed subset of $X$, we conclude
that $\pi_X(F_X - X')$ is a closed subset of $X//K^*$ and hence that $\pi_X^{-1}(\pi_X(F_X - X'))$
is a closed subset of $X$.  Therefore, we conclude that the subsets defined below
$$\align
C &:= F_X \cap X' \\
D &:= \pi_X^{-1}(\pi_X(F_X - X')) \cap X'\\
\endalign$$
are disjoint $K^*$-invariant closed subsets of $X'$.  Therefore, by Corollary 1.2 in
Mumford-Fogarty-Kirwan [1], there exists a $K^*$-invariant function $f \in A(X')^{K^*}$ such that
$$f \equiv 1 \text{\ on\ }C \hskip.1in \& \hskip.1in f \equiv 0 \text{\ on\ }D.$$
Then by construction, setting
$$X'_f = \{x \in X';f(x) \neq 0\},$$
we see that $\eta(p) \in X'_f$ and that $X'_f//K^* \subset X//K^*$ is an open subset.  Set $V'_f =
\eta^{-1}(X'_f)$.  Then since
$$V' \rightarrow X' \times_{X'//K^*}V'//K^*$$
is an isomorphism, so is
$$V'_f \rightarrow X'_f \times_{X'_f//K^*}V'_f//K^*.$$  
Therefore, we
conclude

\ \ (i) $\eta_p = \eta|_{V_p}:V_p = V'_f \rightarrow X_p = X$ is \'etale, and

\ \ (ii) $V_p \rightarrow X'_f \times_{X'_f//K^*}V'_f//K^* = X_p \times_{X_p//K^*}V_p//K^*$ is an
isomorphism.

Thus
$$X_p \overset{\eta_p}\to{\leftarrow} V_p \overset{i_p}\to{=} U_p$$
is Luna's locally toric chart at $p$.

This completes the proof A) for Proposition 1-3-4 using the canonical resolution of
singularities and Luna's Fundamental Lemma.

\vskip.2in

B) We recall Luna's \'Etale Slice Theorem.

\proclaim{Theorem 1-3-7 (Luna's \'Etale Slice Theorem)} Let $G$ be a reductive group.  Let $U$ be
an affine variety with a $G$-action and let $T$ be a closed orbit of $G$ in $U$ along which $U$
is normal.  If $p \in T$, then there exists a locally closed $G_p$-stable ($G_p$ is the
stabilizer of $p$ in $G$) affine subvariety $Z$ of $U$ with $p \in Z$ such that $U_p = G \cdot
Z$ is affine open in $U$ and the $G$-equivariant morphism
$$\iota:G \times_{G_p}Z \rightarrow U_p$$
is strongly \'etale.  (Note that $G \times_{G_p}Z$ is by definition the quotient of $G \times Z$
by the action of $G_p$ given by
$$g \cdot (t,z) = (tg^{-1},g \cdot w) \text{\ for\ }g \in G_p \text{\ and\ }(t,w) \in G \times
Z.)$$
Moreover, if $U$ is nonsingular, then $Z$ can be taken to be also nonsingular and we have the
following commutative diagram of fiber squares
$$\CD
G \times_{G_p}N_p @. \overset{\text{\'etale}}\to{\leftarrow} @. \hskip.25in G \times_{G_p}Z @.
\hskip.25in \overset{\text{\'etale\ surjective}}\to{\rightarrow} @. U_p
\\ 
@VVV \square @.\hskip.25in @VVV \hskip.25in \square @.@VVV \\
N_p//G_p @. \overset{\text{\'etale}}\to{\leftarrow} @. \hskip.25in Z//G_p @.
\hskip.25in \overset{\text{\'etale\ surjective}}\to{\rightarrow} @. U_p//G,\\
\endCD$$
where $N_p$ is the normal vector space to $Z$ at $p$. 
\endproclaim

Let $p \in W$ be a closed point.  By a theorem of Sumihiro, there exists a $K^*$-invariant
affine open neighborhood $U$ of $p$ such that $U$ is also contained in $W - \{\overline{O(p)} -
O(p)\}$.  Now $T = O(p)$ is a closed orbit of $G = K^*$ in $U$, which is nonsingular by
assumption.  Then applying Luna's \'Etale Slice Theorem, we obtain strongly \'etale morphisms
$$\align
i_p:& V_p = G \times_{G_p}Z \rightarrow U_p = G \cdot Z \hskip.1in (\text{surjective}) \\
\eta_p: &V_p = G \times_{G_p}Z \rightarrow X_p = G \times_{G_p}N_p,\\
\endalign$$
where $X_p = G \times_{G_p}N_p = K^* \times_{(K^*)_p}N_p$ is easily seen to be a nonsingular
affine toric variety with $K^*$ acting as a one-parameter subgroup, as it acts on the first
factor of $K^* \times_{(K^*)_p}N_p$ by the usual multiplication.

Thus
$$X_p \overset{\eta_p}\to{\leftarrow} V_p \overset{i_p}\to{\rightarrow} U_p$$
is Luna's locally toric chart at $p$.

Moreover, if $p \in W$ is a fixed point, then $i_p$ as above is an isomorphism, i.e., $V_p =
U_p$ is a $K^*$-invariant affine neighborhood of $p$.  If $p \in W$ is not a fixed point, then
by taking $U$ to be contained in $W - F_W$ where $F_W$ is the fixed point set in $W$, we see
that none of $X_p, V_p$ or $U_p$ has any fixed points by construction. 

\vskip.2in

This completes the proofs A) and B) of Proposition 1-3-4.

\vskip.2in

\S 1-4. Elimination of points of indeterminacy

\vskip.1in

For the weak factorization problem, we start with a general birational map $\phi:X_1
\dashrightarrow X_2$ between complete nonsingular varieties.  In this section, we show, by
the method of elimination of points of indeterminacy, that we only have to deal with the case where
$\phi$ is a projective birational morphism.

\proclaim{Lemma 1-4-1} There is a commutative diagram
$$\CD
X_1' @.\hskip.2in \overset{\phi'}\to{\rightarrow}@.\hskip.2in X_2' \\
@V{g_1}VV @. \hskip.2in @VV{g_2}V \\
X_1@.\hskip.2in \overset{\phi}\to{\dashrightarrow} @.\hskip.2in X_2 \\
\endCD$$
such that $g_1$ and $g_2$ are sequences of blowups with smooth centers disjoint from $U$ and
that $\phi'$ is a projective birational morphism. 
\endproclaim

\demo{Proof}\enddemo By Hironaka's theorem on elimination of points of indeterminacy, there is a
sequence of blowups $g_2:X_2' \rightarrow X_2$ with smooth centers disjoint from $U$ such that the
birational map $h := \phi^{-1} \circ g_2:X_2' \rightarrow X_1$ is a morphism.  By the
same theorem, there is a sequence of blowups $g_1:X_1' \rightarrow X_1$ with smooth centers
disjoint from
$U$ such that the birational map $\phi':X_1' \rightarrow X_2$ is a morphism.  

Since the composite $h \circ \phi'$ is projective, we conclude that $\phi'$ is also projective.

\vskip.1in

We only have to replace $\phi:X_1 \dashrightarrow X_2$ by $\phi':X_1' \rightarrow X_2'$ as above
to see the reduction Step 1 in the strategy for the proof described in Chapter 0. Introduction.

\proclaim{Remark 1-4-2 (Elimination of Points of Indeterminacy)}\endproclaim

We should make a couple of remarks about Hironaka's method of elimination of points of
indeterminacy for a birational map $\phi:X_1 \dashrightarrow X_2$.

First we may assume that $\phi^{-1}$ is a morphism.  Otherwise, replace $W_2$ by the graph of
$\phi$.  (If one wants to stay in the nonsingular category, take a resolution of singularities
of the graph.)

Secondly, if $\phi^{-1}$ is a projective morphism then we can find an ideal $I$, trivial over
$U$, such that $\phi^{-1}$ is the blowup of $X_1$ along $I$.  (See, for example, the proof of
Theorem 2-2-2.)  In general, including the case where $\phi^{-1}$ is not a projective morphism,
Hironaka's version of Chow's Lemma (cf. Hironaka [3]) asserts that there exists an ideal $I$ on $X_1$
such that the blowup of $X_1$ along $I$ factors through $X_2$.  (Remark that although it is not
explicitly stated in Hironaka [3], the ideal $I$ can be taken so that it is trivial over $U$. 
Remark also that although Hironaka [3] works specifically over the field ${\Bbb C}$ of complex
numbers, the canonicity of his method and Lefschetz principle imply that it is actually valid over
any field of characteristic zero (cf. Chapter 4 and Chapter 5).  See also Section 5 in
Abramovich-Karu-Matsuki-W{\l}odarczyk [1].)  We only have to take
$g_1:X_1'
\rightarrow X_1$ to be a principalization of the ideal $I$ by a sequence of blowups with
smooth centers (lying over the support of ${\Cal O}_{X_2}/I$) in order to obtain a morphism $h = \phi
\circ g_1:X_1'
\rightarrow X_2$.

\newpage

$$\bold{CHAPTER\ 2.\ BIRATIONAL\ COBORDISM}$$

\vskip.2in

In this chapter, we study the theory of birational cobordisms by W{\l}odarczyk [2] (See also
Morelli [1].), which is the main tool for our construction of the factorization of a given
birational map
$\phi:X_1 \dashrightarrow X_2$ between nonsingular complete varieties.  The power of the theory
lies in the fact that it allows us to analyze the structure of birational transformations
much like the usual Morse theory allows us to study the structure of homotopy transformations via
the usual cobordisms.  The theory may be considered an algebraic version of
the Morse theory in terms of the action of the multiplicative group $K^*$ through this analogy,
which also justifies the definition we give below.  (See the introduction in Chapter 0.)

It may be worthwhile to note that the connection between the Morse theory and the action of
$K^*$ on an algebraic variety has been known and may even be considered classical, as was
observed by many authors (e.g. Brion-Processi [1] Dolgachev-Hu [1] Frankel [1] Kirwan
[1] Mumford-Fogarty-Kirwan [1] Thaddeus [1] [2]).  When
$K = {\Bbb C}$ is the field of complex numbers, it is also studied in the realm of symplectic
geometry.  Actually we take advantage of the interpretation of the birational cobordisms in terms
of Geometric Invariant Theory to show that, if both $X_1$ and $X_2$ are projective, then we
can choose our factorization so that all the intermediate varieties are also projective.  The
subject of the change of the G.I.T. quotients associated to the change of linearizations has
been found to have a close connection with the factorization problem of certain birational maps,
as studied by Thaddeus and others.

The novelty of the ideas by W{\l}odarczyk, therefore, is in connecting the
longstanding factorization problem of GENERAL birational maps to these classical ideas by
constructing birational cobordisms and hence providing a view point from the Morse theory to
approach the problem.

We cannot probably emphasize too much the importance of the original work of Morelli [1], who
introduced the notion of combinatorial cobordisms to solve the factorization problem of toric
birational maps.  It is the geometric interpretation of Morelli's combinatorial objects that
ultimately led W{\l}odarczyk to the theory of birational cobordism for general birational maps.

\vskip.1in

\S 2-1. Definition of a birational cobordism and the toric main example

\vskip.1in

\proclaim{Definition 2-1-1 (Birational Cobordism)} Let $\phi:X_1 \dashrightarrow X_2$ be a
proper birational map between normal varieties defined over $K$, isomorphic over a common open
subset
$X_1
\supset U \subset X_2$.  A normal variety $B$ is called a $\bold{birational\ cobordism}$ for
$\phi$ and denoted by $B_{\phi}(X_1,X_2)$ if it satisfies the following conditions:

\ \ (i) the multiplcative group $K^*$ acts on $B = B_{\phi}(X_1,X_2)$,  (We denote the action of
$t \in K^*$ on $x \in B$ by $t(x)$ or $t \cdot x$.)

\ \ (ii) the sets
$$\align
B_+ &:= \{x \in B; \lim_{t \rightarrow \infty}t(x) \text{\ does\ NOT\ exist\ in\ }B\} \\ 
B_- &:= \{x \in B; \lim_{t \rightarrow 0}t(x) \text{\ does\ NOT\ exist\ in\ }B\} \\ 
\endalign$$
are nonempty Zariski open subsets of $B$, and

\ \ (iii) there are isomorphisms
$$\align
B_+/K^* &\overset{\sim}\to{\rightarrow} X_2 \\
B_-/K^* &\overset{\sim}\to{\rightarrow} X_1 \\
\endalign$$
so that the birational map induced by the inclusions $B_- \supset B_- \cap B_+ \subset B_+$
and the isomorphisms above
$$X_1 = B_-/K^* \supset (B_- \cap B_+)/K^* \subset B_+/K^* = X_2$$
coincides with $\phi$.

We say that $\bold{B\ respects\ the\ open\ subset\ U}$ if $U$ is contained in $B_- \cap B_+/K^*$.
\endproclaim

We discuss the following fundamental example of a birational cobordism in the toric setting, as
was first observed by Morelli [1].

\proclaim{Main Toric Example 2-1-2}\endproclaim

Let $B = {\Bbb A}^n = \roman{Spec}\ K[z_1, \cdot\cdot\cdot, z_n]$ with a $K^*$-action given by
$$t(z_1, \cdot\cdot\cdot, z_j, \cdot\cdot\cdot, z_n) = (t^{\alpha_1}z_1, \cdot\cdot\cdot,
t^{\alpha_j}z_j, \cdot\cdot\cdot, t^{\alpha_n}z_n).$$
We regard ${\Bbb A}^n = X(N,\sigma)$ as a toric variety defined by a lattice $N \cong {\Bbb
Z}^n$ and a regular cone $\sigma \subset N_{\Bbb R}$ generated by the standard ${\Bbb Z}$-basis of
$N$
$$\sigma = \langle v_1, \cdot\cdot\cdot, v_j, \cdot\cdot\cdot, v_n\rangle.$$
The dual cone $\sigma^{\vee}$ is generated by the dual ${\Bbb Z}$-basis
$$\sigma^{\vee} = \langle v_1^*, \cdot\cdot\cdot, v_j^*, \cdot\cdot\cdot, v_n^*\rangle$$
and we identify
$$z_j = z^{v_j^*}.$$
The $K^*$-action then corresponds to the one parameter subgroup
$$a = (\alpha_1, \cdot\cdot\cdot, \alpha_j, \cdot\cdot\cdot, \alpha_n) \in N.$$
We have then the obvious description of the sets $B_+$ and $B_-$
$$\align
B_+ &= \{(z_1, \cdot\cdot\cdot, z_n);z_j \neq 0 \text{\ for\ some\ }j \text{\ with\ }\alpha_j =
(v_j^*,a) > 0\} \\
B_- &= \{(z_1, \cdot\cdot\cdot, z_n);z_j \neq 0 \text{\ for\ some\ }j \text{\ with\ }\alpha_j =
(v_j^*,a) < 0\}. \\
\endalign$$
We define the upper boundary and lower boundary of $\sigma$ (with respect $a \in N$) to be
$$\align
\partial_+\sigma &= \{x \in \sigma; x + \epsilon \cdot (-a) \not\in \sigma \text{\ for\
}\epsilon > 0\}\\ 
\partial_-\sigma &= \{x \in \sigma; x + \epsilon \cdot a \not\in \sigma \text{\ for\
}\epsilon > 0\}.\\ 
\endalign$$
Then we obtain the description of $B_+$, $B_-$ and $B_+ \cap B_-$ as the toric varieties
corresponding to the fans $\partial_+\sigma$, $\partial_-\sigma$ and $\partial_+\sigma \cap
\partial_-\sigma$, i.e.,
$$\align
B_+ &= X(N,\partial_+\sigma) \\
B_- &= X(N,\partial_-\sigma). \\
B_+ \cap B_- &= X(N,\partial_+\sigma \cap \partial_-\sigma) \\
\endalign$$
Accordingly, if we denote by $\pi:N_{\Bbb R} \rightarrow N_{\Bbb R}/{\Bbb R} \cdot a$ the
projection onto $N_{\Bbb R}/{\Bbb R} \cdot a$ with the lattice $\pi(N)$, then we have the
description of the quotients as the toric varieties
$$\align
B_+/K^* &= X(\pi(N),\pi(\partial_+\sigma)) \\ 
B_-/K^* &= X(\pi(N),\pi(\partial_-\sigma)) \\ 
(B_+ \cap B_-)/K^* &= X(\pi(N),\pi(\partial_+\sigma \cap \partial_-\sigma)) \\
B//K^* &= X(\pi(N),\pi(\sigma)) \\
\endalign$$
where $\pi(\partial_+\sigma), \pi(\partial_-\sigma), \pi(\partial_+\sigma \cap \partial_-\sigma)$
and
$\pi(\sigma)$ are fans in
$N_{\Bbb R}/{\Bbb R} \cdot a$ defined by
$$\align
\pi(\partial_+\sigma) &= \{\pi(\tau);\tau \in \partial_+\sigma\} \\
\pi(\partial_-\sigma) &= \{\pi(\tau);\tau \in \partial_-\sigma\} \\
\pi(\partial_+\sigma \cap \partial_-\sigma) &= \{\pi(\tau);\tau \in \partial_+\sigma \cap
\partial_-\sigma\} = \pi(\partial_+\sigma) \cap \pi(\partial_-\sigma) \\
\pi(\sigma) &= \{\pi(\sigma) \text{\ and\ its\ proper faces\ in\ }N_{\Bbb R}/{\Bbb R} \cdot
a\},\\
\endalign$$
respectively.  Since $\pi(\partial_+\sigma)$ and $\pi(\partial_-\sigma)$ are subdivisions of
$\pi(\sigma)$, we obtain a diagram of toric birational maps

$$\CD
B_-/K^* @.@.\overset{\phi}\to{\dashrightarrow}@.@.B_+/K^* \\ 
@.\searrow \hskip.2in@.
\hskip.6in @.\hskip.2in \swarrow @.\\ 
@| @.B//K^* @.@.@| \\
X(\pi(N),\pi(\partial_-\sigma)) @.@.\dashrightarrow @.@.X(\pi(N),\pi(\partial_+\sigma)) \\
@.\searrow \hskip.2in@. @| \hskip.2in \swarrow @.\\
@.@.X(\pi(N),\pi(\sigma)) @.@.\\
\endCD$$
where $\phi$ coincides with the birational map induced by
$$\CD
B_-/K^* @. \supset @. (B_- \cap B_+)/K^* @. \subset @. B_+/K^* \\
@| @.@| @. @| \\
X(\pi(N),\pi(\partial_-\sigma)) \hskip.2in @. \supset @.\hskip.2in 
X(\pi(N),\pi(\partial_-\sigma
\cap
\partial_+\sigma)) \hskip.2in @. \subset @.\hskip.2in  X(\pi(N),\pi(\partial_+\sigma)).\\
\endCD$$
Therefore, we conclude that $B$ is a birational cobordism for $\phi$.

More generally, one can prove (See Morelli [1] Abramovich-Matsuki-Rashid [1].) that if
$\Sigma$ is a subdivision of a convex polyhedral cone in $N_{\Bbb R}$ with the lower boundary
$\partial_-\Sigma$ and upper boundary $\partial_+\Sigma$ with respect to a one-parameter
subgroup $a \in N$, then the toric variety $X(N,\Sigma)$ corresponding to $\Sigma$ with the
$K^*$-action given by the one-parameter subgroup $a \in N$, is a birational cobordism between
the two toric varieties $X(\pi(N),\pi(\partial_-\Sigma))$ and
$X(\pi(N),\pi(\partial_-\Sigma))$.

\vskip.1in

Next in the pursuit of analogy to the usual Morse theory where the critical points are ``lined
up nicely" according to the levels given by the Morse function, we would like to have the fixed
points of a birational cobordism ``ordered nicely" so that we can study the birational
transformations as we go through the fixed points ``from the bottom to the top", though we may
not see a Morse function explicitly.  This requirement naturally leads us to the notion of
$\bold{``collapsibility"}$ introduced by Morelli [1].  In the toric setting that Morelli [1]
studied, the fixed points of the cobordism, associated to a fan $\Sigma$, correspond to the cones
$\sigma$ of the maximal dimension (bubbles) in $\Sigma$.  Factorization then corresponds to going
from
$\partial_-\Sigma$ to $\partial_+\Sigma$, by replacing $\partial_-\sigma$ with
$\partial_+\sigma$ one at a time (collapsing of the bubble $\sigma$).  (See the main toric
example 2-1-2 and Morelli [1] Abramovich-Matsuki-Rashid [1].)  Thus the question of whether or not
we can order the fixed points nicely corresponds to that of whether or not we can collapse these
bubbles in a nicely ordered manner, thus giving rise to the name ``collapsibility". 
 
\vskip.1in

First we introduce the notations for some specific subsets in $B$ associated to the
fixed point set.

\proclaim{Notation 2-1-3} Let $B = B_{\phi}(X_1,X_2)$ be a birational cobordism, and let $F
\subset  F_B = B^{K^*}$ be a subset of the fixed point set with respect to the $K^*$-action. 
We define
$$\align
F^+ &:= \{x \in B;\lim_{t \rightarrow 0}t(x) \in F\} \\
F^- &:= \{x \in B;\lim_{t \rightarrow \infty}t(x) \in F\} \\
F^{\pm} &:= F^+ \cup F^- \\
F^* &:= F^{\pm} - F.\\
\endalign$$
\endproclaim 

\proclaim{Definition 2-1-4} Let $B = B_{\phi}(X_1,X_2)$ be a birational cobordism.  W define a
relation $\prec$ among the connected components of the fixed point set $B^{K^*}$ as follows:
let $F_1, F_2 \subset B^{K^*}$ be two (not necessarily distinct) connected components, and set
$F_1 \prec F_2$ if there exists a point $p \in B$ with $p \not\in B^{K^*}$ such that
$$\lim_{t \rightarrow 0}t(p) \in F_1 \text{\ and\ }\lim_{t \rightarrow \infty}t(p) \in F_2.$$
That is to say, in terms of Notation 2-1-3, we have the relation $F_1 \prec F_2$ if and only
if $(F_1^+ - F_1) \cap (F_2^- - F_2) \neq \emptyset$. 
\endproclaim

\proclaim{Definition 2-1-5 (Collapsibility)} We say that a birational cobordism $B =
B_{\phi}(X_1,X_2)$ is collapsible if the relation $\prec$ is a strict preorder, namely, there
is no directed cycle of the connected components of the fixed point set
$$F_1 \prec F_2 \prec \cdot\cdot\cdot \prec F_m \prec F_1.$$
Note that collapsibility excludes the possibility of a self-loop $F_1 \prec F_1$ for any
connected component $F_1 \subset B^{K^*}$. 
\endproclaim

In Chapter 3, we will consider the torific ideals on a birational cobordism $B$ where, however,
the existence of a point $p \in B$ with $p \not\in B^{K^*}$ such that both limits $\lim_{t
\rightarrow 0}t(p)$ and $\lim_{t \rightarrow \infty}t(p)$ exist within $B$ would cause a
problem in order to have such ideals well-defined.  The notion of a quasi-elemetary
birational cobordism below is introduced exactly to avoid this problem.

\proclaim{Definition 2-1-6 (Quasi-Elementary Birational Cobordism)} A birational cobordism $B$
is said to be quasi-elementary if there does not exist any point $p \in B$ with $p \not\in B^{K^*}$ such that both limits $\lim_{t
\rightarrow 0}t(p)$ and $\lim_{t \rightarrow \infty}t(p)$ exist within $B$.  In terms of the
relation $\prec$, this is equivalent to saying that for any two (not necessarily distinct)
connected components $F_1, F_2 \subset B^{K^*}$ neither $F_1 \prec F_2$ nor $F_2 \prec F_1$
holds.  That is to say, in terms of Notation 2-1-3, $B$ is quasi-elementary if and only if $(F^+
- F) \cap (F^- - F) = \emptyset$ where $F = B^{K^*}$ is the entire fixed point set in $B$.
\endproclaim

\proclaim{Definition 2-1-7 (Elementary Birational Cobordism)} A quasi-elementary birational
cobordism is said to be elementary if the fixed point set $B^{K^*}$ is connected.
\endproclaim

We will observe that the birational transformation represented by an elementary birational
cobordism corresponds to a (weighted) blowup of a connected center, which is \'etale locally
equivalent to a toric birational transformation.  (See \S 2-4 for the details.)  One might want
to call such a birational transformation ``elementary", as W{\l}odarczyk [2] does, and hence
name the corresponding birational cobordism also ``elementary".

\vskip.2in

\S 2-2. Construction of a (collapsible) birational cobordism

\vskip.1in

After defining a birational cobordism, the natural and important issue is its existence,
which W{\l}odarczyk [2] shows for any proper birational map between (normal) varieties.  Here
we present a simple construction of a collapsible birational cobordism for a projective
birational morphism $\phi:X_1 \rightarrow X_2$ between complete nonsingular varieties, which
suffices for our purposes after reducing the factorization problem of an arbitrary
birational map to that of a projetive morphism via Hironaka's theorem elimination of
indeterminacy.  See \S 2-4 for the details of the reduction step.

\proclaim{Theorem 2-2-2 (Construction of Birational Cobordism)} Let $\phi:X_1 \rightarrow
X_2$ be a projective (in the sense of Grothendieck and not as defined in Hartshorne [1])
birational morphism between complete nonsingular varieties (defined over $K$), which is an
isomorphism over a common open subset $U$.  Then there exists a complete nonsingular variety
$\overline{B}$ with an effective $K^*$-action, (Note that a $K^*$-action is called
$\bold{effective}$ if the action is not induced from that of the nontrivial quotient of $K^*$,
that is to say, if $\cap_{p \in \overline{B}}Stab(p) = \{1\}$.) satisfying the following
properties: 

\ \ (i) there exist closed embeddings
$$\align
i_1:X_1 &\hookrightarrow \overline{B} \\
i_2:X_2 &\hookrightarrow \overline{B} \\
\endalign$$
with disjoint images in $\overline{B}^{K^*}$,

\ \ (ii) there is a coherent sheaf ${\Cal E}$ on $X_2$, with a $K^*$-action (which is
compatible with the trivial action of $K^*$ on $X_2$), and a $K^*$-equivariant closed embedding
$$\overline{B} \hookrightarrow {\Bbb P}({\Cal E}) := Proj_{X_2} \oplus_{m \geq 0}
Sym^m{\Cal E},$$

\ \ (iii) the open subvariety $B := \overline{B} - (i_1(X_1) \cup i_2(X_2))$ is a
collapsible birational cobordism for $\phi$.

We call such a variety $\overline{B}$ a $\bold{compactified\ birational\ cobordism}$ projective
over
$X_2$.
\endproclaim

\proclaim{Remark 2-2-3}\endproclaim

Suppose that a projective birational morphism between complete nonsingular varieties
$\phi:X_1 \rightarrow X_2$ is a sequence of blowups with smooth centers.  Then we can
construct a compatified birational cobordism
$\overline{B}$ projective over $X_2$ in the following simple manner: We start with the
product $W_0 = X_2 \times {\Bbb P}^1$ where $K^*$ acts on the second factor as the
multiplication on $K^* = {\Bbb P}^1 - \{0,\infty\}$.  We take the sequence of blowups $W
\rightarrow X_2
\times {\Bbb P}^1$ with the centers identified with those for $\phi$ but considered as lying in
(the strict transforms of) the $0$-section $X_2 \times \{0\}$ instead of lying in (the
subsequent blowups of) $X_2$.  Now
$\overline{B} = W$ is a complete nonsingular variety projective over $X_2$, with an effective
$K^*$-action, satisfying the properties:

\ \ (i) we have two closed embeddings
$$\align
i_1:X_1 &\overset{\sim}\to{\rightarrow} (\text{the\ strict\ transform\ of\ }X_2 \times \{0\})
\subset \overline{B} \\
i_2:X_2 &\overset{\sim}\to{\rightarrow} X_2 \times \{\infty\} \subset \overline{B} \\
\endalign$$
with disjoint images $i_1(X_1), i_2(X_2) \subset \overline{B}^{K^*}$ 

\ \ (ii) the centers of blowups lie in the fixed point set of the $K^*$-action and hence the
morphism $\overline{B} \rightarrow X_2 \times {\Bbb P}^1 \rightarrow X_2$ is obviously
$K^*$-equivariant, which easily implies the existence of such a coherent sheaf ${\Cal E}$ as
above,

\ \ (iii) $B = \overline{B} - (i_1(X_1) \cup i_2(X_2))$ is a collapsible birational cobordism
for $\phi$.  Remark that the factorization induced by $B$, as will be discussed in \S 2-4,
coincides with the given sequence of blowups with smooth centers.

\vskip.1in

The above birational cobordism can be considered as the standard birational cobordism
associated to a sequence of blowups with smooth centers factoring a projective birational
morphism $\phi$ (if such a sequence exists at all for $\phi$).  The construction below for the
proof of Theorem 2-2-2 may be considered as modeled on that of the standard birational
cobordism discussed above. 

\vskip.1in

\demo{Proof of Theorem 2-2-2}\enddemo We start with the product $W_0 = X_2 \times {\Bbb P}^1$ where $K^*$ acts on the second factor as the
multiplication on $K^* = {\Bbb P}^1 - \{0,\infty\}$.

Since $\phi$ is a projective birational morphism which
is an isomorphism over $U$, there exists an ideal sheaf $J \subset {\Cal O}_{X_2}$ such that
$\phi:X_1
\rightarrow X_2$ is the blowup morphism of $X_2$ along $J$ and that the support of ${\Cal
O}_{X_2}/J$ is disjoint from
$U$.  (This can be verified as follows: Take a divisor $D$ on $X_1$ which is $\phi$-ample. 
Then we see $D - \phi^*\phi_*D$ is also $\phi$-ample and that it is of the form $\Sigma
(-a_i)E_i$ where $f_*D$ is the cycle-theoretic image of $D$ under $\phi$ as a divisor and where
the
$E_i$ are exceptional divisors with
$a_i > 0$.  We only have to take
$J =
\phi_*{\Cal O}_{X_1}(l \cdot \Sigma (-a_i)E_i)$ for some sufficiently large $l \in {\Bbb N}$. 
Note that in our case the existence of such a $\phi$-ample divisor $\Sigma (-a_i)E_i$ follows
from the construction.  (cf. \S 5-1.)  Remark also that in our construction of the
factorization the projective morphism $\phi$ is obtained as $\phi':X_1' \rightarrow
X_2'$ starting from the original birational map $\phi:X_1 \dashrightarrow X_2$ via
the elimination of the points of indeterminacy.  Thus in the notation of Lemma
1-4-1 we have a $g_1$-ample divisor of the form $\Sigma (-b_j)F_j$, where $F_j$ are
$g_1$-exceptional divisors whose centers on $X_1$ lie outside of $U$.  The divisor $\Sigma
(-b_j)F_j$ is hence $\phi$-ample.  We can take $J = \phi_*{\Cal O}_{X_1}(l \cdot
\Sigma (-b_j)F_j)$ for some sufficiently large $l \in {\Bbb N}$.)  Let
$I_0$ be the ideal of the origin
$0
\in {\Bbb P}^1$.  We set
$$I = (p_1^{-1}J + p_2^{-1}I_0) \cdot {\Cal O}_{W_0}$$
where
$$\align
p_1:W_0 = X_2 \times {\Bbb P}^1 &\rightarrow X_1 \\
p_2:W_0 = X_2 \times {\Bbb P}^1 &\rightarrow {\Bbb P}^1 \\
\endalign$$ 
are the projections.

Let $W$ be the blowup of $W_0$ along $I$.

We claim that $X_1$ and $X_2$ are embedded in the nonsingular locus of $W$, as the strict
transform of
$X_2 \times \{0\}$ and $X_2 \times \{\infty\}$, respectively.  For $i_2:X_2
\overset{\sim}\to{\rightarrow} X_2
\times
\{\infty\} \subset W$ this is clear, as the centers of blowup for $X_2 \times {\Bbb P}^1$ only
lie over $X_2 \times \{0\}$.  For $i_2:X_1 \overset{\sim}\to{\rightarrow} (\text{the\ strict\
transform\ of\ }X_2 \times \{0\})$, in order to prove that $X_1$ is in the nonsingular locus
of $W$, it suffices to show that $X_1$ is a Cartier divisor as $X_1$ itself is nonsingular. 
We look at local coordinates.  Let $A(V)$ be the affine coordinate ring for an affine open
subset $V \subset X_2$ and let $y_1, \cdot\cdot\cdot, y_m$ be a set of generators for $J$ on
$V$.  Let $K[x]$ be the affine coordinate ring for ${\Bbb A}^1 = {\Bbb P}^1 - \{\infty\}$. 
Then the affine open subset $V \times {\Bbb A}^1 \subset X_2 \times {\Bbb P}^1$ has the affine
coordiante ring $A(V) \otimes_KK[x]$, where the ideal $I$ is generated by $y_1,
\cdot\cdot\cdot, y_m, x$.  The charts of the blowup containing th strict transform of $V
\times \{0\} = \{x = 0\}$ are of the form
$$\roman{Spec}\ A[\frac{y_1}{y_i}, \cdot\cdot\cdot, \frac{y_m}{y_i},\frac{x}{y_i}] =
\roman{Spec}\ A[\frac{y_1}{y_i}, \cdot\cdot\cdot, \frac{y_m}{y_i}] \times
\roman{Spec}\ K[\frac{x}{y_i}],$$
where $K^*$ acts on the second factor by the multiplication on $x$.  The strict transform of
$\{x = 0\}$ is defined by $\frac{x}{y_i}$ and hence is Cartier.

Let $\overline{B} \rightarrow W$ be the canonical resolution of singularities.  Since the
centers of blowups of the canonical resolution of singularities are taken over the singular
locus (cf. the condition $(\spadesuit^{res}-0)$) and since the (effective) $K^*$-action lifts to
$\overline{B}$ (cf. the condition $(\spadesuit^{res}-1)$), the above analysis of
$W$ immediately implies the properties (i) and (ii), except for the collapsibility of $B$.

Remark that the morphism $\tau:\overline{B} \rightarrow X_2$ is $K^*$-equivariant as well as
projective.  Therefore, there exists a relatively ample line bundle ${\Cal L}$ on
$\overline{B}$ over $X_2$, equipped with a $K^*$-action (compatible with the $K^*$-action on
$\overline{B}$).  (Remark that in our case the existence of such a relatively ample line bundle
follows directly from the construction (cf. \S 5-1).)  In order to see the property (iii), we
only have to set
$${\Cal E} = \tau_*({\Cal L}^{\otimes l}) \text{\ for\ sufficiently\ large\ }l \in {\Bbb N}.$$

\vskip.1in

Finally we show that the birational cobordism $B$ constructed as above is collapsible.

\vskip.1in

Let ${\Cal C}$ be the set of the connected components of $B^{K^*}$.  It suffices to show that
there exists a strictly increasing function $\chi:{\Cal C} \rightarrow {\Bbb Z}$, i.e., a
function such that
$$F_1 \prec F_2 \longrightarrow \chi(F_1) < \chi(F_2).$$
Since $K^*$ acts trivially on $X_2$ and since $K^*$ is reductive, there exists a direct sum
decomposition
$${\Cal E} = \oplus_{b \in {\Bbb Z}}{\Cal E}_b$$
where ${\Cal E}_b$ is the subsheaf on which the $K^*$-action is given by the character $t
\mapsto t^b$.  Denote by 
$$b_0 < b_1 < \cdot\cdot\cdot < b_{h-1} < b_h$$
the characters which
appear in this representation.  Note that there are disjoint embeddings
${\Bbb P}({\Cal E}_b)
\subset {\Bbb P}({\Cal E})$.

Let $p \in B$ be a point lying in the fiber ${\Bbb P}(E_q)$ over $q \in X_2$, where
$E_q = {\Cal E} \otimes {\Cal O}_{X_2,q}/m_{X_2,q}$.  We choose a basis of $E_q$, which are
then considered to be homogeneous coordinates of ${\Bbb P}(E_q)$
$$(x_{b_0,1}, \cdot\cdot\cdot, x_{b_0,d_0}, \cdot\cdot\cdot, x_{b_h,1}, \cdot\cdot\cdot,
x_{b_h,d_h})$$
where
$x_{b_j,k}
\in E_{b_j} = {\Cal E}_{b_j}
\otimes {\Cal O}_{X_2,q}/m_{X_2,q}$.

It is straightforward to observe that $p \in B$ is a fixed point if and only if
$$\exists b_{j_p} \text{\ such\ that\ }x_{b_j,k}(p) = 0 \text{\ whenever\ }b_j \neq b_{j_p}.$$
Therefore, if $F \subset B^{K^*}$ is a connected component, then $F$ is contained in ${\Bbb
P}({\Cal E}_{b_j})$ for some $b_j$.  We define
$$\chi(F) := b_j.$$
Suppose $F_1 \prec F_2$, i.e., there exists a point $p \in B$ but $p \in B^{K^*}$ with
$$\lim_{t \rightarrow 0}t(p) \in F_1 \text{\ and\ }\lim_{t \rightarrow \infty}t(p) \in F_2.$$
It is easy to see that
$$\align
\lim_{t \rightarrow 0}t(p) &\in {\Bbb P}(E_{b_{min}}) \\
\lim_{t \rightarrow 0}t(p) &\in {\Bbb P}(E_{b_{max}}), \\
\endalign$$
where
$$\align
b_{min} &= \min\{b_j;x_{b_j,k}(p) \neq 0 \text{\ for\ some\ }k\} \\
b_{max} &= \max\{b_j;x_{b_j,k}(p) \neq 0 \text{\ for\ some\ }k\}. \\
\endalign$$
Observe also that
$$b_{min} < b_{max}$$
since $p \in B$ is not a fixed point.  Thus we conclude that
$$\chi(F_1) = b_{min} < b_{max} = \chi(F_2),$$
showing $\chi$ is a strictly increasing function.

Therefore, the birational cobordism $B$ is collapsible.

\vskip.1in

This completes the proof of Theorem 2-2-2.

\vskip.2in

\S 2-3. Interpretation by Geometric Invariant Theory

\vskip.1in

In this section, we apply Geometric Invariant Theory to the compactified birational
cobordism $\overline{B}$ projective over $X_2$ and see that the birational cobordism
$B$ decomposes into the quasi-elementary pieces $B_{a_i}$ (See below for the precise definition),
which are the semi-stable loci of the different linearizations of the $K^*$-action on
$\overline{B}$.  Thus we can analyze the birational transformations of $B_{a_i}//K^*$ as the
change of the G.I.T. quotients when we vary the linearizations, a subject which is extensively
studied by Thaddeus [1][2] and others.

\vskip.1in

We continue using the notation of Theorem 2-2-2.

\vskip.1in

Let ${\Cal E}$ be a coherent sheaf on $X_2$, with a $K^*$-action (which is
compatible with the trivial action of $K^*$ on $X_2$), and a $K^*$-equivariant closed embedding
$$\overline{B} \hookrightarrow {\Bbb P}({\Cal E}) := Proj_{X_2} \oplus_{m \geq 0}
Sym^m{\Cal E}$$
and let
$${\Cal E} = \oplus_{b \in {\Bbb Z}}{\Cal E}_b$$
be its decomposition into the eigen-subsheaves according to their characters.  Let
$$a_1 < a_2 < \cdot\cdot\cdot < a_{s-2} < a_{s-1}$$
be the values of $\chi:{\Cal C} \rightarrow {\Bbb Z}$, i.e., $\{a_i\}$ is the subset of
$\{b_j\}$, which are the values associated to the connected components of the fixed point set
$B^{K^*}$.  By using the Veronese embedding ${\Bbb P}({\Cal E}) \hookrightarrow {\Bbb
P}(Sym^2{\Cal E})$ and replacing ${\Cal E}$ with $Sym^2{\Cal E}$, we may assume that all the
$a_i$ are even and that in particular 
$$a_{i-1} < a_i - 1 < a_i < a_i + 1 < a_{i+1} \text{\ for\ } i = 2, \cdot\cdot\cdot, s-2.$$ 
(This is a technical condition which avoids the rational ``twists" and allows us only to deal
with the integral ones.)

Denote by
$$\rho_0(t)$$
the original $K^*$-action on ${\Cal E}$.  For any $r \in {\Bbb Z}$ we consider the ``twisted"
action of $K^*$ on ${\Cal E}$, denoted by
$$\rho_r(t) = t^{-r} \cdot \rho_0(t),$$
so that in the decomposition ${\Cal E} = \oplus_{b \in {\Bbb Z}}{\Cal E}_b$ the $\rho_r(t)$
acts on ${\Cal E}_b$ by the multiplication $t^{b - r}$.  Note that the ``twists" do not change
the induced action of $K^*$ on ${\Bbb P}({\Cal E})$ but only change the linearizations on
${\Cal O}_{{\Bbb P}({\Cal E})}(1)$.

We can now apply Geometric Invariant Theory in its relative form (cf. Hu [1] Pandripande [1])
to the situation.

We denote by
$$({\Bbb P}({\Cal E}),\rho_r)^{ss}$$
the semi-stable locus of ${\Bbb P}({\Cal E})$ with respect to the linearization of ${\Cal
O}_{{\Bbb P}({\Cal E})}(1)$ induced by the twisted action $\rho_r$ of $K^*$ on ${\Cal E}$. 
Recall that a point $p \in {\Bbb P}({\Cal E})$ is semi-stable, i.e., $p \in ({\Bbb P}({\Cal
E}),\rho_r)^{ss}$ if there exists a local section $s \in Sym^n{\Cal E}$ for some positive $n
\in {\Bbb N}$, invariant under the twisted action $\rho_r$ of $K^*$ on $Sym^n{\Cal E}$, such
that $s(p) \neq 0$.  We denote by
$$(\overline{B},\rho_r)^{ss}$$
the semi-stable locus of $\overline{B}$ with respect to the linearization of ${\Cal
O}_{{\Bbb P}({\Cal E})}(1)|_{\overline{B}}$ coming from that of ${\Cal
O}_{{\Bbb P}({\Cal E})}(1)$ induced by the twisted action $\rho_r$ of $K^*$.  Then we have
$$(\overline{B},\rho_r)^{ss} = \overline{B} \cap ({\Bbb P}({\Cal E}),\rho_r)^{ss}$$
by Theorem 1.19 in Kirwan-Fogarty-Mumford [1].

Now Thaddeus and others tell us that we should observe the change of the semi-stable loci when
the twists pass through critical values and that we should also observe the birational
transformations of the G.I.T. quotients.

In our situtaion, these critical values turn out to be just
$$a_1 < a_2 < \cdot\cdot\cdot < a_{s-2} < a_{s-1}$$
and the change of the semi-stable loci provides the decomposition of the birational cobordism
$B$ into the quasi-elementary pieces.

\proclaim{Proposition 2-3-1} Let the situation be as above.  Then

\ \ (o) the critical values for the change of the semi-stable loci are
$$a_1 < a_2 < \cdot\cdot\cdot < a_{s-2} < a_{s-1},$$

\ \ (i) at each critical value $a_i$ the semi-stable locus
$(\overline{B},\rho_{a_i})^{ss} := B_{a_i}$ is a quasi-elementary birational cobordism and we
have set-theoretic descriptions
$$\align
B_{a_i} &= (\overline{B},\rho_{a_i})^{ss} \\
&= B - \{(\cup_{F \in {\Cal C},\chi(F) < a_i}F^-) \cup
(\cup_{F \in {\Cal C},\chi(F) > a_i}F^+)\} \\
(B_{a_i})_+ &= (\overline{B},\rho_{a_i + 1})^{ss} \\
(B_{a_i})_- &= (\overline{B},\rho_{a_i - 1})^{ss}, \\
\endalign$$

\ \ (ii) we have the commutative diagarm of (proper) birational morphsims and maps
between varieties PROJECTIVE over $X_2$

$$\CD
(B_{a_i})_-/K^* @.@.\overset{\varphi_i}\to{\dashrightarrow}@.@.(B_{a_i})_+/K^* \\ 
@.\searrow \hskip.2in@.
\hskip.6in @.\hskip.2in \swarrow @.\\ 
@| @.B_{a_i}//K^* @.@.@| \\
(\overline{B},\rho_{a_i - 1})^{ss}/K^* @.@.\dashrightarrow @.@.(\overline{B},\rho_{a_i +
1})^{ss}/K^*\\ 
@.\searrow \hskip.2in@. @| \hskip.2in \swarrow @.\\
@.@. (\overline{B},\rho_{a_i})^{ss}//K^* @.@.\\
\endCD$$

\ \ (iii) we also have
$$\align
(B_{a_1})_- &= B_- \text{\ and\ }(B_{a_1})_-/K^* = B_-/K^* = X_1 \\
(B_{a_{s-1}})_+ &= B_+ \text{\ and\ }(B_{a_{s-1}})_+/K^* = B_+/K^* = X_2. \\
\endalign$$

\endproclaim

\demo{Proof}\enddemo The proof is an easy consequence of the (relative) G.I.T. and the
definition of the semi-stable locus as above and left to the reader as an exercise.

\vskip.2in

\S 2-4. Factorization into locally toric birational maps

\vskip.1in

As a summary of the results discussed in the previous sections \S 2-1, \S 2-2 and \S 2-3 and
the use of elimination of points of indeterminacy by Hironaka, we obtain the following
factorization theorem of a birational map into locally toric transformations by W{\l}odarczyk
[2].

\proclaim{Theorem 2-4-1 (Factorization into Locally Toric Birational Maps)} Let $\phi:X_1
\rightarrow X_2$ be a birational map between complete nonsingular varieties over $K$, which
is an isomorphism over a common open subset $X_1 \supset U \subset X_2$.  Then there exists a
sequence of birational maps

$$X_1 = W_1 \overset{\varphi_1}\to{\dashrightarrow} W_2
\overset{\varphi_2}\to{\dashrightarrow} \cdot\cdot\cdot
\overset{\varphi_{i-1}}\to{\dashrightarrow} W_i \overset{\varphi_i}\to{\dashrightarrow} W_i
\overset{\varphi_{i+1}}\to{\dashrightarrow} \cdot\cdot\cdot
\overset{\varphi_{s-1}}\to{\dashrightarrow} W_s = X_2$$

such that

\ \ (i) $\phi = \varphi_{s-1} \circ \varphi_{s-2} \circ \cdot\cdot\cdot \varphi_2 \circ
\varphi_1$,

\ \ (ii) $\varphi_i$ is V-locally toric and an isomorphism over $U$, and that

\ \ (iii) there exists an inedx $i_o$ with the property that for all $i \leq i_o$ the map $W_i
\dashrightarrow X_1$ is a projective birational morphism and for all $i \geq i_o$ the map $W_i
\dashrightarrow X_2$ is a projective birational morphism, and hence if $X_1$ and $X_2$ are
projective then all the $W_i$ are projective.
\endproclaim

\demo{Proof}\enddemo By the reduction step as discussed in \S 1-4, we may assume that $\phi:X_1
\rightarrow X_2$ is a projective birational morphism.  Thus we are in a situation to apply
Theorem 2-2-2 to obtain a compactified birational cobordism
$\overline{B} \subset {\Bbb P}({\Cal E})$ projective over $X_2$, where ${\Cal E}$ is a coherent
sheaf on $X_2$ with a $K^*$-action.  The associated birational cobordism
$B = \overline{B} - (i_1(X_1) \cup i_2(X_2))$ is collapsible and the connected components $F$ of
its fixed points are ordered by the values of a strictly increasing function $\chi:{\Cal C}
\rightarrow {\Bbb Z}$
$$a_1 < a_2 < \cdot\cdot\cdot < a_{s-2} < a_{s-1}.$$
By looking at the semi-stable loci with respect to the linearizations induced by the ``twisted"
actions of $K^*$ on ${\Cal E}$ and at the corresponding G.I.T. quotients, we obtain a
factorization
$$\CD
W_i = (B_{a_i})_-/K^* @.@.\overset{\varphi_i}\to{\dashrightarrow}@.@. W_{i+1} = (B_{a_i})_+/K^*
\\  @.\searrow \hskip.2in@.
\hskip.6in @.\hskip.2in \swarrow @.\\ 
@| @.B_{a_i}//K^* @.@.@| \\
(\overline{B},\rho_{a_i - 1})^{ss}/K^* @.@.\dashrightarrow @.@.(\overline{B},\rho_{a_i +
1})^{ss}/K^*\\ 
@.\searrow \hskip.2in@. @| \hskip.2in \swarrow @.\\
@.@. (\overline{B},\rho_{a_i})^{ss}//K^* @.@.\\
\endCD$$
starting with $X_1 = B_-/K^* = (B_{a_1})_-/K^* = W_1$ and ending with $W_l = (B_{a_l})_+/K^* =
B_+/K^* = X_2$, as given in Proposition 2-3-1.  By construction, each $W_i$ is projective over
$X_2$ and $\varphi_i$ is an isomorphism over $U$.  The existence of the index $i_o$ as claimed
in (iii) is also clear, once one takes the reduction in Step 0 into consideration. 

\vskip.1in

Thus it remains to show that each $\varphi_i$ is a V-locally toric birational map.

\proclaim{Lemma 2-4-4} Let $B_{a_i}$ be a quasi-elementary cobordism.  Then for every closed
point $p \in B_{a_i} - F_{B_{a_i}}^*$, where $F_{B_{a_i}} = B_{a_i}^{K^*}$ and $F_{B_{a_i}}^*$ is
the set as defined in Notation 2-1-3, there exists Luna's locally toric chart
$$X_p \overset{\eta_p}\to{\leftarrow} V_p \overset{i_p}\to{\rightarrow} U_p \subset B_{a_i}$$
such that the following condition $(\star)$ is satisfied:

$(\star)$ If a $K^*$-orbit $O(q)$ lies in $U_p$, then its closure $\overline{O(q)}$ in $B_{a_i}$
also lies in $U_p$.
\endproclaim

\demo{Proof}\enddemo We take
Luna's locally toric chart 
$$X_p \overset{\eta_p}\to{\leftarrow} V_p \overset{i_p}\to{\rightarrow} U_p \subset B_{a_i}$$
as constructed in Proposition 1-3-4, which may not satisfy the
condition $(\star)$.  Let $\pi:B_{a_i} \rightarrow B_{a_i}//K^*$ be the quotient map.  

\vskip.1in

Case: $p \in B_{a_i} - F_{B_{a_i}}^*$ is not a fixed point.

\vskip.1in

In this case, $p \in B_{a_i} - \pi^{-1}(\pi(F_{B_{a_i}})$.  Thus by shrinking $U_p$ if
necessary, we may assume $U_p \subset B_{a_i} - \pi^{-1}(\pi(F_{B_{a_i}}))$, in which case the
condition $(\star)$ is automatically satisfied.

\vskip.1in

Case: $p \in B_{a_i} - F_{B_{a_i}}^*$ is a fixed point.

\vskip.1in

In this case, set
$$\align
C &:= F_{B_{a_i}} \cap U_p \\
D &:= \pi^{-1}(\pi(F_{B_{a_i}} - U_p)) \cap U_p = (F_{B_{a_i}} - U_p)^{\pm} \cap U_p. \\
\endalign$$
Since $C$ and $D$ are disjoint $K^*$-invariant closed subsets of an affine $K^*$-invariant
variety $U_p$, Corollary 1.2 in Kirwan-Fogarty-Mumford [1] implies that there exists a
$K^*$-invariant function $f \in A(U_p)^{K^*}$ such that
$$f \equiv 1 \text{\ on\ }C \hskip.1in \& \hskip.1in f \equiv 0 \text{\ on\ }D.$$
We only have to replace the oiginal $U_p$ and $V_p$ with
$$\align
(U_p)_f &= \{q \in U_p;f(q) \neq 0\} \text{\ and\ }\\
 (V_p)_f &= \{v \in V_p;i_p^*f(v) \neq 0\},\\
\endalign$$
respectively, in order for Luna's locally toric chart to satisfy the condition $(\star)$.

\vskip.1in

We resume the proof for $\varphi_i$ being V-locally toric.

\vskip.1in

Remark that $\{U_p\}_{p \in B_{a_i} - F_{B_{a_i}}^*}$ where we have Luna's locally toric chart
$$X_p \overset{\eta_p}\to{\leftarrow} V_p \overset{i_p}\to{\rightarrow} U_p \subset B_{a_i}$$
satisfying the condition $(\star)$ is an open covering of $B_{a_i}$.  By the condition $(\star)$
the morphism $V_p//K^* \rightarrow B_{a_i}//K^*$ is \'etale, as it is a composite of an \'etale
morphism $V_p//K^* \rightarrow U_p//K^*$ with an inclusion $U_p//K^* \subset B_{a_i}//K^*$.  We
also have the commutative diagram

$$\CD
(B_{a_i})_-/K^* @. \rightarrow \hskip.25in @. B_{a_i}//K^* @. \hskip.25in \leftarrow @.
(B_{a_i})_+/K^*
\\ @AAA \square \hskip.25in @.@AAA \hskip.25in \square @. @AAA \\
(U_p)_-/K^* \times_{U_p//K^*}V_p//K^* @. \rightarrow \hskip.25in @. V_p//K^* @.
\hskip.25in \leftarrow @. (U_p)_+/K^*
\times_{U_p//K^*}V_p//K^*
\\ @| \hskip.25in @. @| \hskip.25in @. @| \\
(V_p)_-/K^* @. \rightarrow \hskip.25in @. V_p//K^* @. \hskip.25in\leftarrow @. (V_p)_+/K^* \\
@| \hskip.25in @.@| \hskip.25in @. @| \\
(X_p)_-/K^* \times_{X_p//K^*}V_p//K^* @. \rightarrow \hskip.25in @. V_p//K^* @.
\hskip.25in \leftarrow @. (X_p)_+/K^*
\times_{X_p//K^*}V_p//K^*
\\ @VVV \square \hskip.25in @. @VVV \hskip.25in \square @. @VVV \\
(X_p)_-/K^* @. \rightarrow \hskip.25in @. X_p//K^* @. \hskip.25in \leftarrow @. (X_p)_+/K^*. \\
\endCD$$
Therefore, according to Definition 1-2-10 (i), the proper birational map $\varphi_i$ is
V-locally toric.

\vskip.1in

This completes the proof of Theorem 2-4-1.

\newpage

$$\bold{CHAPTER\ 3.\ TRIFICATION}$$

\vskip.2in

In the previous chapter, we observed that the $K^*$-action on the entire
cobordism $B$ and hence on the quasi-elementary cobordisms $B_{a_i}$ is
locally toric and that this gives rise to the $V$-locally toric birational
map

$$\CD
W_i @.@. \overset{\varphi_i}\to{\dashrightarrow} @.@. W_{i+1} \\
@|@.\hskip.6in @.@. @| \\
(B_{a_i})_-/K^* @.@. \hskip.6in @.@.(B_{a_i})_+/K^* \\
@.\searrow \hskip.2in@. \hskip.6in @.\hskip.2in \swarrow @.\\
@.@.(B_{a_i})/K^* @.@.\\
\endCD$$

The reason why $\varphi_i:W_i = (B_{a_i})_-/K^* \dashrightarrow
W_{i+1} = (B_{a_i})_+/K^*$ is merely locally toric and not toroidal is that
the choices of the local coordinate divisors (for Luna's
locally toric charts for
$B_{a_i}$) are not canonical and hence that they do NOT patch together to
provide a global boundary divisor $B_{a_i} - U_{B_{a_i}}$ necessary for
a desired toroidal structure $(U_{B_{a_i}},B_{a_i})$.

In order to create such a global boundary divisor, we blow up $B_{a_i}$
along an ideal $I$, called the ``$\bold{torific}$" ideal on $B_{a_i}$,
canonically defined in terms of the locally toric $K^*$-action on
$B_{a_i}$.  We show that the divisor $D^{tor}$ defined by the principal ideal
$\mu^{-1}(I) \cdot B_{a_i}^{tor}$, where
$$\mu:B_{a_i}^{tor} \rightarrow B_{a_i}$$
is the normalization of the blowup of $B_{a_i}$ along $I$, provides the
toroidal structure $(U_{B_{a_i}^{tor}} := B_{a_i}^{tor} -
D^{tor},B_{a_i}^{tor})$ with respect to which the induced $K^*$-action on
$B_{a_i}^{tor}$ is toroidal.  We call such a procedure
``$\bold{torification}$".  

The basic idea that if one blows up an ideal
then the divisors defined by the principalized ideal provide the resulting
variety with some useful extra structure, can be traced back at least to
Hironaka.  The torific ideal we construct is closely related to the one
used by Abramovich-DeJong [1] for their proof of Weak Resolution of
Singularities.

After torification, $B_{a_i}^{tor}$ turns out to be a quasi-elementary
cobordism and hence we obtain a $V$-toroidal birational map

$$\CD
(U_{W_{i-}^{tor}},W_{i-}^{tor}) @.@.
\overset{\varphi_i^{tor}}\to{\dashrightarrow} @.@.
(U_{W_{i+}^{tor}},W_{i+}^{tor})
\\ @|@.\hskip.6in @.@. @| \\
(U_{B_{a_i}^{tor}},B_{a_i}^{tor})_-/K^* @.@. \hskip.6in
@.@.(U_{B_{a_i}^{tor}},B_{a_i}^{tor})_+/K^* \\ @.\searrow \hskip.2in@.
\hskip.6in @.\hskip.2in \swarrow @.\\ @.@.(U_{B_{a_i}^{tor}},B_{a_i}^{tor})/K^* @.@.\\
\endCD$$

Remark in general that $W_i, W_{i\pm}^{tor}$ are all singular and that the
induced morphisms $f_{i\pm}:W_{i\pm}^{tor} \rightarrow W_i$ are not
sequences of blowups with smooth centers on nonsingular varieties.  Thus it
seems that through the birational cobordism we decomposed $\phi$ into a
sequence of locally toric birational maps, which were made then toroidal via
torification, but at the fatal cost of introducing singularities.  In Chapter
4, we remedy this situation by applying canonical resolution of
singularities to
$W_i, W_{i\pm}^{tor}$ (and canonical
principalization of some ideals on them) so that we can apply the
strong factorization theorem of toroidal birational maps between
nonsingular toroidal embeddings by Morelli [1] (cf.
W{\l}odarczyk [1] Abramovich-Matsuki-Rashid [1]) to finish the proof of the
Weak Factroization Theorem.

\vskip.1in

The ideal situation would have been that we could define the torific ideal
on the entire cobordism $B$, the blowup of the torific ideal would then torify
the
$K^*$-action and hence we could reduce the factorization of $\phi$ to that
of ONE toroidal birational map.  But at the moment, for a small but
essential reason to be discussed below, we can only define the torific
ideal on each quasi-elementary cobordism $B_{a_i}$ separately but not on the
entire cobordism.  Thus we are forced to break up $\phi$ into SEVEARL toroidal
birational maps and cannot make the toroidal structures compatible with each
other.  This is the main reason why we can only prove the WEAK Factorization
Theorem but not the STRONG one (for the moment).

\vskip.1in

In this chapter, $B = B_{a_i}$ always refers
to a nonsingular quasi-elementary cobordism by abuse of notation with $K^*$
acting effectively.  We also omit the subscript $a_i$ from the corresponding
notations, e.g., $F = F_{a_i}, F^* = F_{a_i}^*$ etc.  This convention is
restricted to this chapter only.

\vskip.2in

\S 3-1.  Construction of the torific ideal

\vskip.1in

\proclaim{Definition 3-1-1} Let $p \in B$ be a closed point with $G_p
\subset K^*$ being the stabilizer of $p$.  Fix an integer $\alpha \in {\Bbb
Z}$.  Then we define
$$J_{\alpha,p} \subset {\Cal O}_{B,p}$$
to be the ideal generated by all the semi-invariant functions $f \in {\Cal
O}_{B,p}$ of
$G_p$-character $\alpha$, that is, 
$$t^*(f) = t^{\alpha}\cdot f \text{\ for\ }t \in G_p \subset K^*.$$
\endproclaim

\proclaim{Theorem 3-1-2} There exists a unique coherent $K^*$-equivariant
ideal sheaf \linebreak
$I_{\alpha} \subset {\Cal O}_B$ such that
$$(I_{\alpha})_p = J_{\alpha,p} \hskip.1in \forall p \in B - F^*.$$
(We note that a $K^*$-equivariant sheaf is nothing but a $K^*$-sheaf in the
sense of Mumford [1], where he considers the notion of a $G$-sheaf when a group $G$ acts on a
variety.  In our case,
${\Cal O}_B$ has the natural structure of a
$K^*$-sheaf induced by the $K^*$-action on $B$ and we require that the ideal
$I_{\alpha}$ has the $K^*$-subsheaf structure.)
\endproclaim

\proclaim{Definition 3-1-3} The coherent sheaf $I_{\alpha}$ as above is
called the
$\alpha$-torific ideal sheaf with respect to the $K^*$-action on $B$.
\endproclaim

\proclaim{Remark 3-1-4}\endproclaim

\vskip.1in

(i) It is absolutely necessary to exclude the points in $F^*$ from the ones
on which we require the condition $(I_{\alpha})_p = J_{\alpha,p}$.  The
collection of ideals $J_{\alpha,p}$ for $p \in B$ does not define a
coherent sheaf in general.  As an example, let $B = {\Bbb A}^2$ and let
$K^*$ act on it by
$$t \cdot (x,y) = (tx,t^{-1}y) \text{\ for\ }t \in K^*.$$ 
We fix $\alpha = 1$.  Then at $p = (0,0)$ the stabilizer is $G_p
= K^*$ and
$J_{1,p} = (x)$.  On the other hand, at any other point $q \in B - \{p\}$,
the stabilizer is trivial $G_q = \{1\} \subset K^*$ and hence $J_{1,q} =
{\Cal O}_{B,q}$. While $J_{1,q} = {\Cal O}_{B,q} \text{\ for\ }q = (x,y)
\text{\ with\ }x
\neq 0$ and $J_{1,p} = (x)$ patch together, $J_{1,q} = {\Cal O}_{B,q} \text{\
for\ }q = (x,y)
\text{\ with\ }x = 0$ and $J_{1,p} = (x)$ do NOT.

(ii) For a point $q \in F^*$, the stalk $(I_{\alpha})_q$ is determined by
the $K^*$-equivariance of $I_{\alpha}$ and by the fact that the closure
$\overline{O(q)}$ of the orbit of $q$ has a unique fixed point $p \in F
\subset B - F^*$, since
$B$ is quasi-elementary.  (See the proof for the uniqueness of the
$\alpha$-torific ideal.)  On the other hand, if $B$ is not
quasi-elementary, the closure $\overline{O(q)}$ may have two distinct fixed
points $p_1$ and $p_2$ and the possible stalks at $q$ induced from those at
$p_1$ and $p_2$, respectively, by the $K^*$-equivariance may not coincide. 
This is why we have trouble defining the torific ideal if $B$ is not
quasi-elementary.

\vskip.1in

\demo{Proof of Theorem 3-1-2}\enddemo

\vskip.1in

($\bold{Uniqueness}$) Suppose we have two coherent ideal sheaves
$I_{\alpha}$ and
$I_{\alpha}'$ satisfying the condition of Theorem 3-1-2.  Then
$$(I_{\alpha})_p = (I_{\alpha}')_p = J_{\alpha,p} \text{\ for\ }p \in B -
F^*.$$
We have to check
$$(I_{\alpha})_q = (I_{\alpha}')_q \text{\ for\ }q \in F^*.$$
Since $B$ is quasi-elementary, the closure $\overline{O(q)}$ contains only
one fixed point $r \in \overline{O(q)}$.  Since $r \in F \subset B -
F^*$, we have
$$(I_{\alpha})_r = (I_{\alpha}')_r = J_{\alpha,r}.$$
Since both $I_{\alpha}$ and $I_{\alpha}'$ are coherent, there exists an open
neighborhood $r \in U$ such that
$$I_{\alpha}|_U = I_{\alpha}'|_U.$$
Take a point $q_o = t_o \cdot q \in O(q) \cap U$.  Then by the
$K^*$-equivariance we have
$$(I_{\alpha})_q = t_o^*\{(I_{\alpha})_{q_o}\} =
t_o^*\{(I_{\alpha}')_{q_o}\} = (I_{\alpha}')_q,$$
proving the uniqueness.

\vskip.1in

($\bold{Existence}$) 

\vskip.1in

Step 1. First observe that it is enough to prove the
following.  

\proclaim{Claim 3-1-5} For each point $p \in B - F^*$, there exist a
$K^*$-invariant open neighborhood $p \in U_p$ satisfying the condition $(*)$:
$$(*) \text{\ If\ a\ }\text{$K^*$-orbit\ } O(q) \text{\ lies\ in\ }U_p,
\text{\ then\ its\ closure\ }\overline{O(q)} \text{\ also\ lies\ in\ }U_p.$$
and the $\alpha$-torific ideal over $U_p$, i.e., a $K^*$-equivariant coherent
ideal sheaf
$I_{\alpha, U_p}$ over
$U_p$ such that
$$(I_{\alpha,U_p})_q = J_{\alpha,q} \hskip.1in \forall q \in U_p \cap (B -
F^*).$$
\endproclaim

In fact, the collection $\{U_p\}_{p \in B - F^*}$ is an open
covering of $B$ and the collection of ideals $\{I_{\alpha,U_p}\}_{p \in B -
F^*}$ patches together by the uniqueness argument, using the condition
$(*)$, to provide the global $\alpha$-torific ideal over $B$ as required. 

\vskip.1in

Step 2. We take Luna's locally toric chart for $p \in B - F^*$
$$X \overset{\eta}\to{\leftarrow} V \overset{i}\to{\rightarrow} Y = U_p
\subset B$$
where $\eta$ and $i$ are strongly \'etale with $i$ being surjective, $X_p$ is a nonsingular
affine toric variety and where $Y = U_p$ satisfies the condition $(*)$. 

\proclaim{Claim 3-1-6} In order to show the existence of the $\alpha$-torific
ideal
$I_{\alpha,U_p}$, it suffices to show the existence of the $\alpha$-torific
ideal $I_{\alpha,X}$ over the toric variety $X$.
\endproclaim

We spend the rest of Step
2 for the verification of this claim.

\proclaim{Lemma 3-1-7} Let $p \in B$ and $J_{\alpha,p} \subset {\Cal
O}_{B,p}$ be the ideal generated by all the semi-invarinat functions in
${\Cal O}_{B,p}$ of
$G_p$-character $\alpha$ as in Definition 3-1-1.  We lift the $G_p$-action
on ${\Cal O}_{B,p}$ to the one on the completion ${\hat {\Cal O}_{B,p}}$. 
Then the completion $\widehat{J_{\alpha,p}}$ of the ideal $J_{\alpha,p}$
$$\widehat{J_{\alpha,p}} = J_{\alpha,p} \otimes \widehat{{\Cal O}_{B,p}}
\subset {\Cal O}_{B,p} \otimes \widehat{{\Cal O}_{B,p}} = \widehat{{\Cal
O}_{B,p}}$$ is characterized as the ideal generated by all the
semi-invariant functions in $\widehat{{\Cal O}_{B,p}}$ of $G_p$-character
$\alpha$ with respect to the induced $G_p$-action on $\widehat{{\Cal O}_{B,p}}$.
\endproclaim

\demo{Proof}\enddemo The affine coordinate ring for $Y = \roman{Spec} A$
splits into the direct sum
$$A = \oplus_{\beta \in {\Bbb Z}} A_{\beta}$$
where the $A_{\beta}$ consist of all the semi-invariant functions of
$K^*$-character
$\beta$.  This splitting induces the splitting of the maximal ideal $m_p$ for the
point $p$
$$m_p = \oplus_{\beta' \in {\Bbb Z}/b_p} m_{p,\beta'}$$
where $m_{p,\beta'}$ are the elements of $m_p$ with $G_p$-character
$\beta'$.  (The stabilizer $G_p \subset K^*$ is identified with the group of
$b_p$-th roots of unity.  We set $b_p = 0$ when $G_p = K^*$.)  In fact,
it is easy to see that any direct summand $A_{\beta}$ must be contained in
$m_p$ if $\beta (\roman{mod}\ b_p) \neq 0$.  Thus taking $m_{p,0} = m_p
\cap \oplus_{\beta (\roman{mod}\ b_p) = 0}A_{\beta}$ and $m_{p,\beta'} =
\oplus_{\beta (\roman{mod}\ b_p) = \beta'}A_{\beta}$ for $\beta' \neq
0$ we obtain the above splitting.  Therefore, we
can choose semi-invariant functions $z_1, \cdot\cdot\cdot, z_n$ (with
respect to the $G_p$-action) which form a basis of the vector space $m_p/m_p^2$
and hence generates $m_p$ by Nakayama's Lemma.  We have a $G_p$-equivariant
isomorphism
$$\widehat{{\Cal O}_{B,p}} \cong K[[z_1, \cdot\cdot\cdot, z_n]]$$
where the action of $G_p$ on the formal power series is induced by the
$G_p$-characters on $z_1, \cdot\cdot\cdot, z_n$.

Since the $G_p$-action respects the natural grading on $K[[z_1,
\cdot\cdot\cdot, z_n]]$ given by the degree and since $K[[z_1,
\cdot\cdot\cdot, z_n]]$ is noetherian, the
ideal $K_{\alpha}$ generated by all the semi-invariant
functions in
$\widehat{{\Cal O}_{B,p}}$ of $G_p$-character $\alpha$ is generated by a
finite number of monomials $m_1, \cdot\cdot\cdot, m_l$ of $G_p$-character
$\alpha$ in $z_1, \cdot\cdot\cdot, z_n$.  Since obviously $m_1,
\cdot\cdot\cdot, m_l \in J_{\alpha,p}$, we have
$$K_{\alpha} = (m_1, \cdot\cdot\cdot, m_l) \subset
\widehat{J_{\alpha,p}}.$$ 
Since $J_{\alpha,p} \subset K_{\alpha}$ by definition and hence
$\widehat{J_{\alpha,p}}
\subset K_{\alpha}$, we conclude
$$K_{\alpha} = \widehat{J_{\alpha,p}}$$
as required.

\vskip.1in

We resume the verification of Claim 3-1-6.

\vskip.1in

For any point $q \in Y$, let $q_V \in i^{-1}(q)$ be a point in the inverse
image of $q$ by $i$ and $q_X = \eta(q_V)$.  (We note that $q_V$ and hence
$q_X$ may not be determined uniquely only by specifying the point $q$.) 
Observe by the condition $(*)$ that we have
$$q \in B - F^* \Longleftrightarrow q_X \in X - F_X^*$$
and in particular
$$p \in B - F^* \Longrightarrow p_X \in X - F_X^*.$$

Let ${\Cal H}_{\alpha} \subset {\Cal O}_Y$ be the ideal sheaf over $Y$
generated by the monomials $m_1, \cdot\cdot\cdot, m_l$ as taken in the
proof of Lemma 3-1-7.
Then we have
$$\widehat{(\eta^*I_{\alpha,X})_{p_V}} =
\eta^*(\widehat{(I_{\alpha,X})_{p_X}})
\overset{\overset{\hskip.05in\eta^*}\to{\sim}}\to{\longleftarrow}
\widehat{(I_{\alpha,X})_{p_X}}$$
and
$$\widehat{(i^*{\Cal H}_{\alpha})_{p_V}} =
i^*(\widehat{({\Cal H}_{\alpha})_p})
\overset{\overset{\hskip.05in i^*}\to{\sim}}\to{\longleftarrow}
\widehat{({\Cal H}_{\alpha})_p} = \widehat{J_{\alpha,p}}.$$
Now since $I_{\alpha,X}$ is the $\alpha$-torific ideal over $X$ and since
$p_X \in X - F_X^*$, Lemma 3-1-7 characterizes both
$\widehat{(\eta^*I_{\alpha,X})_{p_V}}$ and $\widehat{(i^*{\Cal H})_{p_V}}$
as the ideal generated by all the semi-invariant functions in
$\widehat{{\Cal O}_{V,p_V}}$ of $G_{p_V}$-character $\alpha$.  (Remark
that $G_{p_V} = G_{p_X} = G_p$ since $\eta$ and $i$ are both strongly
\'etale.)   Therefore, we
have
$$\widehat{(\eta^*I_{\alpha,X})_{p_V}} = \widehat{(i^*{\Cal
H}_{\alpha})_{p_V}},$$  which implies (See, e.g., Matsumura [1] Theorem 8.15
(3).)
$$(\eta^*I_{\alpha,X})_{p_V} = (i^*{\Cal H}_{\alpha})_{p_V}.$$
Therefore, noting that the equality above holds for all the points $p_V
\in i^{-1}(p)$ and by shrinking
$U_p = Y$ and accordingly
$V$ if necessary, we may assume
$$\eta^*I_{\alpha,X} = i^*{\Cal H}_{\alpha}.$$
Now for any point $q \in Y - F^*$ we have $q_X \in X - F_X^*$ and
$$i^*(\widehat{J_{\alpha,q}}) = \eta^*(\widehat{(I_{\alpha,X})_{q_X}}),$$
both sides being characterized as the ideal generated by all the
semi-invariant functions of $G_{q_V} (= G_q = G_{q_X})$-character $\alpha$
in $\widehat{{\Cal O}_{V,q_V}} (\cong \widehat{{\Cal O}_{B,q}} \cong
\widehat{{\Cal O}_{X,q_X}})$.  On the other hand, we have
$$\align
\eta^*(\widehat{(I_{\alpha,X})_{q_X}}) &=
\widehat{(\eta^*I_{\alpha,X})_{q_V}} \\
&= \widehat{(i^*{\Cal H}_{\alpha})_{q_V}} \\
&= i^*(\widehat{({\Cal H}_{\alpha})_q}).\\
\endalign$$
Since $i$ is \'etale, we conclude
$$\widehat{J_{\alpha,q}} = \widehat{({\Cal H}_{\alpha})_q},$$
which implies
$$J_{\alpha,q} = ({\Cal H}_{\alpha})_q.$$
Since $\eta$ and $i$ are $K^*$-equivariant and since $I_{\alpha,X}$ is
$K^*$-equivariant as it is the $\alpha$-torific ideal over $X$, the
equality $\eta^*I_{\alpha,X} = i^*{\Cal H}_{\alpha}$ implies that ${\Cal
H}_{\alpha}$ is also $K^*$-equivariant.  

Therefore, we conclude that we
can take ${\Cal H}_{\alpha}$ to be the $\alpha$-torific ideal
$I_{\alpha,U_p}$ over $Y = U_p$.

\vskip.1in

Step 3. We show the existence of the $\alpha$-torific ideal on an affine
nonsingular toric variety $X$ where $K^*$-acts as a one-parameter subgroup.

\proclaim{Proposition 3-1-8} Let $X = X(N,\sigma)$ be an affine
nonsingular toric variety where $K^*$ acts as a one-parameter subgroup.  Then the
ideal generated by all the monomials of $K^*$-character $\alpha$ is the
$\alpha$-torific ideal over
$X$.
\endproclaim

\demo{Proof}\enddemo We take a ${\Bbb Z}$-basis $\{v_1, \cdot\cdot\cdot,
v_n\}$ of the lattice $N$ so that
$$\align
\sigma &= \langle v_1, \cdot\cdot\cdot, v_l\rangle \\
\sigma^{\vee} &= \langle v_1^*, \cdot\cdot\cdot v_l^*, \pm v_{l+1}^*,
\cdot\cdot\cdot, \pm v_n^*\rangle.\\
\endalign$$
Setting $z_j = z^{v_j^*}$, we have
$$X = \roman{Spec} K[z_1, \cdot\cdot\cdot, z_l,z_{l+1}^{\pm 1},
\cdot\cdot\cdot, z_n^{\pm 1}].$$
The $K^*$-action corresponds to a one-parameter subgroup
$$a = (\alpha_1, \cdot\cdot\cdot, \alpha_n) \in N,$$
i.e., the action of $t \in K^*$ is given by
$$t \cdot (z_1, \cdot\cdot\cdot, z_n) = (t^{\alpha_1}z_1, \cdot\cdot\cdot,
t^{\alpha_n}z_n).$$

Let $H_{\alpha,X}$ be the ideal generated by all the monomials of
$K^*$-character $\alpha$.  Then we see obviously that $H_{\alpha,X}$ is
$K^*$-equivariant and that
$$(H_{\alpha,X})_q \subset J_{\alpha,q} \text{\ for\ }q \in X - F_X^*.$$

Thus in order to prove $H_{\alpha,X}$ is the $\alpha$-torific ideal over $X$ we
only have to show the opposite inclusion
$$J_{\alpha,q} \subset (H_{\alpha,X})_q.$$

Remark first that $J_{\alpha,q}$ is generated by all the monomials in $z_1,
\cdot\cdot\cdot, z_n$ \linebreak of $G_q$-character $\alpha$.  (Note that
$z_1 - z_1(q), \cdot\cdot, z_n(q) - z_n(q)$ are generators of the maximal
ideal $m_q$ and are all semi-invariant functions with respect to the action of
$G_q$.  Indeed, if $z_i(q) \neq 0$, then $z_i$ has the trivial
$G_q$-character.  Then using the argument in the proof of Lemma 3-1-7, we see
that $J_{\alpha,q}$ is generated by all the monomials in $z_1 - z_1(q),
\cdot\cdot\cdot, z_n(q) - z_n(q)$ of $G_q$-character $\alpha$ and hence generated
by all the monomials in $z_1, \cdot\cdot\cdot, z_n$ of $G_q$-character
$\alpha$). 

Let $f \in {\Cal O}_{X,q}$ be a monomial in $z_1,
\cdot\cdot\cdot, z_n$ of $G_q$-character $\alpha$.  We claim that there exists a
monomial $g$ in $z_1,
\cdot\cdot\cdot, z_n$ invertible at $q$, i.e., $g \in {\Cal O}_{X,q}^*$, such
that $f \cdot g$ is a monomial of $K^*$-character $\alpha$.  

In fact, let $b_q
\in {\Bbb Z}_{\geq 0}$ be the non-negative generator of the subgroup of ${\Bbb
Z}$ generated by the following set
$$A_q := \{\alpha_j;\alpha_j \neq 0 \text{\ with\ }z_j(q) \neq 0\}.$$
The stabilizer $G_q \subset K^*$ consists of the $b_q$-th roots of unity.  Let
$\alpha_f$ be the
$K^*$-character of the monomial
$f$.  Then
$\alpha -
\alpha_f$ is a multiple of $b_q$, say $b_fb_q$.  By the description of the set
$A_q$, we have a monomial $g$ in $z_1, \cdot\cdot\cdot, z_n$, invertible at $q$,
of
$K^*$-character
$b_fb_q$.  Then $f \cdot g$ is a monomial of $K^*$-character $\alpha$ by
construction.

Since $f \cdot g$ is a monomial regular at $q$, it is regular in a
$K^*$-invariant open neighborhood $U_q$ of $q$.  Since $X$ is affine and hence
quasi-elementary, we conclude that $U_q$ must contain $\pi^{-1}(U)$, where
$\pi:X \rightarrow X//K^*$ is the quotient morphism and where $U$ is some open
neighborhood of $\pi(q)$.  Take a global function $h$ on $X//K$ (which is
identified with a global
$K^*$-invariant function on $X$) such that $h \equiv 0$ on $X//K^* - U$ and
$h(\pi(q)) = 1$  (See Mumford-Fogarty-Kirwan [1] Corollary 1.2.).  Then $f \cdot
g
\cdot h^l$ (for
$l >> 0$) is a global regular function of $K^*$-character $\alpha$, hence a
linear combination of monomials in $z_1, \cdot\cdot\cdot, z_n$ of $K^*$-character
$\alpha$.

Since $f$ is an arbitrary monomial of $G_q$-character $\alpha$ and since $g$
and $h$ are invertible at $q$, we conclude that $J_{\alpha,q}$ generated by all
the monomials of $G_q$-character $\alpha$ is also generated by all the
monomials of $K^*$-character $\alpha$.  Therefore, we have
$$J_{\alpha,q} \subset (H_{\alpha,X})_q.$$
This completes the proof for Proposition 3-1-8 and hence the proof for Theorem
3-1-2 via Steps 1, 2 and 3.

\vskip.2in

\S 3-2. The torifying property of the torific ideal  

\vskip.1in

In this section, we discuss the torifying property of the torific ideal (which
is actually the product of several $\alpha$-torific ideals defined in the
previous section, where the characters
$\alpha$ are taken from a certain finite set), i.e., how we can induce the
torification of the $K^*$-action by blowing up the torific ideal, as announced
at the beginning of Chapter 3.

\vskip.1in

First we explain the $\bold{torification\ principle}$ for affine nonsingular
toric varieties $X$, showing the torifying property for
them.  The explanation of the torification principle should also provide the
motivation for the definition of the torific ideal, including that of the
$\alpha$-torific ideal of the previous section.  Secondly,
utilizing on Luna's locally toric charts and based upon the torifying property for toric
varieties, we discuss the torifying property of the torific ideal on 
a quasi-elementary cobordism in general.

\proclaim{Observation 3-2-1 (Torification Principle for Toric
Varieties)}\endproclaim

Let
$X = X(N,\sigma)$ be an affine nonsingular toric variety with $K^*$ acting as a
one-parameter subgroup.  As in Proposition 3-1-8, we take a ${\Bbb Z}$-basis of
$N$ so that
$$\align
\sigma &= \langle v_1, \cdot\cdot\cdot, v_l\rangle \\
\sigma^{\vee} &= \langle v_1^*, \cdot\cdot\cdot, v_l^*, \pm v_{l+1}^*,
\cdot\cdot\cdot, \pm v_n^*\rangle \\
X &= \roman{Spec}K[z_1, \cdot\cdot\cdot, z_l,z_{l+1}^{\pm}, \cdot\cdot\cdot,
z_n^{\pm}] \text{\ with\ }z_j = z^{v_j^*} \\
\endalign$$
and that the $K^*$-action corresponds to
$$a = (\alpha_1, \cdot\cdot\cdot, \alpha_n) \in N.$$
We have the coordinate divisors
$$D_{v_1} = \{z_1 = 0\}, \cdot\cdot\cdot, D_{v_l} = \{z_l = 0\}$$
corresponding to the 1-dimensional faces (generated by the $v_i$) of $\sigma$. 
As
$X$ is a toric variety containing the torus $T_X$, obviously
$$(T_X = X - (D_{v_1} \cup \cdot\cdot\cdot \cup D_{v_l}),X)$$
has the toric and hence toroidal structure with respect to which the
$K^*$-action is toric and hence toroidal.

Let $I$ be an ideal generated by some monomials and let
$$\mu_X:{\tilde X} \rightarrow X$$
be the normalization of the blowup of $X$ along $I$.  Since $I$ is
toric having monomial generators, the torus-action lifts to the blowup of $X$
along $I$ and to the normalization ${\tilde X}$, which is thus again a toric
variety ${\tilde X} = X(N,{\tilde \sigma})$ where ${\tilde \sigma}$ is a fan
obtained by subdividing $\sigma$.  (We define a toric variety to be a NORMAL
variety which contains a torus, as a dense open subset, whose action on itself
extends to the entire variety.)  The
$K^*$-action also lifts to ${\tilde X}$ being a one-parameter subgroup.  

Let ${\tilde D}$
to be the divisor defined by the principal ideal $\mu_X^{-1}(I) \cdot {\Cal
O}_{\tilde X}$.  Then as ${\tilde X}$ is a toric variety containing the torus
$T_X = T_{\tilde X}$, 
$$(T_{\tilde X} = {\tilde X} - \{(D_{v_1} \cup \cdot\cdot\cdot \cup D_{v_l})
\cup {\tilde D}\},{\tilde X})$$
has the toric and hence toroidal structure with respect to which the $K^*$-action
is toric and hence toroidal.  (We use the same notation $D_{v_i}$ for the
strict transform on ${\tilde X}$ of the divisor $D_{v_i}$ on $X$, indicating
both are the closures of the orbit corresponding to the 1-dimensional ray
generated by $v_i$ in ${\tilde X}$ and $X$, respectively.)

We look for a situation on ${\tilde X}$ (and the conditions on $I$ which yields
that situation on ${\tilde X}$) where
$$(U_{\tilde X} = {\tilde X} - {\tilde D},{\tilde X})$$
has a toroidal structure with respect to which the $K^*$-action is toroidal. 
That is to say, after removing those coordinate divisors among $D_{v_1},
\cdot\cdot\cdot, D_{v_l}$ which are NOT contained in ${\tilde D}$ from the
toric boundary divisor $(D_{v_1} \cup \cdot\cdot\cdot \cup D_{v_l}) \cup
{\tilde D}$ we still want to preserve a toroidal stucture on $(U_{\tilde
X},{\tilde X})$ with respect to which the $K^*$-action remains toroidal. 

Recall that the original coordinate divisors $D_{v_1}, \cdot\cdot\cdot,
D_{v_l}$ correspond to a locally toric chart of the quasi-elementary cobordism
$B$ and are not chosen canonically and that we have to remove those which are
NOT contained in the canonically constructed boundary divisor, which
corresponds to ${\tilde D}$ which arises from the torific
ideal defined canonicaly only in terms of the $K^*$-action.

\vskip.1in

\proclaim{Lemma 3-2-2} We consider a fixed coordinate divisor $D_{v_j}$ where $j
= 1,
\cdot\cdot\cdot, l$.  Suppose that we have the following situation
$(\heartsuit_j)$:

$$(\heartsuit_j) \left\{\aligned
&\text{For\ any\ affine\ (open)\ toric\ subvariety\ } Z \subset
{\tilde X},\\  
&\text{we\ have\ } (T_Z,Z) = (T_{{\Bbb A}^1},{\Bbb
A}^1) \times (T_{{Z'}_j},{Z'}_j) \text{as\ toric\ varieties,}\\
&\text{\ provided\ that\ }D_{v_j} \not\subset
{\tilde D} \text{\ and\ }Z \cap D_{v_j} \neq \emptyset,\\  
& \hskip.3in
\text{where\ }{\Bbb A}^1 \text{\ has\ the\ standard\ toric\ structure\ with\
torus\ }T_{{\Bbb A}^1} = {\Bbb A}^1 - \{0\}\\ & \hskip.5in \text{\ and\ }{Z'}_j
\text{\ is\ an\ affine\ toric\ variety\ with\ the\ torus\ }T_{{Z'}_i}\\ 
& \hskip.3in \text{and\ where\
the\ $K^*$-action\ on\ }{\Bbb A}^1 =
\roman{Spec}K[x_j]\text{\ is\ trivial\ } \\
&\hskip.5in \text{\ and\ $K^*$\ acts\ on\ }{Z'}_j \text{\ as\ a\ one-parameter\
subgroup\ in\ }T_{{Z'}_j}\\
&\hskip.5in\text{\ and\ we\ have\ }Z \cap D_{v_j} = \{x_j = 0\}.\\
\endaligned
\right.$$
Then
$$(U_{{\tilde X},j} = {\tilde X} - \{(D_{v_1} \cup \cdot\cdot\cdot
\cup \overset{\vee}\to{D_{v_j}} \cup \cdot\cdot\cdot \cup D_{v_l}) \cup {\tilde
D}\},{\tilde X})$$
has a toroidal structure with respect to which the $K^*$-action is
toroidal.
\endproclaim

\demo{Proof}\enddemo If $D_{v_j} \subset {\tilde D}$, then for any affine
(open) toric subvariety
$Z \subset {\tilde X}$ we have
$$(T_Z,Z) = Z \cap (U_{{\tilde X},j},{\tilde X}) = Z \cap (U_{\tilde X},{\tilde X})$$
has the toric and hence toroidal structure with respect to which the
$K^*$-action is toric and hence toroidal.  Since these $Z$'s cover the entire
${\tilde X}$, the claim follows.

So we may assume $D_{v_j} \not\subset {\tilde D}$.  Let $Z \subset {\tilde X}$
be an affine (open) toric subvariety.  

If $Z \cap D_{v_j} = \emptyset$, then
$$(T_Z,Z) = Z \cap (U_{{\tilde X},j},{\tilde X}) = Z \cap (U_{\tilde X},{\tilde X})$$
has the toric and hence toroidal structure with respect to which the
$K^*$-action is toric and hence toroidal.  
 
If $Z \cap D_{v_j} \neq \emptyset$,
then by the situation $(\heartsuit_j)$ we have
$$\align
Z \cap (U_{{\tilde X},j},{\tilde X}) &= ({\Bbb A}^1,{\Bbb A}^1) \times (U_{{Z'}_j},{Z'}_j) \\
&= (U_{Z_{\{0\}}},Z_{\{0\}}) \cup (U_{Z_{\{1\}}},Z_{\{1\}})\\
& \text{where}\\
&(U_{Z_{\{0\}}},Z_{\{0\}}) := ({\Bbb A}^1 - \{0\},{\Bbb A}^1 - \{0\}) \times
(U_{{Z'}_j},{Z'}_j)
\text{\ and}\\ 
&(U_{Z_{\{1\}}},Z_{\{1\}}) := ({\Bbb A}^1 -
\{1\},{\Bbb A}^1 - \{1\}) \times (U_{{Z'}_j},{Z'}_j).\\
\endalign$$
Since the $K^*$-action on ${\Bbb A}^1 =({\Bbb A}^1 - \{0\},{\Bbb
A}^1 - \{0\}) \cup ({\Bbb A}^1 - \{1\},{\Bbb A}^1 - \{1\})$ is trivial, we have
$$(U_{Z_{\{0\}}},Z_{\{0\}}), (U_{Z_{\{1\}}},Z_{\{1\}})
\overset{\text{$K^*$-equivariant}}\to{\cong} ({\Bbb A}^1 -
\{0\},{\Bbb A}^1 -
\{0\})
\times (U_{{Z'}_j},{Z'}_j).$$
Since $Z$'s with $Z \cap D_{v_j} = \emptyset$ and $Z_{\{0\}},Z_{\{1\}}$'s with
$Z \cap D_{v_j} \neq \emptyset$ cover the entire ${\tilde X}$, the claim
follows.

\vskip.1in

\proclaim{Remark 3-2-3} Let $S = \{j;D_{v_j} \in {\tilde D}\}$ be the set of indices of the
primitive vectors generating the 1-dimensional cones in
${\tilde \sigma}$ which correspond to the divisors in ${\tilde D}$.  Then
$(\heartsuit_j)$ is equivalent to the following $(\heartsuit_{j,fan})$ in terms
of the geometry of the fan ${\tilde \sigma}$:

$$(\heartsuit_{j,fan}) \left\{\aligned
&\text{For\ any\ cone\ } \sigma_Z \in {\tilde \sigma},\\  
&\text{we\ have\ } \sigma_Z = \langle v_j\rangle \oplus \tau_{{Z'}_j} \text{\ as\ cones},\\
&\text{provided\ that\ }j \not\in S \text{\ and\ }v_j \in \sigma_Z,\\
&\hskip.3in \text{where\ the\ cone\ }\tau_{{Z'}_j} = \langle w_k's \rangle \text{\ is\
generated\ by}\\
&\hskip.5in \text{the\ extremal\ rays\
$w_k$'s\ (other\ than\
those\ $v_j$'s\ with\ $j \not\in S$)\ of\ }\sigma_Z \\
&\hskip.5in \text{and\ is\ contained\ in\ a\ hyperplane\ }H_j \text{\ which\ also\
contains\ }a \in N\\ 
&\hskip.3in \text{and\ where\ we\ have\ the\ exact\ sequence}\\
&\hskip.5in 0 \rightarrow H_j \cap N \rightarrow N \rightarrow {\Bbb Z}
\rightarrow 0\\
&\hskip.5in \text{with\ the\ image\ of\ }v_j \in N \text{\ to\ be\ }\pm 1 \in
{\Bbb Z}.\\
\endaligned
\right.$$
\endproclaim

Now it is immediate to see that the blowup of the $\alpha_j$-torific ideal
$I_{\alpha_j,X}$ yields the situation $(\heartsuit_j)$.  (See the proof of
Theorem 3-2-6.)  This justifies and actually motivated the definition of the
$\alpha$-torific ideal in the previous section.

In order to see that we can remove all the coordinate divisors which are not
contained in ${\tilde D}$ simultaneously, not just one $D_{v_j} \not\subset
{\tilde D}$ at a time and  
$$(U_{\tilde X},{\tilde X})$$
has a toroidal structure, we need the lemma below.

\proclaim{Lemma 3-2-4} Suppose we have the following situation $(\heartsuit)$:
 $$(\heartsuit) \left\{\aligned
&\text{For\ any\ affine\ (open)\ toric\ subvariety\ } Z \subset
{\tilde X},\\  
&\text{we\ have\ } (T_Z,Z) = \prod_j(T_{{\Bbb A}_j^1},{\Bbb
A}_j^1) \times (T_{Z'},Z') \text{as\ toric\ varieties,}\\
&\text{with\ the\ indices\ }j \text{\ varying\ among\ those\ with\ }D_{v_j}
\not\subset {\tilde D} \text{\ and\ }Z \cap D_{v_j} \neq \emptyset,\\  
& \hskip.3in
\text{where\ }{\Bbb A}_j^1 \text{\ has\ the\ standard\ toric\ structure\ with\
torus\ }T_{{\Bbb A}_j^1} = {\Bbb A}_j^1 - \{0\}\\ & \hskip.5in \text{\ and\ }Z'
\text{\ is\ an\ affine\ toric\ variety\ with\ the\ torus\ }T_{Z'}\\ 
& \hskip.3in \text{and\ where\
the\ $K^*$-action\ on\ }{\Bbb A}_j^1 =
\roman{Spec}K[x_j]\text{\ is\ trivial\ } \\
&\hskip.5in \text{\ and\ $K^*$\ acts\ on\ }Z' \text{\ as\ a\ one-parameter\
subgroup\ in\ }T_{Z'}\\
&\hskip.5in\text{\ and\ we\ have\ }Z \cap D_{v_j} = \{x_j = 0\}.\\
\endaligned
\right.$$
Then
$$(U_{\tilde X} = {\tilde X} - {\tilde D},{\tilde X})$$
has a toroidal structure with respect to which the $K^*$-action is
toroidal.
\endproclaim

The proof is identical to that of Lemma 3-2-2.

\proclaim{Lemma 3-2-5} The collection of the situations $(\heartsuit_j)$ for
all $j = 1, \cdot\cdot\cdot, l$ implies the situation $(\heartsuit)$.
\endproclaim

\demo{Proof}\enddemo We prove that the collection of the situations
$(\heartsuit_{j,fan})$ implies the following situation $(\heartsuit_{fan})$,
which is equivalent to $(\heartsuit)$ but in terms of the geometry of the fan
${\tilde \sigma}$:

$$(\heartsuit_{fan}) \left\{\aligned
&\text{For\ any\ cone\ } \sigma_Z \in {\tilde \sigma},\\  
&\text{we\ have\ } \sigma_Z = \Sigma_j\langle v_j\rangle \oplus \tau_{Z'} \text{\ as\ cones},\\
&\text{with\ the\ indices\ varying\ among\ those\ }j \not\in S \text{\ with\ }v_j
\in \sigma_Z \\
&\hskip.3in \text{where\ the cone\ }\tau_{Z'} = \langle w_k's \rangle \text{\ is generated\
by}\\
&\hskip.5in \text{the\ extremal\ rays\
$w_k$'s\ (other\ than\ those\
$v_j$'s\ with $i \not\in S$)\ of\ }\sigma_Z \\
&\hskip.5in \text{and\ is\ contained\ in\ a\ linear space\ }L\ (\text{of\
codimension\ }\#\{j;j \in S,v_j \in \sigma_Z\}) \\ 
&\hskip.5in \text{which\ also\ conains\ }a \in N\\
&\hskip.3in \text{and\ where\ we\ have\ the\ exact\ sequence}\\
&\hskip.5in 0 \rightarrow L \cap N \rightarrow N \rightarrow {\Bbb
Z}^{\oplus \#\{j;j
\not\in S,v_j \in \sigma_Z\}}
\rightarrow 0\\
&\hskip.5in \text{with\ the\ image\ of\ }\{v_j;j \not\in S, v_j \in \sigma_Z\} \in N \\
&\text{\ forming\ a\ ${\Bbb Z}$-basis\ of\ }{\Bbb Z}^{\oplus \#\{j;j
\not\in S,v_j \in \sigma_Z\}}.\\
\endaligned
\right.$$
Now having the exact sequences for the indices $j \not\in S$ with $v_j \in \sigma_Z$
$$0 \rightarrow H_j \cap N \rightarrow N \rightarrow {\Bbb Z} \rightarrow 0$$
where $v_j \in N$ maps to $\pm 1 \in {\Bbb Z}$, it is elementary to see that we
have the exact sequence with $L = \cap_{\{j;j \not\in S,v_j \in \sigma_Z\}} \cap H_j$
$$0 \rightarrow L \cap N \rightarrow N \rightarrow {\Bbb Z}^{\oplus \#\{j;j
\not\in S,v_j \in \sigma_Z\}} \rightarrow 0$$
such that $\{v_j;j \not\in S, v_j \in \sigma_Z\}$ maps to a ${\Bbb Z}$-basis of ${\Bbb
Z}^{\oplus
\#\{j;j
\not\in S,v_j \in \sigma_Z\}}$.

This completes the proof for Lemma 3-2-5.

\vskip.1in

Therefore, in order to achieve the situation $(\heartsuit)$, we only have to
verify the collection of the situations $(\heartsuit_j)$ for $j = 1,
\cdot\cdot\cdot, l$, which should follow from taking the blowup
$\mu_X$ to be the ``collection" of the blowups of the $I_{\alpha_j,X}$ for $j =
1,
\cdot\cdot\cdot, l$.  This motivates us to define the torific ideal $I$ to be
th product of these:
$$I := \prod_{j = 1}^l I_{\alpha_j,X}.$$  
(Actually in Definition 3-2-7 we may multiply some more $\alpha$-torific
ideals with $\alpha \neq 1, \cdot\cdot\cdot, l$.  See the definition of the
torific ideal for the quasi-elementary cobordism $B$.) 

\vskip.1in

This finishes the discussion of the torification principle for toric varieties.

\vskip.1in

Now we present the verification of the torifying property of the torific ideal
according the torification principle.

\proclaim{Lemma 3-2-6} There exists a finite set ${\goth C}$ of integers so that
at each point $q \in B$ it contains all the characters of $G_q$ acting on the
tangent space $T_{B,q}$ (modulo $b_q$).
\endproclaim

\demo{Proof}\enddemo In the case of $B = X$ being a affine nonsingular toric
variety, the requirement is equivalent to $\{\alpha_1,
\cdot\cdot\cdot,
\alpha_l\} \subset {\goth C}$ in the notation of Observation 3-2-1.  The claim
for the case of a general quasi-elementary cobordism $B$ follows from this and
the fact that $B$ can be covered by a finite number of Luna's locally toric charts.

\proclaim{Definition 3-2-7} Let ${\goth C}$ be a finite set of integers as
above.  We define the torific ideal
$I_{\goth C}$ with respect to ${\goth C}$ to be
$$I = \prod_{\alpha \in {\goth C}} I_{\alpha}.$$
Note that the torific ideal $I_{\goth C}$ does depend on the choice of ${\goth
C}$ but we usually suppress this, call it the torific ideal and denote it by
$I$ without the subscript.
\endproclaim

\proclaim{Theorem 3-2-8} Let $X$ be an affine nonsingular toric variety with
$K^*$ acting as a one-parameter subgroup.  Let
$$\mu_X:X^{tor} \rightarrow X$$
be the normalization of the blowup of the torific ideal $I_X$ for $X$ and
$D^{tor}$ the divisor defined by the principal ideal $\mu_X^{-1}(I_X) \cdot
{\Cal O}_{X^{tor}}$.  Then

\ \ (i) $(U_{X^{tor}} := X^{tor} - D^{tor},X^{tor})$ has a toroidal structure with respect to
which the induced $K^*$-action is toroidal.  More precisely, we have the situation
$(\heartsuit)$ for $X^{tor}$.

\ \ (ii) $X^{tor}$ is a quasi-elementary cobordism.

\ \ (iii) Denoting by $\{\sigma_{Z_{max}}\}$ the set of the maximal cones of
$\sigma^{tor}$, where $\sigma^{tor}$ is the fan corresponding to the toric variety $X^{tor} =
X(N,\sigma^{tor})$.  Then each of the affine toric varieties $\{X(N,\sigma_{Z_{max}}\}$
satisfies the condition $(*)$: 

$(*)$ If a $K^*$-orbit $O(q)$ lies in $X(N,\sigma_{Z_{max}})$, then its closure
$\overline{O(q)}$ also lies in $X(N,\sigma_{Z_{max}})$.

As $\{X(N,\sigma_{Z_{max}}\}$ cover $X^{tor} = X(N,\sigma^{tor})$, the quotients
$\{X(N,\sigma_{Z_{max}})/K^* = X(\pi(N),\pi(\sigma_{Z_{max}}))\}$ cover $X^{tor}/K^*$. 
Moreover, we have the commutative diagram of Cartesian products

$$\CD
X(\pi(N),\pi(\partial_-\sigma_{Z_{max}})) @.@.@.@.@. X(\pi(N),\pi(\partial_+\sigma_{Z_{max}}))
\\ @. \searrow @.@.@. \swarrow @. \\
@VVV @.@.X(\pi(N),\pi(\sigma_{Z_{max}})) @.@.@VVV\\
(X^{tor})_-/K^* @.@.@.@.@. (X^{tor})_+/K^* \\
@. \searrow @.@.@VVV \swarrow @. \\
@. @.@.X^{tor}/K^*.@.@.\\
\endCD$$

\endproclaim

\demo{Proof}\enddemo We use the same notation as in Observation 3-2-1.  By
Lemma 3-2-4 and Lemma 3-2-5 we only have to show the situations
$(\heartsuit_j)$ for $j = 1, \cdot\cdot\cdot, l$.

First note that $\mu_X$ is characterized by the following properties:

\ \ (o) $\mu_X$ is proper birational,

\ \ (i) $X^{tor}$ is normal,

\ \ (ii) $\mu_X^{-1}(I_{\alpha,X}) \cdot {\Cal O}_{X^{tor}}$ is principal for
each
$\alpha \in S$, and

\ \ (iii) for any $\mu':X' \rightarrow X$ satisfying the properties (o), (i)
and (ii), there exists a unique $\nu':X' \rightarrow X^{tor}$ with $\mu' =
\mu_X
\circ \nu'$. 

\vskip.1in

Fix $j = 1, \cdot\cdot\cdot, l$.

\vskip.1in

First take 
$$\mu_j:X_j \rightarrow X$$
to be the normalization of the blowup of $X$ along $I_{\alpha_j,X}$.  Secondly,
let 
$$\tau_j:X^{tor} \rightarrow X_j$$
to be the normalization of the blowup of $X_j$
along the ideal
$$\prod_{\alpha \in S, \alpha \neq \alpha_j}\mu_i^{-1}(I_{\alpha,X}) \cdot {\Cal O}_{X_j}.$$
Then by the characterization above, we have
$$\mu_X = \mu_j \circ \tau_j.$$

We choose monomial generators of $K^*$-character $\alpha_j$ for $I_{\alpha_j,X}$ to be
$$v_j^*, m_1, \cdot\cdot\cdot, m_s$$
where $m_1, \cdot\cdot\cdot, m_s$ do not contain $v_j^*$, i.e.,
$$m_1, \cdot\cdot\cdot, m_s \in v_j^{\perp} \cap \sigma^{\vee} = \langle v_1^*,
\cdot\cdot\cdot, \overset{\vee}\to{v_j^*}, \cdot\cdot\cdot, v_l^*,\pm v_{l+1}^*,
\cdot\cdot\cdot, v_n^*\rangle.$$ 

\vskip.1in

Case: $v_j^*$ generates $I_{\alpha_j,X}$ (and hence we do not have any of $m_1,
\cdot\cdot\cdot, m_s$).

\vskip.1in

In this case, we have
$$D_{v_j} \subset D^{tor}.$$

\vskip.1in

Case: $v_j^*$ does not generate $I_{\alpha_j,X}$.

\vskip.1in

In this case, $X_j$ is covered by affine toric varieties which are the normalization
of the affine charts of the following types for the blowup of $X$ along
$I_{\alpha_j,X}$.

\ \ $\circ$ the chart obtained by inverting $v_j^*$: Over this chart, there is NO
intersection with $D_{v_j}$.

\ \ $\circ$ the chart(s) obtained by inverting, say $m_1$: The affine coordinate ring
for the chart (before normalization) is
$$K[v_j^* - m_1, m_2 - m_1, \cdot\cdot\cdot, m_s - m_1, v_1^*, \cdot\cdot\cdot,
\overset{\vee}\to{v_j^*}, \cdot\cdot\cdot, v_l^*, \pm v_{l+1}^*, \cdot\cdot\cdot,
\pm v_n^*]$$
where
$$m_2 - m_1, \cdot\cdot\cdot, m_s - m_1, v_1^*, \cdot\cdot\cdot,
\overset{\vee}\to{v_j^*}, \cdot\cdot\cdot, v_l^*, \pm v_{l+1}^*, \cdot\cdot\cdot,
\pm v_n^* \in v_j^{\perp} \cap M.$$
Therefore, the normalization has the form
$$K[v_j^* - m_1] \otimes K[\{{y'}_{h'}\}]$$
where the monomials $\{{y'}_{h'}\}$ are taken from $v_j^{\perp} \cap M$.

\proclaim{Lemma 3-2-9} Over this chart, 
$$\mu_j^{-1}(I_{\alpha,X}) \cdot {\Cal O}_{X_j} \hskip.1in \forall \alpha \in S,
\alpha \neq \alpha_j$$ 
is generated by monomials in $v_j^{\perp} \cap M$.
\endproclaim

\demo{Proof}\enddemo Let $m \in I_{\alpha,X}$ be a monomial in $v_1^*,
\cdot\cdot\cdot, v_l^*, \pm v_{l+1}^*, \cdot\cdot\cdot, \pm v_n^*$ of $K^*$-character
$\alpha$.  Write
$$\align
m &= cv_j^* + m'\\
&= c(v_j^* - m_1) + cm_1 + m' \\
\endalign$$
where $c \in {\Bbb Z}_{\geq 0}$ and $m'$ is a monomial in $v_j^{\perp} \cap M$.  Since
$v_j^* - m_1$ has the trivial $K^*$-character, $cm_1 + m'$ has the
$K^*$-character $\alpha$ and it is a monomial in $v_j^{\perp} \cap M$.  As these
monomials $\{m\}$ generate $\mu_j^{-1}(I_{\alpha,X})$, so do $\{cm_1 + m'\}$ and we
have the claim.

By the above lemma, over this chart $X^{tor}$, which is obtained as the normalization
of the blowup of $X_j$ along $\prod_{\alpha \in S, \alpha \neq
\alpha_j}\mu_j^{-1}(I_{\alpha,X}) \cdot {\Cal O}_{X_j}$, is covered by the affine
toric charts of the form
$$\roman{Spec}\{K[v_j^* - m_1] \otimes K[\{y_h\}]\}$$
where the monomials $\{y_h\}$ are taken from $v_j^{\perp} \cap M$ and where in the exact
sequence
$$0 \rightarrow v_j^{\perp} \cap M \rightarrow M \rightarrow {\Bbb Z} \rightarrow 0$$
the image of $v_j^* - m_1 \in M$ is $\pm 1 \in {\Bbb Z}$.
These affine charts satisfy the condition in the situation $(\heartsuit_j)$.

\vskip.1in

Varying the charts (i.e., inverting the monomial generators $m_1, \cdot\cdot\cdot,
m_s$), we see that the fan $\sigma^{tor}$, where $X^{tor} = X(N,\sigma^{tor})$, is
covered by the cones satisfying the condition in the situation $(\heartsuit_{j,fan})$. 
Therefore, we conclude that any cone in $\sigma^{tor}$ satisfies the condition in the
situation $(\heartsuit_j)$, that is to say, we have the situation $(\heartsuit_j)$.

This concludes the proof of the assertion (i) of Theorem 3-2-8.

\vskip.1in

Now we verify the assertion (ii) that $X^{tor}$ is a quasi-elementary cobordism.

Suppose, on the contrary, that there is a point $q \in X^{tor}$ such that $q$ is
not a fixed point $q
\not\in F_{X^{tor}}$ and that both limits $\lim_{t \rightarrow \infty}t(q)$ and
$\lim_{t \rightarrow 0}t(q)$ exists in $X^{tor}$.  

Then it is necessarily the case that
$\mu_X(q)
\in X$ is not a fixed point.  

In fact, consider the natural finite morphism
$$f:X^{tor} \rightarrow X \times {\Bbb P}^N$$
where the $N + 1$ homogeneous coordinates $(f_0, \cdot\cdot\cdot, f_N)$ of ${\Bbb P}^N$
are given by the monomial generators of the torific ideal $I_X$ of the same
$K^*$-character
$\prod_{\alpha
\in S}\alpha$.  If $\mu_X(q)$ is a fixed point, then $f(q) = \mu_X(q) \times
(f_0(q),
\cdot\cdot\cdot, f_N(q))$ is also a fixed point in $X \times {\Bbb P}^N$.  As $f$ is
a finite morphism, this would imply that $q$ is a fixed point, a contradiction. 

Since $\mu_X$ is proper, both limits $\lim_{t \rightarrow \infty}t(\mu_X(q)) =
\mu_X(\lim_{t
\rightarrow \infty}t(q))$ and \linebreak
$\lim_{t \rightarrow 0}t(\mu_X(q)) =
\mu_X(\lim_{t
\rightarrow 0}t(q))$ exist in $X$, implying that $X$ is not quasi-elementary.  This
is a contradiction, as $X$ is affine and hence quasi-elementary. 

\vskip.1in

Finally we verify the assertion (iii).

Observe that $q \in X^{tor}$, not being a fixed point, has the limit $\lim_{t \rightarrow
\infty}t(q)$ in $X^{tor}$ (resp. $\lim_{t \rightarrow
0}t(q)$) if and only if for the cone $\sigma_Z$ whose corresponding orbit $O(\sigma_Z)$
contains $q$ we have a cone $(\sigma_Z)^{\infty}$ (resp. $(\sigma_Z)^0$) such that there exists
a point $x \in \sigma_Z$ with $x + \epsilon a \in (\sigma_Z)^{\infty}$ (resp. for all $x -
\epsilon a \in (\sigma_Z)^0$ sufficiently small positive number $\epsilon$.  From this
observation and from the fact that $X^{tor}$ is quasi-elementary, it follows that
$X(N,\sigma_{Z_{max}})$ satisfies the condition $(*)$ for each of the maximal cones
$\sigma_{Z_{max}}$ in
$\sigma^{tor}$.  Thus we have the inclusion
$$X(N,\sigma_{Z_{max}})/K^* \hookrightarrow X^{tor}/K^*$$ 
and the commutative diagram follows from the main example.

\vskip.1in

This completes the proof of Theorem 3-2-8.

\vskip.1in

\proclaim{Remark 3-2-10}\endproclaim

(i) By the argument for the assertion (iii) we see that $X^{tor}$ being quasi-elemetary is
equivalent to the condition that NO two maximal cones $\sigma_{Z_{1,max}}$ and
$\sigma_{Z_{2,max}}$ with one being on top of the other: there exists a point $x \in
\sigma_{Z_{1,max}} \cap \sigma_{Z_{2,max}}$ such that $x + \epsilon a \in \sigma_{Z_{1,max}}$
and
$x - \epsilon a \in \sigma_{Z_{2,max}}$ for all sufficiently small positive number $\epsilon$. 

This in turn is equivalent to saying that the subdividion of $\sigma$ to obtain $\sigma^{tor}$
``comes from downstairs": Let $\pi:N_{\Bbb R} \rightarrow N_{\Bbb R}/{\Bbb R} \cdot a$
be the projection.  Then there is a fan $\pi(\sigma)^{tor}$ obtained by
subdividing $\pi(\sigma)$ so that the fan $\sigma^{tor}$ is covered by the
cones \linebreak
$\{\pi^{-1}(\tau)
\cap
\sigma;\tau \in \pi(\sigma)^{tor}\}$.

The last condition follows easily from the observation below, giving an alternative
interpretation for the assertion (iii):

If $f_0, \cdot\cdot\cdot, f_N$ are ther monomial generators of the torific
ideal $I_X$, then the blowup of $X$ along $I_X$ is covered by the affine
charts of the form (say, inverting $f_0$)
$$\roman{Spec}K[\sigma^{\vee} \cap M][\frac{f_1}{f_0}, \cdot\cdot\cdot,
\frac{f_N}{f_0}]$$
where $\frac{f_1}{f_0}, \cdot\cdot\cdot,
\frac{f_N}{f_0}$ all have the trivial $K^*$-character, i.e., they are all
$K^*$-invariants.

\vskip.1in

(ii) For a special choice of $S$ so that it contains a nonzero character $\beta$
which is divisible by all $\alpha_j$ for $j = 1, \cdot\cdot\cdot, l$  we can show
that
$$\pi(\partial_-\sigma^{tor}) = \pi(\partial_+\sigma^{tor})$$
and hence the two toroidal embeddings, obtained by excluding those
coordinate divisors $D_{\pi(v_i)}$ with $i \not\in S$ from the boundaries,
also coincide
$$(U_{X^{tor}},X^{tor})_-/K^* = (U_{X^{tor}},X^{tor})_+/K^*.$$
As we do not use this fact in our proof, we only refer th reader to
Abramovich-Karu-Matsuki-W{\l}odarczyk [1] for a proof.

\vskip.1in

\proclaim{Example 3-2-11}\endproclaim

Here we give an easy example to demonstrate the torification principle for
toric varieties.

\vskip.1in

Take $X = X(N,\sigma) = {\Bbb A}^3$ where $\sigma$ is the regular
cone 
$$\sigma = <v_1,v_2,v_3> \subset N_{\Bbb R}$$
and where $t \in K^*$ acts on the coordinates as
$$t(z_1, z_2, z_3) = (t^2z_1, tz_2, t^{-1}z_3),$$
i.e., it corresponds to
$$a = (2,1,-1) \in N.$$

We have the following monomial generators of the $\alpha$-torific ideals
 for $\alpha = 2, 1, -1$
$$\align
I_{2,X} &= \{z_1, z_1^2\} \\
I_{1,X} &= \{z_2, z_1z_3\} \\
I_{-1,X} &= \{z\}. \\
\endalign$$

Set the torific ideal $I_X$ (with respect to the set ${\goth C} = \{2,1,-1\}$
to be
$$I_X = I_{2,X}I_{1,X}I_{-1,X}.$$

Then the
normalization $({\Bbb A}^3)^{tor} = X(N,\sigma^{tor})$ of the blowup of $X$
along the torifific ideal $I_X$ is described by the fan $\sigma^{tor}$
consisting of the following three maximal cones
$$\align
\sigma_1 &= <v_1, v_3, 2v_1 + v_2> \\
\sigma_2 &= <v_2, v_1 + v_2, v_2 + v_3> \\
\sigma_3 &= <v_3, 2v_1 + v_2, v_1 + v_2, v_2 + v_3>. \\
\endalign$$

We have
$$D^{tor} = D_{v_3} \cup D_{v_1 + v_2} \cup D_{v_2 + v_3} \cup D_{2v_1 +
v_2}.$$
Thus the toroidal structure $(U_{X^{tor}},X^{tor})$ obtained from
$$(T_{X^{tor}} = X^{tor} - \{(D_{v_1} \cup D_{v_2}) \cup (D_{v_3} \cup D_{v_1 + v_2} \cup D_{v_2 + v_3} \cup D_{2v_1 +
v_2}),X^{tor})$$
by removing $D_{v_1}$ and $D_{v_2}$ from the boundary divisor.  In fact, we
observe the following.

\vskip.1in

The cone $\sigma_1$ has the product structure
$$\sigma_1 = \langle v_1\rangle \oplus \langle v_3, 2v_1 + v_2\rangle,$$
from which we have
$$X(N,\sigma_1) = \roman{Spec}\{K[v_1^* - 2v_2^*] \otimes K[v_2^*,v_3^*]\} \cong
{\Bbb A}_1 \times Z'.$$
with $D_{v_1} \cap X(N,\sigma_1)$ defined by $\{v_1^* - 2v_2^* = 0\}$. 

Thus we
see that even after removing $D_{v_1}$ from the boundary divisor
$$(U_{X^{tor}} \cap X(N,\sigma_1),X(N,\sigma_1)) = ({\Bbb A}_1 - \{0\},{\Bbb
A}_1) \times (T_{Z'},Z') \cup ({\Bbb A}_1 - \{1\},{\Bbb
A}_1) \times (T_{Z'},Z')$$
has the toroidal structure with respect to which $K^*$-action is toroidal.

The cone $\sigma_2$ has the product structure
$$\sigma_1 = \langle v_2\rangle \oplus \langle v_1 + v_2, v_2 + v_3\rangle,$$
from which we have
$$X(N,\sigma_1) = \roman{Spec}\{K[v_2^* - (v_1^* + v_3^*)] \otimes
K[v_1^*,v_3^*]\}
\cong {\Bbb A}_1 \times Z''.$$
with $D_{v_2} \cap X(N,\sigma_1)$ defined by $\{v_2^* - (v_1^* + v_3^*) = 0\}$. 

Thus we
see that even after removing $D_{v_2}$ from the boundary divisor
$$(U_{X^{tor}} \cap X(N,\sigma_2),X(N,\sigma_2)) = ({\Bbb A}_1 - \{0\},{\Bbb
A}_1) \times (T_{Z''},Z'') \cup ({\Bbb A}_1 - \{1\},{\Bbb
A}_1) \times (T_{Z''},Z'')$$
has the toroidal structure with respect to which $K^*$-action is toroidal.

The cone $\sigma_3$ does not contain $v_1$ or $v_2$ and hence
$$(U_{X^{tor}} \cap X(N,\sigma_3),X(N,\sigma_3)) =
(T_{X(N,\sigma_3)},X(N,\sigma_3))$$
has the original toric and hence structure with respect to which the
$K^*$-action is toric and hence toroidal.

\vskip.2in

Once the torification principle for toric varieties is established, it is
straightforward to verify the torifying property of the torific ideal for a
quasi-elementary cobordism via Luna's locally toric charts.

\proclaim{Theorem 3-2-12} Let $B$ be a nonsingular quasi-elementary cobordism
with $K^*$ acting effectively.  Let
$$\mu:B^{tor} \rightarrow B$$
be the normalization of the blowup of $B$ along the torific ideal $I$ and
$D^{tor}$ the divisor defined by the principal ideal $\mu^{-1}(I) \cdot {\Cal
O}_B$.  Then
$$(U_{B^{tor}} = B^{tor} - D^{tor},B^{tor})$$
has a toroidal structure with respect to which the induced $K^*$-action is
toroidal.

Moreover, $B^{tor}$ is a quasi-elementary cobordism, which induces a
$V$-toroidal birational map $\varphi$ as below:

$$\CD
(U_{B^{tor}},B^{tor})_-/K^* @.@.@.\overset{\varphi}\to{\dashrightarrow}@.@.
(U_{B^{tor}},B^{tor})_+/K^*
\\ @. \searrow @.@.@. \swarrow @. \\
@. @.@.(U_{B^{tor}},B^{tor})/K^*@.@.\\
\endCD$$

\endproclaim

\demo{Proof}\enddemo We take Luna's locally toric chart for $p \in B - F^*$
$$X \overset{\eta}\to{\leftarrow} V \overset{i}\to{\rightarrow} Y = U_p
\subset B$$
where $\eta$ and $i$ are strongly \'etale with $i$ being surjective and where
$Y = U_p$ satisfies the condition $(*)$.  (Note that the collection
$\{U_p;p \in B - F^*\}$ provides an open covering of $B$.)

\vskip.1in

We use the same notation as in the proof of Claim 3-1-6.

\vskip.1in

Since for any point $q \in Y - F_Y^* = Y \cap (B - F^*)$ (which implies $q_V \in V - F_V^*$
and $q_X \in X - F_X^*$), we have
$$\widehat{(\eta^*I_{\alpha,X})_{q_V}} =
\widehat{(I_{\alpha,V})_{q_V}} = \widehat{(i^*I_{\alpha})_{q_V}}$$ 
with all sides being
characterized as the ideal in $\widehat{{\Cal O}_{V,q_V}}$ generated by all the semi-invariant
functions of
$G_{q_V} (= G_q = G_{q_X})$-character $\alpha$ (See Lemma 3-1-7.), which implies
$$(\eta^*I_{\alpha,X})_{q_V} = (I_{\alpha,V})_{q_V} = (i^*I_{\alpha})_{q_V}.$$
Since $\eta^*I_{\alpha,X}, I_{\alpha,V}$ and $i^*I_{\alpha}$ are $K^*$-equivariant,
this implies that the above equality holds not only for $q \in Y - F_Y$ but also for all $q
\in Y$, i.e.,
$$\eta^*I_{\alpha,X} = I_{\alpha,V} = i^*I_{\alpha}.$$
Thus we conclude that the torific ideal on $X$ pulls back to the torific ideal on $V$,
coinciding with the torific ideal on $B$ pulled back to $V$, i.e.,
$$\eta^*I_X = I_V = i^*I.$$
This gives rise to the commutative diagram of Cartesian products
$$\CD
(U_{X^{tor}},X^{tor}) @<{\eta^{tor}}<< (U_{V^{tor}},V^{tor}) @>{i^{tor}}>>
(U_{Y^{tor}},Y^{tor}) \hskip.1in @.\subset \hskip.1in (U_{B^{tor}},B^{tor})
\\ 
@V{\mu_X}VV @V{\mu_V}VV @V{\mu_Y}VV @.\\
X @<\eta<< V @>i>> Y @.\subset \hskip.1in {\ \ \ \ \ B.\ \ \ \ \ \ \ \ \ }\\
\endCD$$

For any point $r \in Y^{tor}$, we can find by the torification principle for the toric
variety $X$ an affine open set
$Z_{\prod \{\delta_j\}}$ such that
$$r_{X^{tor}} \in (U_{Z_{\prod \{\delta_j\}}},Z_{\prod \{\delta_j\}}) := Z_{\prod
\{\delta_j\}}
\cap (U_{X^{tor}},X^{tor})
\cong
\prod ({\Bbb A}^1 -
\{\delta_j\},{\Bbb A}^1 - \{\delta_j\}) \times (T_{Z'},Z')$$
where $Z'$ is an affine toric variety with the torus $T_{Z'}$ and the $\delta_i$ are either 0
or 1.

In order to obtain Luna's toroidal chart for $q \in Y^{tor} \subset B^{tor}$
$$(U_{X_r},X_r) \leftarrow (U_{V_r},V_r) \rightarrow (U_{Y_r},Y_r) \subset
(U_{Y^{tor}},Y^{tor}) \subset (U_{B^{tor}},B^{tor}),$$
we set
$$\align
(U_{X_r},X_r) &:= (U_{Z_{\prod \{\delta_j\}}},Z_{\prod \{\delta_j\}}) \\
(U_{V_r},V_r) &:= (U_{Z_{\prod \{\delta_j\}}},Z_{\prod \{\delta_j\}}) \times_{X//K^*}V//K^* \\
(U_{Y_r},Y_r) &:= (U_{V_r},V_r) \\
\endalign$$
when $p \in Y$ is a fixed point and hence $i$ and $i^{tor}$ are isomorphisms by construction
of Luna's locally toric charts, and
$$\align
(U_{X_r},X_r) &:= (U_{Z_{\prod \{\delta_j\}}},Z_{\prod \{\delta_j\}}) \\
(U_{V_r},V_r) &:= (U_{Z_{\prod \{\delta_j\}}},Z_{\prod \{\delta_j\}}) \times_{X//K^*}V//K^*
\cap (i^{tor})^{-1}(Y_r)\\ 
(U_{Y_r},Y_r) &:= Y_r \cap (U_Y,Y) \\
\endalign$$
where $r \in Y_r$ is a small affine open neighborhood so that $Y_r \subset i^{tor}(V_r)$
when $p \in Y$ is not a fixed point and hence all the orbits in $Y, V, X$ and $Y^{tor},
V^{tor}, X^{tor}$ are closed by construction of Luna's locally toric charts.

Therefore, $(U_{B^{tor}},B^{tor})$ has a toroidal structure with respect to which the induced
$K^*$-action is toroidal.

\vskip.1in

It remains to prove that $B^{tor}$ is a quasi-elementary cobordism.

\vskip.1in

Suppose that $B^{tor}$ is not quasi-elementary, i.e., there is a point $q \in B^{tor}$ such
that $q$ is not a fixed point and that both limits $\lim_{t \rightarrow \infty}t(q)$ and
$\lim_{t \rightarrow 0}t(q)$ exist in $B^{tor}$.  Take $p \in B - F^*$ so that $Y = U_p$ (in
Luna's locally toric chart for $p$) cotains $\mu(q)$.  Note that since $U$ satisfies the
condition $(*)$, $Y^{tor} = U^{tor} \subset B^{tor}$ also satisfies the condition $(*)$.  In
particular, both limits sit inside of $Y^{tor}$.  On the other hand, since both $\eta^{tor}$
and $i^{tor}$ are strongly \'etale as so are $\eta$ and $i$, this implies that $q_{X^{tor}}
\in X^{tor}$ is not a fixed point and that both limits $\lim_{t \rightarrow
\infty}t(q_{X^{tor}})$ and
$\lim_{t \rightarrow 0}t(q_{X^{tor}})$ exist in $X^{tor}$.  This contradicts the fact by
Theorem 3-2-8 that $X^{tor}$ is quasi-elementary.

\vskip.1in

The commutative diagram indicating $\varphi$ is $V$-toroidal is an immediate consequence of
the toroidal structure $(U_{B^{tor}},B^{tor})$ (with respect to which the induced
$K^*$-action is toroidal) and $B^{tor}$ being quasi-elementary.

\vskip.1in

This completes the proof of Theorem 3-2-12.

\newpage

$$\bold{CHAPTER\ 4.\ RECOVERY\ OF\ NONSINGULARITY}$$

\vskip.2in

Let us recall what we have established so far.

\vskip.1in

In Chapter 2, we constructed a birational cobordism $B$ for a given birational map $\phi:X_1
\dashrightarrow X_2$ (after replacing it with some projective birational morphism by
eliminating points of indeterminacy by sequences of blowups with smooth centers).  The
birational cobordism gives rise to the decomposition of $\phi$ into a sequence of (V-)locally
toric birational maps

$$\phi:X_1 = W_1 \overset{\varphi_1}\to{\dashrightarrow} W_2
\overset{\varphi_2}\to{\dashrightarrow} \cdot\cdot\cdot W_i
\overset{\varphi_i}\to{\dashrightarrow} W_{i+1} \cdot\cdot\cdot
\overset{\varphi_{s-2}}\to{\dashrightarrow} W_{s-2} \overset{\varphi_{s-1}}\to{\dashrightarrow}
W_s = X_2$$ 
where each $\varphi_i$ corresponds to the commutative diagram

$$\CD
W_i @.@. \overset{\varphi_i}\to{\dashrightarrow} @.@. W_{i+1} \\
@|@.\hskip.6in @.@. @| \\
(B_{a_i})_-/K^* @.@. \hskip.6in @.@.(B_{a_i})_+/K^* \\
@.\searrow \hskip.2in@. \hskip.6in @.\hskip.2in \swarrow @.\\
@.@.(B_{a_i})/K^* @.@.\\
\endCD$$

induced by the quasi-elementary cobordism $B_{a_i}$, associated to the value $a_i$ of the strictly
increasing function from the set of the connected components of the fixed point set of $B$
to their characters in ${\Bbb Z}$.

\vskip.1in

The main object of Chapter 3 was to transform locally toric birational maps into toroidal
birational transformations.  By blowing up the torific ideal, we obtain a toroidal structure
$(U_{B_{a_i}^{tor}},B_{a_i}^{tor})$ with respect to which the $K^*$-action is toroidal and
hence which induces the toroidal birational map $\varphi_i^{tor}$ as in the commutative diagram
below

$$\CD
(U_{W_{i-}^{tor}},W_{i-}^{tor}) @.@.
\overset{\varphi_i^{tor}}\to{\dashrightarrow} @.@.
(U_{W_{i+}^{tor}},W_{i+}^{tor})
\\ @|@.\hskip.6in @.@. @| \\
(U_{B_{a_i}^{tor}},B_{a_i}^{tor})_-/K^* @.@. \hskip.6in
@.@.(U_{B_{a_i}^{tor}},B_{a_i}^{tor})_+/K^* \\ @.\searrow \hskip.2in@.
\hskip.6in @.\hskip.2in \swarrow @.\\ @.@.(U_{B_{a_i}^{tor}},B_{a_i}^{tor})/K^* @.@.\\
\endCD$$

Since the blowup morphism $\mu:B_{a_i}^{tor} \rightarrow B_{a_i}$ gives rise to a
commutative diagram

$$\CD
(U_{B_{a_i}^{tor}},B_{a_i}^{tor})_-/K^* @.@. \hskip.6in
@.@.(U_{B_{a_i}^{tor}},B_{a_i}^{tor})_+/K^* \\ 
@.\searrow \hskip.2in@.
\hskip.6in @.\hskip.2in \swarrow @.\\ 
@VVV @.(U_{B_{a_i}^{tor}},B_{a_i}^{tor})/K^* @.@.@VVV \\
(B_{a_i})_-/K^* @.@.\hskip.6in @.@.(B_{a_i})_+/K^* \\
@.\searrow \hskip.2in@. @VVV \hskip.2in \swarrow @.\\
@.@.(B_{a_i})/K^* @.@.\\
\endCD$$

we actually factor $\varphi_i$ into a sequence of birational maps
$$\varphi_i:W_i \overset{f_{i-}}\to{\leftarrow} (U_{W_{i-}^{tor}},W_{i-}^{tor})
\overset{\varphi_i^{tor}}\to{\dashrightarrow} (U_{W_{i+}^{tor}},W_{i+}^{tor})
\overset{f_{i+}}\to{\rightarrow} W_{i+1}.$$

If $W_i,W_{i-}^{tor},W_{i+}^{tor}$ and $W_{i+1}$ were all nonsingular and both
$f_{i-}$ and $f_{i+}$ were sequences of blowups with smooth centers, then we would only have to
apply the strong factorization theorem of Morelli [1,2] or
Abramovich-Matsuki-Rashid] [1] (cf.W{\l}odarczyk [1]) to the toroidal birational map
$\varphi_i^{tor}$ between nonsingular toroidal embeddings to finish the proof of the weak
factorization theorem.  But in general these varieties are all singular, while $f_{i-}$ and
$f_{i+}$ are not sequences of blowups with smooth centers.

\vskip.1in

The purpose of Chapter 4 is to bring this SINGULAR situation back to a NONSINGULAR one by an
application of the canonical resolution of singularities.

\vskip.1in

More precisely, we construct the following commutative diagram

$$\CD
W_i^{res} @<{p_{i-}}<< W_{i-}^{can} @.
\hskip.2in \overset{\varphi_i^{can}}\to{\dashrightarrow} \hskip.2in @. W_{i+}^{can}
@>{p_{i+}}>> W_{i+1}^{res} \\
@V{r_{i-}}VV @V{h_{i-}}VV @. @VV{h_{i+}}V @VV{r_{i+}}V \\
W_i @<{f_{i-}}<< W_{i-}^{tor} @.
\hskip.2in \overset{\varphi_i^{tor}}\to{\dashrightarrow} \hskip.2in @. W_{i+}^{tor}
@>{f_{i+}}>> W_{i+1} \\
\endCD$$
where 

\ \ (i) $W_i^{res},W_{i-}^{can},W_{i+}^{can}$ and $W_{i+1}^{res}$ are all nonsingular, 

\ \ (ii) the
morphisms $p_{i-}$ and $p_{i+}$ are sequences of blowups with smooth centers, and

\ \ (iii) by setting
$$\align
U_{W_{i-}^{can}} &= h_{i-}^{-1}(U_{W_{i-}^{tor}}) \\
U_{W_{i+}^{can}} &= h_{i+}^{-1}(U_{W_{i+}^{tor}}) \\
\endalign$$
the birational map
$$\varphi_i^{can}:(U_{W_{i-}^{can}},W_{i-}^{can}) \dashrightarrow
(U_{W_{i+}^{can}},W_{i+}^{can})$$
is V-toroidal via the commutative diagram

$$\CD
(U_{W_{i-}^{can}},W_{i-}^{can}) @.@. \overset{\varphi_i^{can}}\to{\dashrightarrow} @.@.
(U_{W_{i+}^{can}},W_{i+}^{can})
\\ @V{h_{i-}}VV @.@.@.@VV{h_{i+}}V \\
(U_{W_{i-}^{tor}},W_{i-}^{tor}) @.@. \overset{\varphi_i^{tor}}\to{\dashrightarrow} @.@.
(U_{W_{i+}^{can}},W_{i+}^{tor})
\\
@.\searrow @. \hskip.6in @. \swarrow @. \\
@.@. (U_{B_{a_i}^{tor}},B_{a_i}^{tor})/K^*. @.@.\\
\endCD$$

After the construction, we only have to apply the strong factorization theorem of toroidal
birational maps to $\varphi_i^{can}$ to finish the proof of the weak factorization theorem. 

\vskip.2in

\S 4-1. Application of the canonical resolution of singularities

\vskip.1in

\proclaim{Remark 4-1-1}\endproclaim
 
(i) By the $\bold{canonical\ resolution\ of\ singularities}$,
we mean an algorithm which, given a variety $X$, provides a uniquely determined sequence of
blowups with smooth centers $r:X^{res} \rightarrow X$ from a nonsingular variety, satisfying the
conditions:

($\spadesuit^{res}-0$) The centers of blowups are taken only over the singular locus $Sing(X)$.

($\spadesuit^{res}-1$) If $Y \rightarrow X$ is a smooth morphism of varieties, then the canonical
resolution of $Y$ is obtained by pulling back that of $X$, i.e., we have a Cartesian product
$$\CD
Y @<<< Y \times_X X^{res} @.= Y^{res} \\
@VVV @VVV @.\\
X @<<< X^{res}.@.\\
\endCD$$

In particular, the canonical resolution can be pulled back by an open immersion or an \'etale
morphism.  

The condition ($\spadesuit^{res}-1$) also implies that any automorphism of $X$ (not necessarily
over $K$) lifts to that of the canonical resolution $X^{res}$.

Moreover, the condition $(\spadesuit^{res}-1$) implies that a family of automorphisms
$\theta_g$ of a variety $X$ parametrized by a smooth variety $g \in G$ lifts to the family of
automorphisms of the canonical resolution $X^{res}$ (for example an action of an algebraic group
$G$ on a variety $X$).  In fact, applying
the condition ($\spadesuit^{res}-1$) to the map parametrizing the automorphisms of $X$
$$\theta:G
\times X
\rightarrow X \text{ given\ by\ }\theta(g,x) = \theta_g(x)$$
and applying it also to the projection
onto the second factor $pr_2:G
\times X
\rightarrow X$ we obtain the map parametrizing the automorphisms of $X^{res}$
$$\theta^{res}:G \times X^{res} = (G \times X) {\times_{pr_2}}_X\ X^{res} = (G \times X)
{\times_{\theta}}_X\ X^{res} = Y^{res} \rightarrow X^{res}.$$  

(ii) By the $\bold{canonical\ principalization\ of\ ideals}$, we mean an algorithm which, given an
ideal sheaf $I$ on a nonsingular variety $X$, provides a uniquely determined sequence of blowups
with smooth (and admissible) centers $p:X^{can} \rightarrow X$ from a nonsingular variety with
the property $p^{-1}(I) \cdot {\Cal O}_{X^{can}}$ being principal, satisfying the conditions:

($\spadesuit^{can}-0$) The centers of blowups are taken over the support of ${\Cal O}_X/I$. 

($\spadesuit^{can}-1$) If $f:Y \rightarrow X$ is a smooth morphism of nonsingular varieties (not
necessarily over the base field $K$), then the canonical principalization of the ideal sheaf
$f^*I$ is obtained by pulling back that of
$I$, i.e., we have a Cartesian product
$$\CD
Y @<<< Y \times_X X^{can} @.= Y^{can} \\
@VVV @VVV @.\\
X @<<< X^{can}.@.\\
\endCD$$

In particular, the canonical principalization (of the corresponding ideals) can be pulled back by
an open immersion or an
\'etale morphism.

Also the condition ($\spadesuit^{can}-1$) implies that any automorphism of $X$ stabilizing
the ideal $I$ (not necessarily over $K$) lifts to that of the canonical principalization
$X^{can}$. 

Moreover, the condition $(\spadesuit^{can}-1$) implies that a family of automorphisms
$\theta_g$ of a variety $X$ parametrized by a smooth variety $g \in G$ lifts to the family of
automorphisms of the canonical pricipalization $X^{can}$ (for example an action of an algebraic
group
$G$ on a variety $X$), if the
ideal $I$ is stable under the auomorphisms, i.e.,
$\theta_g^*I = I \hskip.1in \forall g \in G$.  In fact, applying
the condition ($\spadesuit^{can}-1$) to the map parametrizing the automorphisms of $X$
$$\theta:G \times X \rightarrow X \text{\ given\ by\ }\theta(g,x) = \theta_g(x)$$
and applying it to
the projection onto the second factor $pr_2:G \times X \rightarrow X$ and also noting that the
stability implies $\theta^*I = pr_2^*I$, we obtain the map parametrizing the automorphisms of
$X^{can}$ as before
$$\theta^{can}:G \times X^{can} = (G \times X) {\times_{pr_2}}_X\ X^{can} = (G \times X)
{\times_{\theta}}_X\ X^{can} = Y^{can} \rightarrow X^{can}.$$ 

(iii) All the existing algorithms for the canonical resolution of singularities and/or canonical
principalization of ideals (cf. Bierstone-Milman [1] Villamayor [1] \linebreak
Encinas-Villamayor [1] and Hironaka [2]) satisfy the requirements in (i) and (ii) as above.  As
these are the only properties we utilize in our proof, one may choose any one of the algorithms
and we denote the resulting varieties $X^{res}$ and $X^{can}$ without specifying the choice.

\vskip.1in

We go back to the construction of the commutative diagram for the purpose of recovery from
singular to nonsingular.

\vskip.1in

We take $W_i^{res}$ and $W_{i+1}^{res}$ to be the canonical resolution of singularities of
$W_i$ and $W_{i+1}$, respectively.

\vskip.1in

\proclaim{Lemma 4-1-2} The morphism $f_{i-}:W_{i-}^{tor} \rightarrow W_i$
is a projective morphism (in the sense of Grothendieck [1] but not in the sense of Hartshorne
[1]) with an effective Cartier divisor $E_{i-}^{tor}$ such that

\ \ (i) the support of $E_{i-}^{tor}$ lies in $\{D^{tor} \cap
(B_{a_i}^{tor})_-\}/K^*$, where $D^{tor}$ is the divisor defined by the principal ideal
$\mu^{-1}(I) \cdot {\Cal O}_{B_{a_i}^{tor}}$ with $\mu:B_{a_i}^{tor} \rightarrow B_{a_i}$ being
the normalization of the blowup of the torific ideal $I$, 

\ \ (ii) ${\Cal O}_{W_{i-}^{tor}}(- E_{i-}^{tor})$ is $f_{i-}$-ample, and

\ \ (iii) the morphism $f_{i-}$ is the blowup of the ideal
$$I_{i-} = (f_{i-})_*{\Cal O}_{W_{i-}^{tor}}(- E_{i-}^{tor}).$$

Similarly, the morphism $f_{i+}:W_{i+}^{tor} \rightarrow W_{i+1}$ is a projective morphism with
an effective Cartier divisor satisfying the conditions (i), (ii) and (iii) with the negative
sign $-$ replaced by the positive one $+$.
\endproclaim

\demo{Proof}\enddemo Since $\mu:B_{a_i}^{tor} \rightarrow B_{a_i}$ is the normalization of the
blowup of
$B_{a_i}$ along the torific ideal $I$, we can take an effective Cartier divisor $E_i^{tor}$ such
that the principal ideal $\mu^{-1}(I) \cdot {\Cal O}_{B_{a_i}^{tor}} = {\Cal O}_{B_{a_i}^{tor}}(-
E_i^{tor})$ is $\mu$-ample and hence $\mu|_{(B_{a_i}^{tor})_-}$-ample when restricted to
$\mu|_{(B_{a_i}^{tor})_-}:(B_{a_i}^{tor})_- \rightarrow (B_{a_i})_-$.

Now if we look at Luna's locally toric chart 
$$X = X(N,\sigma) = X \overset{\eta}\to{\leftarrow}
V \overset{i}\to{\rightarrow} Y = U \subset B_{a_i},$$
the morphisms
$$\align
(B_{a_i}^{tor})_- &\rightarrow (B_{a_i})_- \\
(B_{a_i}^{tor})_-/K^* &\rightarrow (B_{a_i})_-/K^* \\
\endalign$$
correspond to the toric morphisms
$$\align
X(N,\partial_-\sigma^{tor}) &\rightarrow X(N,\partial_-\sigma) \\
X(\pi(N),\pi(\partial_-\sigma^{tor})) &\rightarrow X(\pi(N),\pi(\partial_-\sigma)),\\
\endalign$$
where $\pi:{\Bbb N}_{\Bbb R} \rightarrow {\Bbb N}_{\Bbb R}/a \cdot {\Bbb R}$ is the projection
(cf. the notation in Proposition 3-1-8 or Observation 3-2-1).

Remark also that since the torific ideal $I$ corresponds to the toric torific ideal $I_X$, i.e.,
$\eta^*I_X = i^*I$, the Cartier divisor $E_i^{tor}$ corresponds to the toric Cartier divisor
$(E_i^{tor})_{X^{tor}}$ with $\eta^*(E_i^{tor})_{X^{tor}} = i^*E_i^{tor}$.

Say, the toric Cartier divisor $(E_i^{tor})_{X^{tor}}$ is associated to a piecewise linear
function
$\psi_{(E_i^{tor})_{X^{tor}}}$ on the fan $\sigma^{tor}$ (cf. Fulton [1] or Oda [1]).  Since
$\pi:\partial_-\sigma^{tor} \rightarrow \pi(\partial_-\sigma^{tor})$ is a linear isomorphism of
fans (without taking the lattices $N$ and $\pi(N)$ into consideration), there exists a piecewise
linear function $\psi_{(E_{i-}^{tor})_{(X^{tor})_-/K^*}}$ on the fan $\pi(\partial_-\sigma^{tor})$
such that
$$\psi_{(E_i^{tor})_{X^{tor}}|_{\partial_-\sigma^{tor}}} =
\pi^*(\psi_{(E_{i-}^{tor})_{(X^{tor})_-/K^*}}).$$ 
Take $(E_{i-}^{tor})_{(X^{tor})_-/K^*}$ to be the effective ${\Bbb Q}$-Cartier divisor
associated to the piecewise linear function $\psi_{(E_{i-}^{tor})_{(X^{tor})_-/K^*}}$.  (Recall
that $(E_{i-}^{tor})_{(X^{tor})_-/K^*}$ is only guaranteed to be ${\Bbb Q}$-Cartier, since in
general
$\pi:N|_{\partial_-\sigma^{tor}} \rightarrow \pi(N)|_{\pi(\partial_-\sigma^{tor})}$ is not an
isomorphism.)

\vskip.1in

Observe the following implications:

\vskip.1in

[${\Cal O}_{X^{tor}}(- (E_i^{tor})_{X^{tor}})$ is
$\mu_X$-ample]

\vskip.1in

$\Longrightarrow$

\vskip.1in

[${\Cal O}_{X^{tor}}(- (E_i^{tor})_{X^{tor}})|_{(X^{tor})_-}$ is
$\mu_X|_{(X^{tor})_-}$-ample]

\vskip.1in

$\Longleftrightarrow$

\vskip.1in

[$- \psi_{(E_i^{tor})_{X^{tor}}|_{\partial_-\sigma^{tor}}}$ is strictly convex when restricted
to each cone in $\partial_-\sigma$ with respect to $\partial_-\sigma^{tor}$]

\vskip.1in

$\Longleftrightarrow$

\vskip.1in

[$- \psi_{(E_{i-}^{tor})_{(X^{tor})_-/K^*}}$ is strictly convex when restricted
to each cone in $\pi(\partial_-\sigma)$ with respect to $\pi(\partial_-\sigma^{tor})$]

\vskip.1in

$\Longleftrightarrow$

\vskip.1in

[${\Cal O}_{(X^{tor})_-/K^*}(- (E_{i-}^{tor})_{(X^{tor})_-/K^*})$ is $(f_{i-})_X$-ample where
$(f_{i-})_X$ is the morphism $(f_{i-})_X:(X^{tor})_-/K^* \rightarrow X_-/K^*$.]

\vskip.2in

Therefore, if we take $E_{i-}^{tor}$ to be the effective ${\Bbb Q}$-Cartier divisor on
$W_{i-}^{tor}$ corresponding to $(E_{i-}^{tor})_{X^{tor}}$, then the conditions (i) and (ii) are
satisfied.

By replacing it with a sufficiently divisible and high multiple, we see that there exists an
effective Cartier divisor $E_{i-}^{tor}$ satisfying all the conditions (i), (ii) and (iii).

\vskip.1in

This completes the proof of Lemma 4-1-2.

\vskip.1in

We take $W_{i-}^{can}$ (resp. $W_{i+}^{can}$) to be canonical principalization of the ideal
$(r_{i-})^{-1}(I_{i-}) \cdot {\Cal O}_{W_i^{res}}$ (resp. $(r_{i+})^{-1}(I_{i+}) \cdot {\Cal
O}_{W_{i+1}^{res}}$).

\proclaim{Proposition 4-1-3} Let
$$\align
h_{i-}:W_{i-}^{can} &\rightarrow W_{i-}^{tor} \\
h_{i+}:W_{i+}^{can} &\rightarrow W_{i+}^{tor} \\
\endalign$$
be the induced morphisms and set
$$\align
U_{W_{i-}^{can}} &= h_{i-}^{-1}(U_{W_{i-}^{tor}}) \\
U_{W_{i+}^{can}} &= h_{i+}^{-1}(U_{W_{i+}^{tor}}). \\
\endalign$$

Then

\ \ (i) both $(U_{W_{i-}^{can}},W_{i-}^{can})$ and $(U_{W_{i+}^{can}},W_{i+}^{can})$ have the
toroidal structures,

\ \ (ii) the induced birational map $\varphi_i^{can}:(U_{W_{i-}^{can}},W_{i-}^{can})
\dashrightarrow (U_{W_{i+}^{can}},W_{i+}^{can})$ is V-toroidal via the commutative diagram

$$\CD
(U_{W_{i-}^{can}},W_{i-}^{can}) @.@.\hskip.6in @.@. (U_{W_{i+}^{can}},W_{i+}^{can}) \\
@V{h_{i-}}VV @.@.@.@VV{h_{i+}}V \\
(U_{W_{i-}^{tor}},W_{i-}^{tor}) @.@.@.@. (U_{W_{i+}^{tor}},W_{i+}^{tor}) \\
@| @.@.@. @| \\
(U_{B_{a_i}^{tor}},B_{a_i}^{tor})_-/K^* @.@. \hskip.6in
@.@.(U_{B_{a_i}^{tor}},B_{a_i}^{tor})_+/K^* \\ 
@.\searrow \hskip.2in@.
\hskip.6in @.\hskip.2in \swarrow @.\\ 
@.@.(U_{B_{a_i}^{tor}},B_{a_i}^{tor})/K^*. @.@. \\
\endCD$$

\endproclaim

\demo{Proof}\enddemo We take Luna's locally toric chart for $p \in B_{a_i} - F_{a_i}^*$ as before
$$X \overset{\eta}\to{\leftarrow} V \overset{i}\to{\rightarrow} Y = U_p \subset B_{a_i}$$
where $\eta$ and $i$ are strongly \'etale with $i$ being surjective and where $Y = U_p$
satisfies the condition ($*$).

This gives rise to the commutative diagram of Cartesian products
$$\CD
X^{tor} @<{\eta^{tor}}<< V^{tor} @>{i^{tor}}>> Y^{tor} \hskip.1in @.\subset
\hskip.1in B_{a_i}^{tor} \\ 
@V{\mu_X}VV @V{\mu_V}VV @V{\mu_Y}VV @.\\
X @<\eta<< V @>i>> Y @.\subset \hskip.1in B_{a_i}.\\
\endCD$$

Since the canonical resolution is pulled-back by \'etale morphisms and since so does the
canonical principalization of ideals, noting
$$\{\eta^{tor}\}^*\{\mu_X^{-1}I_X \cdot {\Cal O}_{X^{tor}}\} = \{i^{tor}\}^*\{\mu_Y^{-1}I_Y
\cdot {\Cal O}_{Y^{tor}}\},$$
which implies
$$\align
\{\eta_-/K^*\}^*(I_{i-})_{X_-/K^*} &=
\{i_-/K^*\}^*(I_{i-})_{Y_-/K^*}
\\
\{\eta_+/K^*\}^*(I_{i+})_{X_+/K^*} &=
\{i_+/K^*\}^*(I_{i+})_{Y_+/K^*},
\\ 
\endalign$$
we have the commutative diagram of Cartesian products

$$\matrix
X_-/K^* & \leftarrow & V_-/K^* & \rightarrow &
Y_-/K^* &
\subset & W_i \\
\nearrow \hskip.2in && \nearrow \hskip.2in && \nearrow \hskip.2in && \nearrow \hskip.2in \\
(X_-/K^*)^{res} \hskip.3in &\leftarrow \hskip.3in & (V_-/K^*)^{res} \hskip.3in & \rightarrow
\hskip.3in & (Y_-/K^*)^{res}
\hskip.3in & \subset \hskip.3in & W_{i-}^{res} \hskip.3in
\\
\uparrow &  & \uparrow &  & \uparrow & 
& \uparrow \\
\uparrow \hskip.5in && \uparrow \hskip.5in && \uparrow \hskip.5in && \uparrow \hskip.5in \\
(X^{tor})_-/K^* & \leftarrow & (V^{tor})_-/K^* & \rightarrow & (Y^{tor})_-/K^* & \subset
& W_{i-}^{tor} \\
\nearrow \hskip.2in && \nearrow \hskip.2in && \nearrow \hskip.2in && \nearrow \hskip.2in \\
(X_-/K^*)^{can} \hskip.3in &\leftarrow \hskip.3in & (V_-/K^*)^{can} \hskip.3in & \rightarrow
\hskip.3in & (Y_-/K^*)^{can}
\hskip.3in & \subset \hskip.3in & W_{i-}^{can} \hskip.3in
\\
\downarrow &  & \downarrow &  & \downarrow & 
& \downarrow \\
X^{tor}/K^* & \leftarrow & V^{tor}/K^* & \rightarrow & Y^{tor}/K^* & \subset
& B_{a_i}/K^* \\
\uparrow &  & \uparrow &  & \uparrow & 
& \uparrow \\
(X^{tor})_+/K^* & \leftarrow & (V^{tor})_+/K^* & \rightarrow & (Y^{tor})_+/K^* & \subset
& W_{i+}^{tor} \\
\nearrow \hskip.2in && \nearrow \hskip.2in && \nearrow \hskip.2in && \nearrow \hskip.2in \\
(X_+/K^*)^{can} \hskip.3in &\leftarrow \hskip.3in & (V_+/K^*)^{can} \hskip.3in & \rightarrow
\hskip.3in & (Y_+/K^*)^{can}
\hskip.3in & \subset \hskip.3in & W_{i+}^{can} \hskip.3in
\\
\downarrow &  & \downarrow &  & \downarrow & 
& \downarrow \\
\downarrow \hskip.5in && \downarrow \hskip.5in && \downarrow \hskip.5in && \downarrow \hskip.5in \\
X_+/K^* & \leftarrow & V_+/K^* & \rightarrow & Y_+/K^* & \subset
& W_{i+} \\
\nearrow \hskip.2in && \nearrow \hskip.2in && \nearrow \hskip.2in && \nearrow \hskip.2in \\
(X_+/K^*)^{res} \hskip.3in & \leftarrow \hskip.3in & (V_+/K^*)^{res}
\hskip.3in &
\rightarrow
\hskip.3in & (Y_+/K^*)^{res}
\hskip.3in & \subset \hskip.3in & W_{i+}^{res} \hskip.3in
\\
\endmatrix$$

Therefore, we only have to show the assertion assuming that the quasi-elementary cobordism
$B_{a_i} = X$ is a nonsingular affine toric variety.

\vskip.1in

Note first that the morphisms
$$\align
f_{i-}:W_{i-}^{tor} = (X^{tor})_-/K^* &\rightarrow W_i = X_-/K^* \\ 
f_{i+}:W_{i+}^{tor} = (X^{tor})_+/K^* &\rightarrow W_i = X_+/K^* \\ 
\endalign$$
are equivariant (toric) ones between toric varieties by construction, that so are the morphisms
$$\align
r_{i-}:W_i^{res} = (X_-/K^*)^{res} &\rightarrow W_i = X_-/K^* \\
r_{i+}:W_i^{res} = (X_+/K^*)^{res} &\rightarrow W_{i+1} = X_+/K^* \\
\endalign$$
since the canonical resolution of singularities preserves the action of the embedded torus (See
the condition ($\spadesuit^{res}-1$).), that the morphisms
$$\align
p_{i-}:W_{i-}^{can} = (X_-/K^*)^{can} &\rightarrow W_i^{res} = (X_-/K^*)^{res} \\
p_{i+}:W_{i+}^{can} = (X_+/K^*)^{can} &\rightarrow W_i^{res} = (X_+/K^*)^{res} \\
\endalign$$
since the ideals $I_{i-}$ and $I_{i+}$ (and their pull-backs by $r_{i-}$ and $r_{i+}$) are both
toric by construction (cf. Lemma 4-1-2) and since the canonical principalization also preserves
the action of the embedded torus (See the condition ($\spadesuit^{can}-1$).) and that hence so
are the morphisms
$$\align
h_{i-}:W_{i-}^{can} = (X_-/K^*)^{can} &\rightarrow W_{i-}^{tor} = (X_-/K^*)^{tor} \\
h_{i+}:W_{i+}^{can} = (X_+/K^*)^{can} &\rightarrow W_{i+}^{tor} = (X_+/K^*)^{tor}. \\
\endalign$$

We denote the corresponding fans and associated troic varieties by

$$\align
W_i &= X(\pi(N),\pi(\partial_-\sigma)) \\
W_i^{res} &= X(\pi(N),\Sigma_i^{res}) \\
W_{i-}^{tor} &= X(\pi(N),\Sigma_-^{tor}) \\
W_{i-}^{can} &= X(\pi(N),\Sigma_-^{can}) \\
W_{i+1} &= X(\pi(N),\pi(\partial_+\sigma)) \\
W_{i+1}^{res} &= X(\pi(N),\Sigma_{i+1}^{res}) \\
W_{i+}^{tor} &= X(\pi(N),\Sigma_+^{tor}) \\
W_{i+}^{can} &= X(\pi(N),\Sigma_+^{can}) \\
\endalign$$

Let $\{\sigma_{Z_{max}}\}$ be the set of the maximal cones of $\sigma^{tor}$, where
$\sigma^{tor}$ is the fan corresponding to the toric variety $X^{tor} = X(N,\sigma^{tor})$.  By
Theorem 3-2-8 (i) we have the situation $(\heartsuit_{fan})$.  Therefore, for any
$\sigma_{Z_{max}}$ we have the decomposition
$$\sigma_{Z_{max}} = \Sigma_{\{j;j \not\in S, v_j \in \sigma_{Z_{max}}\}} \langle v_j \rangle
\oplus
\tau_{Z'_{max}}
\text{\ as\ cones},$$ where the cone
$\tau_{Z'_{max}}$ is generated by the extremal rays
$w_k's$ (other than those
$v_j$'s with $j \not\in S$) of $\sigma_{Z_{max}}$ and is contained in a linear space $L$ (of
codimension $\#\{j;j \not\in S, v_j \in \sigma_{Z_{max}}\}$) which also contains $a \in N$ and
where we have the exact sequence
$$0 \rightarrow L \cap N \rightarrow N \rightarrow {\Bbb Z}^{\oplus \#\{j;j \not\in S,v_i \in
\sigma_{Z,max}\}} \rightarrow$$
with the image of $\{v_j;j \not\in S,v_i \in \sigma_{Z_{max}}\}$ forming a ${\Bbb Z}$-basis of
${\Bbb Z}^{\#\{j;j \not\in S,v_j \in
\sigma_{Z_{max}}\}}$.

Therefore, under the projection map $\pi:N_{\Bbb R} \rightarrow N/a \cdot N_{\Bbb R}$ we have
$$\pi(N) = \Sigma_{\{j;j \not\in S, v_j \in \sigma_{Z,max}\}} {\Bbb Z} \cdot \pi(v_j) \oplus \pi(N
\cap L)$$ and
$$\align
\pi(\partial_-\sigma_{Z_{max}}) &= \Sigma_{\{j;j \not\in S, v_j \in \sigma_{Z_{max}}\}} \langle
\pi(v_j)
\rangle \oplus \pi(\partial_-\tau_{Z'_{max}}) \\
\pi(\sigma_{Z_{max}}) &= \Sigma_{\{j;j \not\in S, v_j \in \sigma_{Z_{max}}\}} \langle
\pi(v_j)
\rangle \oplus \pi(\tau_{Z'_{max}}) \\
\pi(\partial_+\sigma_{Z_{max}}) &= \Sigma_{\{j;j \not\in S, v_j \in \sigma_{Z_{max}}\}} \langle
\pi(v_j)
\rangle \oplus \pi(\partial_+\tau_{Z'_{max}}). \\
\endalign$$

By Theorem 3-2-8 (iii) we have the commutative diagram of Cartesian products

$$\CD
X(\pi(N),\pi(\partial_-\sigma_{Z_{max}})) @.@.@.@.@. X(\pi(N),\pi(\partial_+\sigma_{Z_{max}})) \\
@. \searrow @.@.@. \swarrow @. \\
@VVV @.@.X(\pi(N),\pi(\sigma_{Z_{max}})) @.@.@VVV\\
(X^{tor})_-/K^* @.@.@.@.@. (X^{tor})_+/K^* \\
@. \searrow @.@.@VVV \swarrow @. \\
@. @.@.X^{tor}/K^*.@.@.\\
\endCD$$

Since the affine toric varieties $\{X(\pi(N),\pi(\sigma_{Z_{max}}))\}$ cover $X^{tor}/K^*$,
we conclude that $\{X(\pi(N),\pi(\partial_-\sigma_{Z_{max}}))\}$ cover $(X^{tor})_-/K^* =
X(\pi(N),\Sigma_-^{tor}))$.  That is to say, the cones
$\{\pi(\partial_-\sigma_{Z_{max}})\}$ cover the fan
$\Sigma_-^{tor}$.  Thus the fan $\Sigma_-^{can}$, which is a refinement of $\Sigma_-^{tor}$, is
the union of the refinements $\pi(\partial_-\sigma_{Z_{max}})^{can}$ of
$\pi(\partial_-\sigma_{Z_{max}})$.  Also we have the identical statement replacing the $-$ sign
with the $+$ sign.

\vskip.1in

We prove the following situations $(\heartsuit_{\pi(\partial_-fan)})$ and
$(\heartsuit_{\pi(\partial_+fan)})$ hold:
$$(\heartsuit_{\pi(\partial_-fan)}) \left\{
\aligned
&\pi(\partial_-{\sigma_{Z_{max}}})^{can} \text{has\ the\ product\ structure\ } \\
&\pi(\partial_-{\sigma_{Z_{max}}})^{can} = \Sigma_{\{j;j \not\in S, v_j \in \sigma_{Z_{max}}\}}
\langle
\pi(v_j)
\rangle \oplus \pi(\partial_-\tau_{Z'_{max}})^{can} \\
&\text{inherited\ from\ the\ product\ structure\ of\ }\pi(\partial_-\sigma_{Z_{max}}) \text{\ and\
}\pi(\sigma_{Z_{max}})\\ 
&\pi(\partial_-\sigma_{Z_{max}}) = \Sigma_{\{j;j \not\in S, v_j \in
\sigma_{Z_{max}}\}}
\langle
\pi(v_j)
\rangle \oplus \pi(\partial_-\tau_{Z'_{max}}) \\
&\pi(\sigma_{Z_{max}}) = \Sigma_{\{j;j \not\in S, v_j \in
\sigma_{Z_{max}}\}}
\langle
\pi(v_j)
\rangle \oplus \pi(\tau_{Z'_{max}}) \\
&\text{where\ }\pi(\partial_-\tau_{Z'_{max}})^{can}\text{\ is\ a\ refinement\ of\
}\pi(\partial_-\tau_{Z'_{max}}).\\ 
\endaligned
\right.$$

The situation $(\heartsuit_{\pi(\partial_+fan)})$ is idetical to
$(\heartsuit_{\pi(\partial_-fan)})$, replacing the $-$ sign with the $+$ sign.

\vskip.1in

In fact, since the toroidal structure
$$X(\pi(N),\pi(\sigma_{Z_{max}})) \cap (U_{X^{tor}},X^{tor})
= (U_{X(\pi(N),\pi(\sigma_{Z_{max}})},X(\pi(N),\pi(\sigma_{Z_{max}}))$$
is covered by
the open sets of the form
$$\align
\prod_{\{j;j \in S, v_j \in
\sigma_{Z_{max}}\}} ({\Bbb A}^1 - \{\delta_j\},{\Bbb A}^1 - \{\delta_j\}) &\times
(T_{X(\pi(N \cap L),\pi(\partial_-\tau_{Z'_{max}}))},X(\pi(N
\cap L),\pi(\partial_-\tau_{Z'_{max}})) \\
\cong \prod_{\{j;j \in S, v_j \in
\sigma_{Z_{max}}\}} ({\Bbb A}^1 - \{0\},{\Bbb A}^1 - \{0\}) &\times
(T_{X(\pi(N \cap L),\pi(\partial_-\tau_{Z'_{max}}))},X(\pi(N
\cap L),\pi(\partial_-\tau_{Z'_{max}})),
\\
\endalign$$
the situations $(\heartsuit_{\pi(\partial_-fan)})$ and
$(\heartsuit_{\pi(\partial_+fan)})$ would imply that the morphisms $h_{i-}$ and $h_{i+}$ over
these open sets are toric morphisms of the form
$$\align
h_{i-}:&\prod ({\Bbb A}^1 - \{0\},{\Bbb A}^1 - \{0\}) \times
(T_{X(\pi(N \cap
L),\pi(\partial_-\tau_{Z'_{max}}))},X(\pi(N \cap L),\pi(\partial_-\tau_{Z'_{max}})^{can})
\\
\rightarrow &\\
&\prod ({\Bbb A}^1 - \{0\},{\Bbb A}^1 - \{0\}) \times
(T_{X(\pi(N \cap L),\pi(\partial_-\tau_{Z'_{max}}))},X(\pi(N
\cap L),\pi(\partial_-\tau_{Z'_{max}}))
\\ h_{i+}:&\prod ({\Bbb A}^1 - \{0\},{\Bbb A}^1 - \{0\}) \times
(T_{X(\pi(N \cap
L),\pi(\partial_+\tau_{Z'_{max}}))},X(\pi(N \cap L),\pi(\partial_+\tau_{Z'_{max}})^{can})
\\
\rightarrow &\\
&\prod ({\Bbb A}^1 - \{0\},{\Bbb A}^1 - \{0\}) \times
(T_{X(\pi(N \cap L),\pi(\partial_+\tau_{Z'_{max}}))},X(\pi(N
\cap L),\pi(\partial_+\tau_{Z'_{max}}))
\\
\endalign$$
where the morphisms on the first factors are identities and where the morphisms on the second
factors are inducde by the refinements of the cones and that hence $h_{i-}$ and $h_{i+}$ are
toroidal.  Since we have the diagram of toric morphisms

$$\CD
X(\pi(N \cap L),\pi(\partial_-\tau_{Z'_{max}})^{can})
@.\overset{\varphi_i^{can}}\to{\dashrightarrow} @. X(\pi(N \cap
L),\pi(\partial_+\tau_{Z'_{max}})^{can}) \\ @VVV @.@VVV \\
X(\pi(N \cap L),\pi(\partial_-\tau_{Z'_{max}})) @.@. X(\pi(N
\cap L),\pi(\partial_+\tau_{Z'_{max}})) \\
\hskip.5in \searrow @.@. \swarrow \hskip.5in \\
@.X(\pi(N \cap L),\pi(\tau_{Z'_{max}})), @.\\
\endCD$$
we also conclude that $\varphi_i^{can}$ is V-toroidal.

\vskip.1in

Now we only discuss the verification of the situation $(\heartsuit_{\pi(\partial_-fan)})$, as
that of $(\heartsuit_{\pi(\partial_+fan)})$ is identical replacing the $-$ sign with the $+$
sign.

\vskip.1in

Note that, induced by the product structure
$$\sigma_{Z_{max}} = \langle v_j \rangle \oplus \tau_{{Z'}_{max,j}},$$
where the cone $\tau_{{Z'}_{max,j}}$ is contained in a hyper plane $H_j$ which also contains
$a \in N$ and where we have the exact sequence
$$0 \rightarrow H_j \cap N \rightarrow N \rightarrow {\Bbb Z} \rightarrow 0$$
with the image of $v_j \in N$ to be $\pm 1 \in {\Bbb Z}$, we have the product structure 
$$\pi(\partial_-\sigma_{Z_{max}}) = \langle
\pi(v_j)
\rangle \oplus \pi(\partial_-\tau_{{Z'}_{max,j}}).$$

Since
$$\pi(L) = \cap_{j \not\in S, v_j \in \sigma_{Z_{max}}} \pi(H_j),$$
it is straightforward to observe that in order to prove the
situation $(\heartsuit_{\pi(\partial_-fan)})$ it suffices to show the collection of the
situations $(\heartsuit_{j,\pi(\partial_-fan)})$ for all $j = 1, \cdot\cdot\cdot, l$:

$$(\heartsuit_{j,\pi(\partial_-fan)}) \left\{
\aligned
&\text{Provided\ }j \not\in S \text{\ and\ }v_j \in \sigma_{Z_{max}}, \\
&\pi(\partial_-{\sigma_{Z_{max}}})^{can} \text{has\ the\ product\ structure\ } \\
&\pi(\partial_-{\sigma_{Z_{max}}})^{can} = \langle
\pi(v_j)
\rangle \oplus \pi(\partial_-\tau_{{Z'}_{max,j}})^{can} \\
&\text{inherited\ from\ the\ product\ structure\ of\ }\pi(\partial_-\sigma_{Z_{max}}) \text{\ and\
}\pi(\sigma_{Z_{max}})\\ 
&\pi(\partial_-\sigma_{Z_{max}}) = \langle
\pi(v_j)
\rangle \oplus \pi(\partial_-\tau_{{Z'}_{max,j}}) \\
&\pi(\sigma_{Z_{max}}) = \langle
\pi(v_j)
\rangle \oplus \pi(\tau_{{Z'}_{max,j}}) \\
&\text{where\ }\pi(\partial_-\tau_{{Z'}_{max,j}})^{can}\text{\ is\ a\ refinement\ of\
}\pi(\partial_-\tau_{{Z'}_{max,j}}).\\ 
\endaligned
\right.$$

Now we prove the situation $(\heartsuit_{j,\pi(\partial_-fan)})$ $(j = 1, \cdot\cdot\cdot, l)$
in the form of the following proposition.

\proclaim{Proposition 4-1-4} Let $\pi(\partial_-{\sigma_{Z_{max}}})^{can}$ be the subfan of
$\Sigma_-{can}$ lying over the cone $\pi(\partial_-\sigma_{Z_{max}})$ of $\Sigma_-^{tor}$.  Let
$v_j \hskip.1in (j = 1, \cdot\cdot\cdot, l)$ be such primitive vector with $v_j \not\in
D^{tor}$ and $v_j \in \sigma_{Z_{max}}$.  Then
any new extremal ray $w \in \pi(\partial_-{\sigma_{Z_{max}}})^{can}$ to subdivide
$\pi(\partial_-\sigma_{Z_{max}})$ sits inside of the cone $\pi(\partial_-\tau_{Z'_{max,j}})$,
where $\pi(\partial_-\sigma_{Z_{max}})$ has the product structure
$$\pi(\partial_-\sigma_{Z_{max}}) = \langle \pi(v_j) \rangle \oplus
\pi(\partial_-\tau_{{Z'}_{max,j}})$$
induced from the product structure of $\sigma_{Z_{max}}$
$$\sigma_{Z_{max}} = \langle v_j \rangle \oplus \tau_{{Z'}_{max,j}}.$$
Therefore, we have the situation $(\heartsuit_{j,\pi(\partial_-fan)})$.
\endproclaim

\demo{proof}\enddemo We use the same notation as in Proposition 3-1-8 or Observation 3-2-1.

We take a ${\Bbb Z}$-basis $\{v_1, \cdot\cdot\cdot,
v_n\}$ of the lattice $N$ so that
$$\align
\sigma &= \langle v_1, \cdot\cdot\cdot, v_l\rangle \\
\sigma^{\vee} &= \langle v_1^*, \cdot\cdot\cdot v_l^*, \pm v_{l+1}^*,
\cdot\cdot\cdot, \pm v_n^*\rangle.\\
\endalign$$
Setting $z_i = z^{v_i^*}$, we have
$$X = \roman{Spec} K[z_1, \cdot\cdot\cdot, z_l,z_{l+1}^{\pm 1},
\cdot\cdot\cdot, z_n^{\pm 1}].$$
The $K^*$-action corresponds to a one-parameter subgroup
$$a = (\alpha_1, \cdot\cdot\cdot, \alpha_n) \in N,$$
i.e., the action of $t \in K^*$ is given by
$$t \cdot (z_1, \cdot\cdot\cdot, z_n) = (t^{\alpha_1}z_1, \cdot\cdot\cdot,
t^{\alpha_n}z_n).$$

Then by the proof of Theorem 3-2-8 the affine toric variety associated to a maximal cone
$\sigma_{Z,max} \in \sigma^{tor}$ has the form
$$\roman{Spec}\{K[v_j^* - m_1] \otimes K[\{y_h\}]\}$$
where the monomials $\{y_h\}$ are taken from $v_j^{\perp} \cap M$ and where in the exact sequence
$$0 \rightarrow v_j^{\perp} \cap M \rightarrow M \rightarrow {\Bbb Z} \rightarrow 0$$
the image of $v_j^* - m_1 \in M$ is $\pm 1 \in {\Bbb Z}$.

As the dual statement in terms of the cones, we obtain the product structure
$$\sigma_{Z,max} = \langle v_j \rangle \oplus \tau_{Z'_{max,j}}$$
where the cone $\tau_{Z'_{max,j}}$ is defined by the equation $(v_j^* - m_1, \cdot) = 0$. 
Therefore, in the induced product structure
$$\Sigma^{tor}_- \supset \partial_-\sigma_{Z,max} = \langle \pi(v_j) \rangle \oplus
\partial_-\tau_{Z'_{max,j}}$$
the cone $\partial_-\tau_{Z'_{max,j}}$ is also definde by the equation $(v_j^* - m_1, \cdot) =
0$.  (Note that since $(v_j^* - m_1,a) = 0$, we have $v_j^* - m_1 \in M \cap a^{\perp} =
Hom_{\Bbb Z}(\pi(N),{\Bbb Z})$.)

Therefore, in order to show that a new extremal ray $w \in \pi(\partial_-{\sigma_{Z_{max}}})^{can}$ to subdivide
$\pi(\partial_-\sigma_{Z_{max}})$ sits inside of
$\partial_-\tau_{Z'_{max,j}}$ it suffices to show
$$(v_j^* - m_1,w) = \roman{ord}_{E_w}(z_j \cdot z^{- m_1}) \leq 0,$$
where $\roman{ord}_{E_w}(r)$ denotes the order of the rational function $r$ on the divisor
(necessarily exceptional for $h_{i-}$) corresponding to the ray $w$.

Now consider the family of automorphisms $\theta_c$ parametrized by $c \in K$ of the ring $K[z_1, \cdot\cdot\cdot, z_l,z_{l+1}^{\pm 1},
\cdot\cdot\cdot, z_n^{\pm 1}]$ (or
equivalently the action of $K = {\Bbb A}^1$ on $\roman{Spec} K[z_1, \cdot\cdot\cdot, z_l,z_{l+1}^{\pm 1},
\cdot\cdot\cdot, z_n^{\pm 1}]$) defined by
$$\align
\theta_c^*(z_j) &= z_j + c \cdot z^{m_1}\\
\theta_c^*(z_k) &= z_k \text{\ for\ }k \neq j.\\
\endalign$$
Observe that the these automorphisms preserve the torific ideals because $z_j$ and $m_1$ has the
same character.  It follows that $\theta_c$ induces a family of automorphisms of $W_{i-}^{tor}$
mapping the divisor $E_{i-}^{tor}$ to itself.  Therefore, the automorphisms
also preserve $I_{i-}$.  (See Lemma 4-1-2.)

Therefore, by the properties $(\spadesuit^{res}-1$) and $(\spadesuit^{can}-1)$ of the canonical
resolution and canonical principaliztion, we conclude that $\theta_c$ induces an automorphism of
$W_{i-}^{can}$ which preserves the exceptional divisor $E_w$.  (Note that the set of
exceptional divisors for $h_{i-}$ are discrete, while the family is parametrized by a
continuous family ${\Bbb A}^1$.  See Remark 4-1-1.)

Then, denoting the induced automorphism by the same $\theta_c$ by abuse of notation, we have
$$\align
\roman{ord}_{E_w}(z_1 \cdot z^{- m_1}) &= \roman{ord}_{E_w}(\theta_c^*(z_1 \cdot z^{- m_1}) \\
&= \roman{ord}_{E_w}((z_1 + c \cdot z^{m_1}) \cdot z^{- m_1}) \\
&= \roman{ord}_{E_w}(z_1 \cdot z^{- m_1} + c).\\
\endalign$$
Since the above equality holds for any $c \in K$, we conclude
$$\roman{ord}_{E_w}(z_1 \cdot z^{- m_1}) \leq 0\ (\text{and\ hence\ }= 0)$$
as required.

This completes the proof of Proposition 4-1-4.

\vskip.1in

\S 4-2. Conclusion of the proof for the Weak Factorization Theorem

\vskip.1in

Since $X_1 = W_1$ and $W_s = X_2$ are nonsingular, we have 
$$X_1 = W_1 \overset{\overset{p_{1-}}\to{\sim}}\to{\leftarrow} W_1^{res} \text{\ and\
}W_l^{res} \overset{\overset{p_{s-1,-}}\to{\sim}}\to{\rightarrow} W_s = X_2.$$
For each $i = 1, \cdot\cdot\cdot, s-1$ we have constructed the commutative diagram
$$\CD
W_i^{res} @<{p_{i-}}<< W_{i-}^{can} @.
\hskip.2in \overset{\varphi_i^{can}}\to{\dashrightarrow} \hskip.2in @. W_{i+}^{can}
@>{p_{i+}}>> W_{i+1}^{res} \\
@V{r_{i-}}VV @V{h_{i-}}VV @. @VV{h_{i+}}V @VV{r_{i+}}V \\
W_i @<{f_{i-}}<< W_{i-}^{tor} @.
\hskip.2in \overset{\varphi_i^{tor}}\to{\dashrightarrow} \hskip.2in @. W_{i+}^{tor}
@>{f_{i+}}>> W_{i+1} \\
\endCD$$
where 

\ \ (i) $W_i^{res},W_{i-}^{can},W_{i+}^{can}$ and $W_{i+1}^{res}$ are all nonsingular, 

\ \ (ii) the
morphisms $p_{i-}$ and $p_{i+}$ are sequences of blowups with smooth centers, and

\ \ (iii) $W_{i-}^{can}$ and $W_{i+}^{can}$ have toroidal structures so that
$$\varphi_i^{can}:(U_{W_{i-}^{can}},W_{i-}^{can}) \dashrightarrow
(U_{W_{i+}^{can}},W_{i+}^{can})$$
is V-toroidal.

\proclaim{Theorem 4-2-1 (Strong Factorization Theorem for Toroidal Birational Maps)} Let
$$\varphi:(U_{W_1},W_1) \dashrightarrow (U_{W_2},W_2)$$
be a proper and toroidal biraional map between nonsingular toroidal embeddings.  Then
$\varphi$ can be factored into a sequence of toroidal blowups immediately followed by toroidal
blowdowns with smooth centers.
\endproclaim

\demo{Proof}\enddemo If $\varphi:(U_{W_1},W_1) \rightarrow (U_{W_2},W_2)$ is a proper and
toroidal birational morphism between nonsingular toroidal
embeddings without self-intersections (which is necessarily
allowable in the sense of KKMS [1] (cf. Abramovich-Karu [1] or Chapter 2)), then the
theorem is exactly what is proved in Abramovich-Matsuki-Rashid [1]. 

In general, by definition of a proper birational map $\varphi$ being toroidal there exists a
toroidal embedding $(U_V,V)$ which dominates $(U_{W_1},W_1)$ and $(U_{W_2},W_2)$ by proper and
toroidal birational morphisms
$$(U_{W_1},W_1) \overset{\varphi_1}\to{\leftarrow} (U_V,V) \overset{\varphi_2}\to{\rightarrow}
(U_{W_2},W_2).$$
(Note that in case $\varphi$ is $V$-toroidal we can take $(U_V,V)$ to consist of the
normalization $V$ of the graph $\Gamma_{\varphi} \subset W_1 \times W_2$ and the open set $U_V =
p_1^{-1}(U_{W_1}) = p_2^{-1}(U_{W_2})$ where $p_1:V \rightarrow W_1$ and $p_2:V \rightarrow W_2$
are the projections.)

By Lemma 1-2-6 there exists a sequence of toroidal blowups of $(U_{W_1},W_1)$ with smooth
centers so that the result is a toroidal embedding without self-intersection.  Since all the
blowups are toroidal and since
$\varphi_1$ is also toroidal, we can take the corresponding sequence of toroidal blowups of
$(U_V,V)$ so that the result dominates the blowup of $(U_{W_1},W_1)$.  We can apply the same
procedure to $\varphi_2$.  By eliminating the self-intersection of and resolving the
singularities of $(U_V,V)$ by a sequence of toroidal blowups, we finally conclude that we may
assume that $(U_{W_1},W_1), (U_V,V)$ and
$(U_{W_2},W_2)$ are all nonsingular toroidal embeddings without self-intersections.  (Note that
in case $\varphi$ is $V$-toroidal, then $\varphi_1$ and $\varphi_2$ are automatically allowable by
constructruction, without referring to the result of Abramovich-Karu [1] (cf. Lemma 1-2-7).)

Now we are in the situation described as above to apply the strong factorization theorem to factor
$\varphi_1$ as
$(U_{W_1},W_1) \overset{\mu_1}\to{\leftarrow} (U_{\tilde V},{\tilde V})
\overset{\mu_V}\to{\rightarrow} (U_V,V)$ where $\mu_1$ and $\mu_V$ are toroidal blowups with
smooth centers.  Now we only have to apply the strong factorization theorem to $\varphi_2 \circ
\mu_V$ to complete the proof for Theorem 4-2-1. 

\vskip.1in

We apply Theorem 4-2-1 to the proper and V-toroidal birational map $\varphi_i^{can}$ to finish
the proof of the weak factorization theorem.

\vskip.1in

($\bold{Preservation\ of\ the\ open\ set\ where\ \phi\ is\ an\ isomorphism}$)

\vskip.1in

Suppose $\phi:X_1 \dashrightarrow X_2$ induces an isomorphism on an open set $\phi:X_1 \supset
U \overset{\sim}\to{\rightarrow} U \subset X_2$.  Then

Step 1: In the process of eliminating the points of indeterminacy to be reduced to the case
where $\phi$ is a projective birational morphism, all the centers of blowups are taken outside of
$U$.

Step 2: In the process of constructing the cobordism $B_{\phi}(X_1,X_2)$ out of the product
$X_2 \times {\Bbb P}^1$, all the centers of blowups are taken outside of $U \times {\Bbb P}^1$.

Step 3: In the process of torification, since the torific ideals are trivial over $U \times
{\Bbb P}^1 \cap B_{a_i}$ for all $i$, all the modifications are made outside of $U \times
{\Bbb P}^1 \cap B_{a_i}$.

\vskip.1in

Therefore, through Steps 1,2 and 3, the open set $U$ remains untouched inside of all
the $W_i$ and $U_{W_{i \pm}^{tor}} \subset W_{i \pm}^{tor}$.

\vskip.1in

Step 4: In the process of obtaining the $W_i^{res}$ or $W_{i+1}^{res}$ from $W_{i \pm}^{tor}$ all
the centers of blowups are taken over the singular locus and hence outside of $U$.  

In the process of obtaining the $W_{i \pm}^{can}$ from $W_i^{res}$ or $W_{i+1}^{res}$
all the centers are taken over the quotient (by the $K^*$-action) of the locus where the torific
ideals are not trivial and hence outside of $U$.

Step 5: Finally in the process of factoring the toroidal birational maps
$\varphi_i^{can}:(U_{W_{i-}^{can}},W_{i-}^{can}) \dashrightarrow (U_{W_{i+}^{can}},W_{i+}^{can})$
all the centers of blowups are taken within the boundary divisors and hence outside of $U$.

\vskip.1in

Thus we conclude that all through the construction of our factorization the open set $U$ is
preserved without being modified.

\vskip.1in

($\bold{Preservation\ of\ the\ projectivity}$)

\vskip.1in

Step 1: In the process of eliminating the points of indeterminacy to be reduced to the case
where $\phi$ is a projective morphism via Lemma 1-4-2, by replacing $\phi:X_1 \rightarrow X_2$
with $\phi':X_1' \rightarrow X_2'$, $X_1'$ as well as all the intermediate varieties in the
sequence of blowups starting from $X_1$ are projective over $X_1$, while $X_2'$ as well as all the
intermediate varieties in the sequence of blowups starting from $X_2$ are projective over $X_2.$

Step 2: In the proces of constructing the birational cobordism
$B_{\phi'}(X_1',X_2')$, all the $W_i$ are projective over $X_2'$ via the Geometric Invariant
Theory interpretation in Chapter 2.  

Steps 3, 4: All the $W_{i \pm}^{tor}, W_i^{res}$ and $W_{i \pm}^{can}$ are also
projective over $X_2'$, as they are obtained from
$W_i$ by sequences of blowing ups.  

Step 4: Now by the use of the STRONG factorization theorem applied
to the toroidal birational maps $\varphi_i^{can}:(U_{W_{i-}^{can}},W_{i-}^{can}) \dashrightarrow
(U_{W_{i-}^{can}},W_{i+}^{can})$ all the varieties between $W_{i-}^{can}$ and $W_{i+}^{can}$ are
also projective over $X_2'$ in the factorization.

\vskip.1in

Thus we conclude that in the factorization we construct there is an index $i_o$ such that all the
intermediate varieties
$V_i$ for $i \leq i_o$ are prjective over $X_1$ and that $V_i$ for $i_o \leq i$ are projective
over $X_2$.  In particular, if both $X_1$ and $X_2$ are projective, they are all projective.

\vskip.1in

This completes the proof of the Weak Factorization Theorem and Chapter 4. 

\newpage

$$\bold{CHAPTER\ 5.\ GENERALIZATIONS}$$

\vskip.2in

In this chapter we provide several generalizations of the weak factorization
theorem, specifying the modifications we have to make to the arguments in Steps 1 through 5 of
the strategy of the proof (cf. \S 0-1. Outline of the strategy) according to the
different purposes.

In \S 5-1, we establish the weak factorization theorem for bimeromorphic maps between compact
complex manifolds.  Except for a couple of places where we have to avoid the use of theorems only
applicable in the algebraic category and replace them with some arguments valid in the analytic
category, the proof goes with little change and in fact is even simpler.  For example, instead of
discussing with Luna's locally toric charts where we have to require the property of being
strongly \'etale in the algebraic category, we can take ${\Bbb C}^*$-equivariant analytic
isomorphisms as local charts where the induced analytic isomorphisms for the quotients are
automatic.

In \S 5-2, we consider the case where a birational map $\phi:X_1 \dashrightarrow X_2$
is equivariant with respect to the actions of a group on $X_1$ and $X_2$ and establish
the equivariant weak factorization theorem.  (Remark that the action does NOT have to be over
the base field $K$.)  We note here that in order to preserve the equivariance we may have to blow
up several smooth irreducible centers simultaneously in the equivariant factorization.  Since most
of the constructions in our proof are canonical, the equivariance comes almost free under the
group action in the main steps of the strategy and the proof is reduced to the equivariant
version of the factorization of toroidal birational maps, for which we give a proof utilizing an
idea of Abramovich-Wang [1] and hence generalizing the result of Morelli [1][2]
Abramovich-Matsuki-Rashid [1]. 

In \S 5-3, we establish the weak factorization theorem over a base field $K$ which may NOT be
algebraically closed.  By taking the algebraic closure $\overline{K}$ of $K$ and considering
the action of the Galois group $G = Gal(\overline{K}/K)$ after base change, it is almost
a direct corollary of the arguments for the equivariant case, except for the technical
difficulty of dealing with the possibility that after base change $X_1 \times_{\roman{Spec}K}
\roman{Spec}\overline{K}$ and $X_2 \times_{\roman{Spec}K}
\roman{Spec}\overline{K}$ may not stay irreducible but split into several smooth irreducible
components.
 
In \S 5-4, we establish the weak factorization theorem in the logarithmic category, factoring a
proper birational map between nonsingular toroidal embeddings into blowups and blowdowns with
smooth and ADMISSIBLE centers.  This has an application to the study of the behavior of the Hodge
structures under birational transformations, e.g., it provides an easy proof for a theorem of
Batyrev [1] claiming the invariance of the Betti numbers among birational nonsingular minimal
models.

In \S 5-5, we discuss the toroidalization and strong factorization conjectures in
the vicinity of the circle of ideas involving the factorization problem, resolution of
singularities, semi-stable reductions and the log category of Kato [1].

\vskip.2in

\S 5-1. Factorization of bimeromorphic maps

\vskip.1in

\proclaim{Theorem 5-1-1 (Weak Factorization Theorem of Bimeromorphic Maps)} Let $\phi:X_1
\dashrightarrow X_2$ be a bimeromorphic map between compact complex manifolds.  Let $X_1 \supset
U \subset X_2$ be a common open subset over which $\phi$ is an isomorphism.  Then $\phi$ can be factored into a
sequence of blowups and blowdowns with smooth irreducible centers disjoint from $U$.  That is to
say, there exists a sequence of bimeromorphic maps between compact complex manifolds
$$X_1 = V_1 \overset{\psi_1}\to{\dashrightarrow} V_2 \overset{\psi_2}\to{\dashrightarrow}
\cdot\cdot\cdot \overset{\psi_{i-1}}\to{\rightarrow} V_i \overset{\psi_i}\to{\rightarrow}
V_{i+1} \overset{\psi_{i+1}}\to{\rightarrow} \cdot\cdot\cdot
\overset{\psi_{l-2}}\to{\dashrightarrow} V_{l-1}
\overset{\psi_{l-1}}\to{\rightarrow} V_l = X_2$$ where

\ \ (i) $\phi = \psi_{l-1} \circ \psi_{l-2} \circ \cdot\cdot\cdot \circ \psi_2 \circ \psi_1$,

\ \ (ii) $\psi_i$ are isomorphisms over $U$, and

\ \ (iii) either $\psi_i:V_i \dashrightarrow V_{i+1}$ or $\psi_i^{-1}:V_{i+1} \dashrightarrow V_i$
is a morphism obtained by blowing up a smooth irreducible center disjoint from $U$.

Moreover, if both $X_1$ and $X_2$ are projective, then we can choose a factorization so that all
the intermediate complex manifolds $V_i$ are projective.
\endproclaim

\demo{Proof}\enddemo  We only specify the modifications we have to make at each step of the
strategy for the proof described in Chapter 0. Introduction.  (Note that if $X_1$ and $X_2$ are
projective, then both are algebraic and hence so is $\phi$.  Thus the ``Moreover" part is already
proved but is stated explicitly again only for the sake of completeness.)

\vskip.1in

Step 1. Elimination of points of indeterminacy

\vskip.1in

Thanks to Hironaka [3], Lemma 1-4-1 holds without any change in the category of compact complex
manifolds.  Note that by construction there exists a line bundle of the form ${\Cal A} = {\Cal
O}_{X_1'}(\Sigma -a_iE_i)$, where the $E_i$ are $g_1$-exceptional divisors with $a_i > 0$, such
that ${\Cal A}$ is relatively ample for $\phi':X_1' \rightarrow X_2'$, i.e., we have a
projective embedding over $X_2'$
$$\phi':X_1' \hookrightarrow {\Bbb P}(\phi'_*({\Cal A}^{\otimes l})) \rightarrow X_2' \text{\
for\ some\ sufficiently\ large\ }l \in {\Bbb N}.$$
Note that $J = \phi'_*({\Cal A}^{\otimes l}) = \phi'_*{\Cal O}_{X_1'}(l \cdot (\Sigma - a_iE_i))
\subset \phi'_*{\Cal O}_{X_1'} = {\Cal O}_{X_2'}$ is an ideal sheaf of ${\Cal O}_{X_2'}$, which
coincides with ${\Cal O}_{X_2'}$ outside of $U$.

Therefore, by replacing $\phi:X_1 \dashrightarrow X_2$ with $\phi':X_1' \rightarrow X_2'$, we may
assume as in the algebraic case that $\phi:X_1 \rightarrow X_2$ is a projective morphism
which is the blowup of
$X_2$ along an ideal sheaf $J \subset {\Cal O}_{X_2}$ with the support of ${\Cal O}_{X_2}/J$ being
disjoint from $U$.

\vskip.1in

Step 2. Construction of a birational cobordism

\vskip.1in

Theorem 2-2-2 holds for a projective bimeromorphic morphism $\phi:X_1 \rightarrow X_2$ described
as in Step 1.  Though the proof goes almost without any change, we show the existence of such a
coherent sheaf ${\Cal E}$ as described in the assertion (ii) with the splitting ${\Cal E} =
\oplus_{b \in {\Bbb Z}}{\Cal E}_b$ into the eigenspaces inductively as follows (since we would
like to avoid the complete reducibility argument for a coherent sheaf ${\Cal E}$ under a
$K^*$-action, which is directly applicable only to an algebraic coherent sheaf under an
algebraic
$K^*$-action):

Recall that we obtain $\overline{B} = W_l$ by a sequence of successive blowups with ${\Bbb C}^*$
-invariant centers:
$$\overline{B} = W_l \rightarrow \cdot\cdot\cdot \rightarrow W_{i+1}
\overset{\mu_i}\to{\rightarrow} W_i \rightarrow \cdot\cdot\cdot \rightarrow W_1 = W
\rightarrow W_0 = X_2 \times {\Bbb P}^1.$$ 
We denote the induced projection onto
$X_2$ by
$$\tau_i:W_i \rightarrow X_2.$$

We start with the product $W_0 = X_2 \times {\Bbb P}^1$.  Here we obviously have a relatively
(very) ample line bundle ${\Cal L}_0 = p_2^*{\Cal O}_{{\Bbb P}^1}$ and a coherent sheaf ${\Cal
E}_0 = (\tau_0)_*{\Cal L}_0$ with the splitting
$$\align
{\Cal E} &= (\tau_0)_*{\Cal L}_0 = (p_1)_*{\Cal L}_0 \\
&= \oplus_{b \in {\Bbb Z}}({\Cal E}_0)_b = p_1^*{\Cal O}_{X_2} \cdot T_0 \oplus p_1^*{\Cal
O}_{X_2} \cdot T_1 \\
\endalign$$
where $T_0$ and $T_1$ are the homogeneous coordinates of ${\Bbb P}^1$, on which $t \in K^*$ acts
as $t^*(T_0) = T_0$ and $t^*(T_1) = t \cdot T_1$ by definition.

Suppose $W_i$ is embedded equivariantly into ${\Bbb P}({\Cal E}_i)$ where ${\Cal E}_i =
(\tau_i)_*{\Cal L}_i$ is a coherent sheaf with the splitting ${\Cal E}_i = \oplus_{b \in {\Bbb
N}}({\Cal E}_i)_b$ into the eigenspaces for some relatively (very) ample line bundle ${\Cal
L}_i$ with a ${\Bbb C}^*$-action.

Since the center of the blowup $W_{i+1} \overset{\mu_i}\to{\rightarrow} W_i$ is
${\Bbb C}^*$-invariant (See the property
$(\spadesuit^{res}-1)$ of the canonical resolution.), we conclude the following: The defining
ideal sheaf $I_i$ of the center has the property that there exists a finite open covering
$\{U\}$ of
$X_2$ and a sufficiently large integer $k_i
\in {\Bbb N}$ such that
$$(\tau_i)^*(\tau_i)_*\{I_i \otimes {\Cal L}_i^{\otimes k_i}\} \rightarrow I_i \otimes {\Cal
L}_i^{\otimes k_i}$$
is surjective and that over each open set $U$ of the covering the sheaf $(\tau_i)^*(\tau_i)_*\{I_i
\otimes {\Cal L}_i^{\otimes k_i}\}$ is generated by a finite number of sum of products, each of
which has the same character as the others, of elements in the eigenspaces $\Gamma(U,({\Cal
E}_i)_b)$ with coefficients in
$\Gamma(U,{\Cal O}_{X_2})$.  It follows immediately that if we set
$$\align
{\Cal L}_{i+1} &= \mu_i^*I_i \cdot {\Cal O}_{W_{i+1}} \otimes \mu_i^*{\Cal L}_i^{\otimes
k_i}
\\
{\Cal E}_{i+1} &= (\tau_{i+1})_*{\Cal L}_{i+1} \\
\endalign$$
then ${\Cal L}_{i+1}$ has a natural $K^*$-action with the splitting ${\Cal E}_{i+1} = \oplus_{b
\in {\Bbb N}}({\Cal E}_{i+1})_b$ and the equivariant embedding $W_{i+1} \hookrightarrow {\Bbb
P}({\Cal E}_{i+1})$.

The rest of the proof for Theorem 2-2-2 is identical to the one in the algebraic setting.

\vskip.1in

The GIT interpretation of the cobordism $B$, obtained from the compactified birational
(bimeromorphic) cobordism $\overline{B}$ projective over $X_2$ constructed as above, holds without
any change in the analytic setting, leading to the definition of the quasi-elementary pieces
$B_{a_i}$. 

\vskip.1in

In order to reduce the analysis in general to the one on toric varieties, we take Luna's locally
toric (toroidal) charts in the algebraic setting.  In order to show the existence of these
charts, we used some theorems of Sumihiro (Equivariant completion theorem $\&$ Covering by
$K^*$-invariant affine spaces), Luna's Fundamental Lemma and (or) Luna's \'Etale Slice Theorem,
all of which are directly applicable only in the algebraic setting (cf. Proposition 1-3-4 and
Lemma 2-4-4).  We repalce Luna's locally toric charts by ${\Bbb C}^*$-equivariant analytic
isomormorphisms via the following lemma.

\proclaim{Lemma 5-1-2} Let $B_{a_i}$ be a quasi-elementary cobordism associated to the
birational (bimeromorphic) cobordism $B$ obtained from the compactified birational
(bimeromorphic) cobordism $\overline{B}$ projective over $X_2$ constructed as above.  Then for
every point $p \in B_{a_i} - F_{a_i}^*$, where $F_{a_i} = B_{a_i}^{{\Bbb C}^*}$ is the fixed
point set and
$F_{a_i}^*$ is the set as defined in Notation 2-1-3, there exists a ${\Bbb C}^*$-equivariant
analytic isomorphism
$$X_p \supset V_p \overset{\overset{\eta_p}\to{\sim}}\to{\leftarrow} U_p \subset B_{a_i}$$
where $X_p$ is a nonsingular affine toric variety with ${\Bbb C}^*$-acting as a one-parameter
subgroup and where $U_p \subset B_{a_i}$ (resp. $V_p \subset X_p$) is a ${\Bbb C}^*$-invariant
open subset (with respect to the usual topology) satisfying the following condition $(\star)$:

$(\star)$ If a ${\Bbb C}^*$-orbit $O(q)$ lies in $U_p$ (resp. $V_p$), then its closure
$\overline{O(q)}$ in $B_{a_i}$ (resp. in $X_p$) also lies in $U_p$ (resp. in $V_p$).

(As a consequence, we obtain
$$X_p//{\Bbb C}^* \supset V_p//{\Bbb C}^* \overset{\overset{\eta_p//{\Bbb
C}^*}\to{\sim}}\to{\leftarrow} U_p//{\Bbb C}^* \subset B_{a_i}//K^*$$
providing not only ${\Bbb C}^*$-equivarinat analytic local isomorphisms of $B_{a_i}$ to
(nonsingular) toric varieties bot also analytic local isomorphisms of the quotient $B_{a_i}//{\Bbb
C}^*$ to (not necessarily nonsingular) toric varieties.)
\endproclaim

\demo{Proof}\enddemo It follows from Theorem 2-2-2 that, locally over $X_2$, $\overline{B}$ can be
embedded equivariantly into $X_2 \times {\Bbb P}^N$, i.e., there exists a open neighborhood
$U_{X_2}$ of $\tau_l(p)$ such that
$$\tau_l^{-1}(U_{X_2}) \hookrightarrow U_{X_2} \times {\Bbb P}^N,$$
where ${\Bbb C}^*$ acts only on the second
factor ${\Bbb P}^N$ with $t \in {\Bbb C}^*$ acting on the homogeneous coordinates $T_0,
\cdot\cdot\cdot, T_N$ of ${\Bbb P}^N$ by characters 
$$t^*(T_j) = t^{b_j} \cdot T_j.$$

Case: $p \in B_{a_i} - F_{a_i}^*$ is a fixed point, i.e., $p \in F_{a_i}$.

\vskip.1in

In this case, $T_j(p) = 0$ if $b_j \neq a_i$.  Take $j_o$ with $T_{j_o}(p) \neq 0$ and $b_{j_0}
= a_i$.  Since $p \in B_{a_i} \cap \tau_l^{-1}(U_{X_2}) \subset U_{X_2} \times {\Bbb P}^N$ is a
nonsingular point, we may choose a regular coordinate system at $p$ consisting of a regular
coordinate system at
$\tau_l(p)$ of $U_{X_2}$ and $\{T_j/T_{j_0};j \in J\}$ for some subset $J \subset \{0, 1,
\cdot\cdot\cdot, N\}$.  Consider the ${\Bbb C}^*$-equivariant morphism
$$\eta:B_{a_i} \cap \tau_l^{-1}(U_{X_2}) \cap \{T_{j_o} \neq 0\} (\hookrightarrow U_{X_2} \times
{\Bbb P}^N \cap \{T_{j_o} \neq 0\}) \rightarrow {\Bbb A}^{\dim X_2 + \# J} = X_p$$ 
where ${\Bbb A}^{\dim X_2 + \# J}$ is the nonsingular affine toric variety $X_p$ with affine
coordinates corresponding to the regular coordinate system at $p$ choesn as above and the action
of
$t
\in {\Bbb C}^*$ is given by
$$t^*(T_j/T_{j_0}) = t^{b_j - b_{j_o}} \cdot T_j/T_{j_0},$$  
leaving others ${\Bbb C}^*$-invariant.  By construction, there exist open neighborhoods
$$\align
p &\in U'''_p \subset B_{a_i} \cap
\tau_l^{-1}(U_{X_2}) \cap \{T_{j_o} \neq 0\} \\
\eta(p) &\in V'''_p \subset X_p \\
\endalign$$
such that $\eta$ induces an analytic isomorphism between them
$$\eta|_{U'''_p}:U'''_p \overset{\sim}\to{\rightarrow} V'''_p.$$ 
Let ${\Bbb C}^*_1 = \{z \in {\Bbb C}^*;|z| = 1\}$ be the unit circle in ${\Bbb C}^*$.  Since $p$
is a fixed point for the action of ${\Bbb C}^*$, we may choose a small open neighborhood (with
respect to the usual topology) $p \in U''_p \subset U'''_p$ such that $t(q) \in U'''_p$ for all $q
\in U''_p$ and $t \in {\Bbb C}^*_1$.  This implies that
$$Stab(q) = Stab(\eta(q)) \text{\ for\ all\ }q \in U''_p.$$
Let
$$\align
U'_p &= \cup_{t \in {\Bbb C}^*}t(U''_p) \\
V'_p &= \cup_{t \in {\Bbb C}^*}t(\eta(U''_p)) \\
\endalign$$
be the two ${\Bbb C}^*$-invariant neighborhood of $p$ and $\eta(p)$, respectively.  From the
construction it follows that
$$\eta|_{U'_p}:U'_p \overset{\sim}\to{\rightarrow} V'_p$$
is a ${\Bbb C}^*$-equivariant analytic isomorphism.

Let $J' = \{j \in J;b_j - b_{j_o} \neq 0\}$ be the subset of $J$.  Then $F_{X_p}$ has the
description
$$F_{X_p} = \{x \in X_p;T_j/T_{j_o}(x) = 0 \text{\ for\ }j \in J'\}.$$
Let $\sigma:X_p = {\Bbb A}^{\dim X_2 + \# J} \rightarrow {\Bbb A}^{\dim X_2 + \# J'}$
be the obvious projection given by $\prod_{j \in J} T_j/T_{j_o}(x) \mapsto
\prod_{j \in J'}T_j/T_{j_o}(x)$, leaving the rest of the coodinates coming from $U_{X_2}$
unchanged.  Take a small open neighborhood
$W
\subset {\Bbb A}^{\dim X_2 + \# J'}$ of
$\sigma \circ \eta (p)$ so that $\sigma^{-1}(W) \cap F_{X_p} \subset V'_p$.  It follows that
$$\align
U_p &= U'_p \cap \eta^{-1}(V'_p \cap \sigma^{-1}(W)) \subset B_{a_i}\\
V_p &= V'_p \cap \sigma^{-1}(W) \subset X_p\\
\endalign$$
are ${\Bbb C}^*$-invariant open subsets satisfying the condition $(\star)$ with a ${\Bbb
C}^*$-equivariant analytic isomorphism
$$\eta_p:U_p \rightarrow V_p.$$

\vskip.1in

Case: $p \in B_{a_i} - F_{B_{a_i}}^*$ is NOT a fixed point, i.e., $p \not\in F_{a_i}$.

\vskip.1in

Take a limit point $q = \lim_{t \rightarrow 0}t(p) \text{\ or\ }\lim_{t \rightarrow \infty}t(p)
\in \overline{B}$ so that $q \in F_{\overline{B}}$ is a fixed point of $\overline{B}$.  By the
same argument as in the previous case, we can find a ${\Bbb C}^*$-equivariant analytic
isomorphism
$$\eta_q:U'_q \rightarrow V'_q$$
where $U'_q \subset \overline{B}$ and $V'_q \subset X_q$ are ${\Bbb C}^*$-invariant open
neighborhoods of $q \in \overline{B}$ and $\eta_q(q) \in X_q$, respectively, and where $X_q$ is a
nonsingular affine toric variety with ${\Bbb C}^*$ acting as a one-parameter subgroup.

Let $p_2:U_{X_2} \times {\Bbb P}^N \rightarrow {\Bbb P}^N$ be the second projection.  Since $p
\in B_{a_i} - F_{B_{a_i}}^*$ and $p \not\in F_{B_{a_i}}$, we conclude $p_2(p) \in ({\Bbb
P}^N)_{a_i} - F_{({\Bbb P}^N)_{a_i}}^*$ and $p_2(p) \not\in F_{({\Bbb P}^N)_{a_i}}$, noting that
$$\align
\tau_l^{-1}(U_{X_2}) \cap B_{a_i} &= \tau_l^{-1}(U_{X_2}) \cap \overline{B} \cap
p_2^{-1}(({\Bbb P}^N)_{a_i}) \text{\ and\ }\\
\tau_l^{-1}(U_{X_2}) \cap F_{B_{a_i}} &= \tau_l^{-1}(U_{X_2})
\cap B_{a_i} \cap p_2^{-1}(F_{({\Bbb P}^N)_{a_i}}).\\
\endalign$$

In the proof of Proposition 1-3-5, we already saw that there exists an affine ${\Bbb
C}^*$-invariant open neighborhood
$p \in U_{({\Bbb P}^N)_{a_i}} \subset ({\Bbb P}^N)_{a_i}$ such that $U_{({\Bbb P}^N)_{a_i}}$
satisfies the condition $(\star)$ and that $U_{({\Bbb P}^N)_{a_i}}$ has no fixed points.  

Note also that since $\eta_q(p)$ is not a fixed point in $X_q$, there exists a nonsingular
affine toric open subvariety $X_p \subset X_q$ such that $\eta_q(p) \in X_p$.

It
follows that
$$\align
U_p &= U_p' \cap p_2^{-1}(U_{({\Bbb P}^N)_{a_i}}) \cap \eta_q^{-1}(X_p) \subset B_{a_i}\\
V_p &= \eta_q(U_p) \subset X_p \\
\endalign$$ 
are ${\Bbb C}^*$-invariant open subsets satisfying the condition $(\star)$ with a ${\Bbb
C}^*$-equivariant analytic isomorphism $\eta_p =
\eta_q|_{U_p}:U_p
\overset{\sim}\to{\rightarrow} V_p$.

This completes the proof of Lemma 5-1-2.

\vskip.1in

It is worthwhile noting that, in the notation of the proof of Lemma 5-1-2,
$\tau_l^{-1}(U_{X_2})
\cap B_{a_i}//{\Bbb C}^*
\hookrightarrow U_{X_2} \times ({\Bbb P}^N)_{a_i}//{\Bbb C}^*$ is embedded as a closed analytic
subvariety.  In fact, let ${\Cal I}$ be the ideal sheaf defining $\tau_l^{-1}(U_{X_2})$ inside of
$U_{X_2} \times {\Bbb P}^N$.  Let ${\Cal O}_{{\Bbb P}^N}(1)$ be the very ample sheaf on ${\Bbb
P}^N$.  For sufficiently large $l \in {\Bbb N}$, we have an exact sequence
$$\align
0 &\rightarrow {p_1}_*\{{\Cal I} \otimes p_2^*{\Cal O}_{{\Bbb P}^N}(l)\} \\
&\rightarrow
(p_1)_*\{p_2^*{\Cal O}_{{\Bbb P}^N}(l)\} = {\Cal O}_{U_{X_2}} \otimes H^0({\Bbb P}^N,{\Cal
O}_{{\Bbb P}^N}(l)) \\
&\rightarrow {p_1}_*\{{\Cal O}_{\tau_l^{-1}(U_{X_2})}\} \rightarrow 0\\
\endalign$$
where
$$\align
{p_1}_*\{{\Cal I} \otimes p_2^*{\Cal O}_{{\Bbb P}^N}(l)\} \otimes {\Bbb C}(y) &\rightarrow
H^0(p_2^{-1}(y),{\Cal I} \otimes p_2^*{\Cal O}_{{\Bbb P}^N}(l)|_{p_2^{-1}(y)}) \\
(p_1)_*\{p_2^*{\Cal O}_{{\Bbb P}^N}(l)\} \otimes {\Bbb C}(y) &\rightarrow
H^0(p_2^{-1}(y), p_2^*{\Cal O}_{{\Bbb P}^N}(l)|_{p_2^{-1}(y)}) \cong H^0({\Bbb P}^N,{\Cal
O}_{{\Bbb P}^N}(l))
\\
\endalign$$
are isomorphisms for any $y \in U_{X_2}$.  It is easy to see inductively as before that the
construction can be carried out so that ${p_1}_*\{{\Cal I} \otimes p_2^*{\Cal O}_{{\Bbb
P}^N}(l)\}$ splits into the eigenspaces with respect to the action of ${\Bbb C}^*$.  We may also
assume that $H^0({\Bbb P}^N,{\Cal O}_{{\Bbb P}^N}(l))^{{\Bbb C}^*}$ generates the homogeneous
invariant ring $\oplus_{m \geq 0}H^0({\Bbb P}^N,{\Cal O}_{{\Bbb P}^N}(ml))^{{\Bbb C}^*}$.  Now
by Fogarty-Mumford-Kirwan [1] Chapter 1 \S 2, in each fiber $p_2^{-1}(y) \cong {\Bbb P}^N$,
the invariant homogeneous polynomials $H^0(p_2^{-1}(y),{\Cal I} \otimes p_2^*{\Cal O}_{{\Bbb
P}^N}(l)|_{p_2^{-1}(y)})^{{\Bbb C}^*}$ defines the image of $B_{a_i}$ under the quotient map
$U_{X_2} \times ({\Bbb P}^N)_{a_i} \rightarrow U_{X_2} \times ({\Bbb P}^N)_{a_i}//{\Bbb C}^*$. 
That is to say, the invariant part ${p_1}_*\{{\Cal I} \otimes p_2^*{\Cal O}_{{\Bbb
P}^N}(l)\}^{{\Bbb C}^*}$ defines the image of $B_{a_i}$ under the quotient map globally inside
of $U_{X_2} \times ({\Bbb P}^N)_{a_i}//{\Bbb C}^* \hookrightarrow U_{X_2} \times {\Bbb
P}(H^0({\Bbb P}^N,{\Cal O}_{{\Bbb P}^N}(l))^{{\Bbb C}^*}$, verifying the assertion.

It is also worthwhile noting that the quotient $\pi:B_{a_i} \rightarrow B_{a_i}//{\Bbb C}^*$ is
characterized by the following three properties: 

\ \ (i) set-theoretically $B_{a_i}//{\Bbb C}^*$ is the set of equivalence classes of ${\Bbb
C}^*$-orbits in $B_{a_i}$ where two orbits $O(p)$ and $O(q)$ are equivalent if and only if
$\overline{O(p)} \cap \overline{O(q)} \neq \emptyset$ (Note that $B_{a_i}$ is
quasi-elementary.),

\ \ (ii) $U \subset B_{a_i}//{\Bbb C}^*$ is open if and only if $\pi^{-1}(U) \subset B_{a_i}$
is open (with respect to the usual topology), and

\ \ (iii) $\Gamma(U,{\Cal O}_{B_{a_i}//{\Bbb C}^*}) = \Gamma(\pi^{-1}(U),{\Cal
O}_{B_{a_i}})^{{\Bbb C}^*}$. 

For each $p \in B_{a_i} - F_{B_{a_i}}^*$ we choose
$\eta_p:U_p
\rightarrow V_p$ as constructed in Lemma 5-1-2.  Then these $U_p$ form an open covering of
$B_{a_i}$ such that
$U_p =
\pi^{-1}(\pi(U_p))$ by the condition $(\star)$ and that $\pi(U_p) = U_p//{\Bbb C}^*
\overset{\sim}\to{\rightarrow} V_p//{\Bbb C}^*$ is an open subset of an affine (not necessarily
nonsingular) toric variety $X_p//{\Bbb C}^*$.  

\vskip.1in

Step 3. Torification

Step 4. Recovery from Singular to Nonsingular

\vskip.1in

The arguments for these two steps work verbatim if we replace Luna's locally toric charts with the
${\Bbb C}^*$-equivariant anlytic isomorphisms $\eta_p:U_p \rightarrow V_p \subset X_p$ as
described in Lemma 5-1-2. 

\vskip.1in

This completes the proof for the weak factorization theorem for bimeromorphic maps.

\vskip.1in

\S 5-2. Equivariant factorization under group action

\vskip.1in

\proclaim{Theorem 5-2-1 (Equivariant Weak Factorization Theorem)} Let $X_1$ and $X_2$ be
complete nonsingular varieties over an algebraically closed field $K$ of characteristic zero. 
Let $G$ be a group acting on $X_1$ and $X_2$.  (Note that the action does NOT have to be over the
base field $K$.)  Let $\phi:X_1 \dashrightarrow X_2$ be a birational map which is equivariant
under the action of $G$.  Let $X_1 \supset U \subset X_2$ be a common open set over which
$\phi$ is an isomorphism.  Then
$\phi$ can be factored into a sequence of equivariant blowups and blowdowns with smooth
$G$-invariant centers disjoint from
$U$.  (The center may be reducible, i.e., a collection of several disjoint smooth irreducible
components, which is
$G$-invariant as a whole.)  That is to
say, there exists a sequence of $G$-equivarint birational maps between complete nonsingular
varieties with $G$-actions
$$X_1 = V_1 \overset{\psi_1}\to{\dashrightarrow} V_2 \overset{\psi_2}\to{\dashrightarrow}
\cdot\cdot\cdot \overset{\psi_{i-1}}\to{\rightarrow} V_i \overset{\psi_i}\to{\rightarrow}
V_{i+1} \overset{\psi_{i+1}}\to{\rightarrow} \cdot\cdot\cdot
\overset{\psi_{l-2}}\to{\dashrightarrow} V_{l-1}
\overset{\psi_{l-1}}\to{\rightarrow} V_l = X_2$$ where

\ \ (i) $\phi = \psi_{l-1} \circ \psi_{l-2} \circ \cdot\cdot\cdot \circ \psi_2 \circ \psi_1$,

\ \ (ii) $\psi_i$ are isomorphisms over $U$, and

\ \ (iii) either $\psi_i:V_i \dashrightarrow V_{i+1}$ or $\psi_i^{-1}:V_{i+1} \dashrightarrow V_i$
is a morphism obtained by blowing up smooth $G$-invariant center disjoint from $U$.

Moreover, if both $X_1$ and $X_2$ are projective, then we can choose a factorization so that all
the intermediate varieties $V_i$ are projective.
\endproclaim

\demo{Proof}\enddemo Again we only specify the modifications we have to make at each step of the
strategy for the proof described in Chapter 0. Introduction.

\vskip.1in

Step 1. Elimination of points of indeterminacy

\vskip.1in

Though not explicitly stated in Hironaka [2,3] (cf.Bierstone-Milman [1] Encinas-Villamayor [1]),
its method for elimination of points of indeterminacy works equivariantly under a group action. 
Thus Lemma 1-4-2 holds in the equivariant case as well and hence we may assume that $\phi:X_1
\rightarrow X_2$ is a projective birational morphism which is the blowup of $X_2$ along a
$G$-invariant ideal sheaf $J \subset {\Cal O}_{X_2}$ (and hence $\phi$ is $G$-equivariant) with
the support of ${\Cal O}_{X_2}/J$ being disjoint from $U$.

\vskip.1in

Step 2. Construction of a birational cobordism

\vskip.1in

Since the canonical resolution of singularities is $G$-equivariant by the property
$(\spadesuit^{res}-1)$, the construction in Theorem 2-2-2 is also $G$-equivariant.  We remark
that the morphism $\tau:\overline{B} \rightarrow X_2$ is projective and $G$-equivariant as well
as $K^*$-equivariant.  We also remark that a relatively ample line bundle ${\Cal L}$ on
$\overline{B}$ over $X_2$ with a $K^*$-linearization (We set ${\Cal E} = \tau_*({\Cal
L}^{\otimes l})$ in the proof of Theorem 2-2-2.) can be chosen so that it comes with a
$G$-linearization which commutes with the $K^*$-linearization.  In fact, ${\Cal L}$ can be
chosen to be of the form ${\Cal O}_{\overline{B}}(n \cdot p_2^*(\infty) - E)$ where
$p_2^*(\infty)$ is the pull-back of the hyperplane divisor $\infty \in {\Bbb P}^1$ and where $E$
is an exceptional divisor for the birational morphism $\overline{B} \rightarrow W_0 = X_2 \times
{\Bbb P}^1$ and is $G$-invariant.  We identify the local sections with the elements in the
function field
$$\Gamma(U,{\Cal L}) = \{s_f = f \in K(\overline{B});div(f) + n \cdot p_2^*(\infty) - E|_{U} \geq
0\}$$
and choose the linearizations of $K^*$ and $G$ to be the one induced from the actions on the
function field $K(W_0 = X_2 \times {\Bbb P}^1) = K(\overline{B})$, i.e.,
$$\align
t^*s_f &:= t^*f \in \Gamma(t^{-1}(U),{\Cal L}) \text{\ for\ }t \in K^*\\
g^*s_f &:= g^*f \in \Gamma(g^{-1}(U),{\Cal L}) \text{\ for\ }g \in G.\\
\endalign$$
As the actions of $K^*$ and $G$ on $X_2 \times {\Bbb P}^1$ obviously commute with each other, so
do the linearizations defined as above.  It also follows then that the linearization of $K^*$
on ${\Cal L}$ commutes with the twisted linearizations discussed in \S 2-3.  This implies that
the quasi-elememtary pieces $B_{a_i}$ are $G$-invariant.  The rest of the argument goes without
any change.   

\vskip.1in

Step 3. Torification

\vskip.1in

Since the definition of the torific ideals is canonical in terms of the action of $K^*$, which
commutes with that of $G$, the torific ideals are $G$-invariant.  Thus the torification is
$G$-equivariant.  The argument in this part goes without any change.

\vskip.1in

Step 4. Recovery from Singular to Nonsingular

\vskip.1in

First we look at Lemma 4-1-2.  Note that $f_{i-}:W_i^{tor} \rightarrow W_i$ is the blowup of the
$G$-invariant ideal sheaf $I_{i-} = (f_{i-})_*{\Cal O}_{W_i^{tor}}(- E_i^{tor})$, where
$E_i^{tor}$ is a $G$-invariant effective Cartier divisor, since the torific ideal is
$G$-invariant by Step 3.  (The same claim holds replacing the negative sign $-$ with the
positive sign for the subscripts.)  Secondly we look at Proposition 4-1-3.  Since the canonical
resolution of singularities and the canonical principalization of ($G$-invariant) ideals are
$G$-equivariant, we see that
$h_{i-}:W_{i-}^{can} \rightarrow W_i^{tor}$ and $h_{i+}:W_{i+}^{can} \rightarrow W_i^{tor}$ are
also $G$-equivariant.  Therefore, we also conclude that
$\varphi_i^{can}:(U_{W_{i-}^{can}},W_{i-}^{can}) \dashrightarrow
(U_{W_{i+}^{can}},W_{i+}^{can})$ is $G$-equivariant V-toroidal birational map. 

Thus the proof of the equivariant weak factorization theorem is reduced to the following:

\proclaim{Theorem 5-2-1 (Equivariant Strong Factorization Theorem for Toroidal Birational Maps)}
Let $(U_{W_1},W_1)$ and $(U_{W_2},W_2)$ be nonsingular toroidal embeddings.  Let $G$ be
a group acting on them as automorphisms of toroidal embeddings.  Let
$$\varphi:(U_{W_1},W_1) \dashrightarrow (U_{W_2},W_2)$$
be a proper and toroidal birational map which is $G$-equivariant.  Then $\phi$ can be factored
into a sequence of toroidal blowups immediately followed by toroidal blowdowns with
$G$-invariant smooth centers (which may be reducible).
\endproclaim

\demo{Proof}\enddemo First note that the process of Lemma 1-2-6, making a general nonsingular
toroidal embedding into the one without self-intersection, is $G$-equivariant.  Thus we may
assume both $(U_{W_1},W_1)$ and $(U_{W_2},W_2)$ are without self-intersection. 

Since $\varphi$ is toroidal, by Definition 1-2-3 there exists a toroidal embedding $(U_Z,Z)$
which dominates $(U_{W_1},W_1)$ and $(U_{W_2},W_2)$ by proper and toroidal birational morphisms
$$(U_{W_1},W_1) \overset{f_1}\to{\leftarrow} (U_Z,Z) \overset{f_2}\to{\rightarrow}
(U_{W_2},W_2).$$
First by the process of Lemma 1-2-6, we may assume that $(U_Z,Z)$ is a nonsingular toroidal
embedding without self-intersection.  By Kempf-Knudsen-Mumford-SaintDonat [1] Abramovich-Karu
[1], the morphisms $f_1$ and $f_2$ correspond to the refinements of the conical complexes
$$\align
\Delta_{f_1}: \Delta_Z &\rightarrow \Delta_{W_1} \\
\Delta_{f_2}: \Delta_Z &\rightarrow \Delta_{W_2}. \\
\endalign$$
Remark that we have a natural homomorphism
$$h:G \rightarrow Aut(\Delta_{W_1}),$$
sending $g \in G$, considered as an automorphism $g:(U_{W_1},W_1) \rightarrow (U_{W_1},W_1)$, to
an automorphism of the conical complex $h(g) = \Delta_g:\Delta_{W_1} \rightarrow
\Delta_{W_1}$.  Since $Aut(\Delta_{W_1})$ is a finite group, so is $\overline{G} :=
G/Ker(h)$.  Let $\Delta_{\widetilde{f_1}}: \Delta_{\tilde Z} \rightarrow \Delta_{W_1}$ be the
smallest common refinement of $\{\Delta_{\overline{g}} \circ \Delta_{f_1}:\Delta_Z
\rightarrow \Delta_{W_1};\overline{g} \in \overline{G}\}$.  Though $G$ may not act on
the original $(U_Z,Z)$, by replacing $(U_Z,Z)$ with $(U_{\tilde Z},{\tilde Z})$ associated to the
refinement
$\Delta_{\widetilde{f_1}}:\Delta_{\tilde Z} \rightarrow \Delta_{W_1}$ (and then taking the
canonical resolution of singularities, if necessary), we may assume from the beginning that $G$
acts on the toroidal embedding $(U_Z,Z)$ and $f_1:(U_Z,Z) \rightarrow (U_{W_1},W_1)$ is
$G$-equivariant and hence that $f_2:(U_Z,Z) \rightarrow (U_{W_2},W_2)$ is $G$-equivariant as
well.

Thus we only have to provide factorization for a $G$-equivariant toroidal birational morphism
$f:(U_Z,Z) \rightarrow (U_W,W)$ between nonsingular toroidal embeddings without
self-intersection, associated to the refinement $\Delta_f:\Delta_Z \rightarrow \Delta_W$ (cf.
the argument in the proof of Theorem 4-2-1).  Let
$\overline{G} \subset Aut(\Delta_W)$ be the subgroup of the automorphism group of the conical
complex $\Delta_W$ induced by the action of $G$.  Note that $\overline{G}$ acts on $\Delta_Z$
as well as on $\Delta_W$ and that $\Delta_f$ is $\overline{G}$-equivariant. 

Now we use an idea of Abramovich-Wang [1].  Take the barycentric star subdivision
$\Delta_b:\Delta_{\hat W} \rightarrow \Delta_W$ (cf. Definition 2.1 in Abramovich-Matsuki-Rashid
[1]), which is $\overline{G}$-equivariant.  It is straightforward to see that the corresponding
morphism $b:(U_{\hat W},{\hat W})
\rightarrow (U_W,W)$ between nonsingular toroidal embeddings without self-intersection is the
blowup of a $G$-invariant toroidal ideal $I$ on $(U_W,W)$ (cf. the argument in the proof of
Theorem 2-2-2).  Let $\lambda:{\hat Z} \rightarrow Z$ be the canonical
principalization of $f^{-1}I \cdot {\Cal O}_Z$, which is a $G$-invariant and toroidal ideal
since $f$ is $G$-equivariant and toroidal.  Thus $\lambda:(U_{\hat Z} = \lambda^{-1}(U_Z),{\hat
Z}) \rightarrow (U_Z,Z)$ is a $G$-equivariant and toroidal morphism between nonsingular
toroidal embeddings without self-intersection by the property
$(\spadesuit^{can}-1)$ of the canonical principalization, and hence so is the induced morphism
${\hat f}:(U_{\hat Z},{\hat Z}) \rightarrow (U_{\hat W},{\hat W})$.  Let $\Delta_{\hat
f}:\Delta_{\hat Z} \rightarrow
\Delta_{\hat W}$ be the associated map of the conical complexes, which is
$\overline{G}$-equivariant.  Observe that, thanks to the process of taking the barycentric star
subdivision, the conical complex $\Delta_{\hat W}$ has the following property $(\natural)$:

$(\natural)$ If a cone $\sigma \in \Delta_{\hat W}$ is mapped to itself by some $\overline{g}
\in \overline{G}$, then $\overline{g}$ acts as an identity on the cone $\sigma$.

Observe also that the conical complex $\Delta_{\hat Z}$ satisfies the same property
$(\natural)$. 

It follows that the quotient $\Delta_{\hat W}/\overline{G}$ (resp. $\Delta_{\hat
Z}/\overline{G}$) is again a nonsingular conical complex and the projection $\Delta_{\hat W}
\rightarrow \Delta_{\hat W}/\overline{G}$ (resp. $\Delta_{\hat Z} \rightarrow \Delta_{\hat
Z}/\overline{G}$) maps the cones in $\Delta_{\hat W}$ (resp. $\Delta_{\hat Z}$) isomorphically to
the cones in $\Delta_{\hat W}/\overline{G}$ (resp. $\Delta_{\hat Z}/\overline{G}$) and we
have a map of nonsingular conical complexes
$\Delta_{\hat f}/\overline{G}:\Delta_{\hat Z}/\overline{G}
\rightarrow
\Delta_{\hat W}/\overline{G}$.

Now apply the combinatorial algorithm of Morelli [1][2] Abramovich-Matsuki-Rashid [1] to have a
conical complex $\Delta_V/\overline{G}$ which is a comon refinement of $\Delta_{\hat
W}/\overline{G}$ and $\Delta_{\hat Z}/\overline{G}$
$$\Delta_{\hat W}/\overline{G} \hskip.1in \overset{\Delta_{\mu_{\hat
W}}/\overline{G}}\to{\leftarrow} \hskip.1in 
\Delta_V/\overline{G} \hskip.1in \overset{\Delta_{\mu_{\hat Z}}/\overline{G}}\to{\rightarrow}
\hskip.1in
\Delta_{\hat Z}/\overline{G}$$
where $\mu_{\hat W}/\overline{G}$ (resp. $\mu_{\hat Z}/\overline{G}$) is a sequence of smooth star
subdivisions.  We can pull back the refinement, i.e., set
$$\Delta_V = \Delta_{\hat W} \times_{\Delta_{\hat W}/\overline{G}}\Delta_V/\overline{G} =
\Delta_{\hat Z} \times_{\Delta_{\hat Z}/\overline{G}}\Delta_V/\overline{G}$$ 
to obtain the
$\overline{G}$-equivariant maps
$$\Delta_{\hat W} \overset{\Delta_{\mu_{\hat W}}}\to{\leftarrow}
\Delta_V \overset{\Delta_{\mu_{\hat Z}}}\to{\rightarrow} 
\Delta_{\hat Z}$$
where $\mu_{\hat W}$ (resp. $\mu_{\hat Z}$) is a sequence of smooth star
subdivisions, possibly taking several disjoint smooth star subdivisions simultaneously. 
Geometrically interpreted, we have the assertion of the theorem.

\vskip.1in

This completes the proof of Theorem 5-2-1 and hence the proof of the equivariant weak
factorization theorem.

\vskip.1in

\proclaim{Remark 5-2-2}\endproclaim

Equivariant Weak Factorization Theorem for Bimeromorphic Maps can be proved in almost an
identical manner, combining the modifications in \S 5-1 and \S 5-2.  (Note that the relatively
(very) ample line bundle ${\Cal L} = {\Cal L}_l$ constructed inductively in \S 5-1 has not only
a ${\Bbb C}^*$-lnearization but also a natural $G$-linearization which commutes with the ${\Bbb
C}^*$-linearization.)  Details are left to the reader as an exercise.

\vskip.1in

\S 5-3. Factorization over a non-algebraically closed field

\vskip.1in

\proclaim{Theorem 5-3-1 (Weak Factorization Theorem over a NON-algebraically closed field)} Let
$\phi:X_1 \dashrightarrow X_2$ be a birational map between varieties smooth and proper over a
field $K$, which may not be algebraically closed, of characteristic zero.  Let $X_1 \supset U
\subset X_2$ be a common open subset over which $\phi$ is an isomorphism.  Then $\phi$ can be factored into a
sequence of blowups and blowdowns with irreducible centers smooth over $K$ and disjoint from $U$. 
That is to say, there exists a sequence of bbirational maps between varieties smooth and proper
over $K$
$$X_1 = V_1 \overset{\psi_1}\to{\dashrightarrow} V_2 \overset{\psi_2}\to{\dashrightarrow}
\cdot\cdot\cdot \overset{\psi_{i-1}}\to{\rightarrow} V_i \overset{\psi_i}\to{\rightarrow}
V_{i+1} \overset{\psi_{i+1}}\to{\rightarrow} \cdot\cdot\cdot
\overset{\psi_{l-2}}\to{\dashrightarrow} V_{l-1}
\overset{\psi_{l-1}}\to{\rightarrow} V_l = X_2$$ where

\ \ (i) $\phi = \psi_{l-1} \circ \psi_{l-2} \circ \cdot\cdot\cdot \circ \psi_2 \circ \psi_1$,

\ \ (ii) $\psi_i$ are isomorphisms over $U$, and

\ \ (iii) either $\psi_i:V_i \dashrightarrow V_{i+1}$ or $\psi_i^{-1}:V_{i+1} \dashrightarrow V_i$
is a morphism obtained by blowing up an irreducible center smooth over $K$ and disjoint from $U$.

Moreover, if both $X_1$ and $X_2$ are projective, then we can choose a factorization so that all
the intermediate complex manifolds $V_i$ are projective.
\endproclaim 

\demo{Proof}\enddemo We specify the modifications we have to make at each step of the proof of
the strategy for the proof described in Chapter 0. Introduction.

\vskip.1in

Step 1. Elimination of points of indeterminacy

\vskip.1in

Hironaka's method for elimination of points of indeterminacy works over ANY field of
characteristic zero (cf. Remark 1-4-2).  Hence again via Lemma 1-4-2, we may assume that $\phi:X_1
\rightarrow X_2$ is a projective birational morphism which is the blowup of of $X_2$ along an
ideal sheaf $J
\subset {\Cal O}_{X_2}$ defined over $K$ with the support of ${\Cal O}_{X_2}/J$ being disjoint
from $U$.

\vskip.1in

Step 2. Construction of a birational cobordism

\vskip.1in

Since the canonical resolution of singularities also works over any field of characteristic
zero, the construction of a birational cobordism goes without any change.

\vskip.1in

Step 3. Torification

Step 4. Recovery from Singular to Nonsingular

\vskip.1in

Here we carry the argument by taking the base change from the original field $K$ to its algebraic
closure $\overline{K}$.

As $X_1$ and $X_2$ share a common open subset $U$, after base change

$$\align
\overline{X_1} &:= X_1 \times_{\roman{Spec}\ K}\roman{Spec}\ \overline{K} \\
(\text{resp.\ }\overline{X_2} &:= X_2 \times_{\roman{Spec}\ K}\roman{Spec}\ \overline{K},\\
\overline{\{B\}} &:= B \times_{\roman{Spec}\ K}\roman{Spec}\ \overline{K})\\
\overline{U} &:= U \times_{\roman{Spec}\ K}\roman{Spec}\ \overline{K})\\
\endalign$$
it splits into a finite number of disjoint components smooth over $\overline{K}$
$$(\overline{X_1})_{g_j} \hskip.1in (\text{resp}.\ (\overline{X_2})_{g_j}, \hskip.1in
(\overline{\{B\}})_{g_j}, (\overline{U})_{g_j})$$
such that they are conjugate to each other under the action of the Galois group $G =
Gal(\overline{K}/K)$ where we set $G_e$ to be the decomposition group of some component
$(\overline{X_1})_e$, i.e.,
$$G_e := \{g \in G;g((\overline{X_1})_e) = (\overline{X_1})_e\}$$
and where $C = \{g_o = e, \cdot\cdot\cdot, g_j, \cdot\cdot\cdot\}$ is a complete representative
modulo $G_e$, and that $(\overline{\{B\}})_{g_j}$ is a birational cobordism for the birational
map $(\overline{\varphi})_{g_j}:(\overline{X_1})_{g_j} \dashrightarrow (\overline{X_2})_{g_j}$
respecting the open subset
$(\overline{U})_{g_j}$.  (Note that we use the notation $\overline{\{B\}}$ in order to
distingush it from the compactified birational cobordism $\overline{B}$.)

Now we apply the argument for equivariant weak factorization in \S 5-2 to the birational
map $(\overline{\varphi})_e:(\overline{X_1})_e \dashrightarrow (\overline{X_2})_e$ over
an algebraically closed field $\overline{K}$, induced from the birational cobordism
$(\overline{\{B\}})_e$, to obtain a $G_e$-equivariant factorization

$$\align
(\overline{X_1})_e &= (\overline{V_1})_e \overset{(\overline{\psi_1})_e}\to{\dashrightarrow}
(\overline{V_2})_e \overset{(\overline{\psi_2})_e}\to{\dashrightarrow} \\
&\cdot\cdot\cdot
\overset{(\overline{\psi_{i-1}})_e}\to{\rightarrow} (\overline{V_i})_e
\overset{(\overline{\psi_i})_e}\to{\rightarrow} (\overline{V_{i+1}})_e
\overset{(\overline{\psi_{i+1}})_e}\to{\rightarrow} \cdot\cdot\cdot \\
&\overset{(\overline{\psi_{l-2}})_e}\to{\dashrightarrow} (\overline{V_{l-1}})_e
\overset{(\overline{\psi_{l-1}})_e}\to{\rightarrow} (\overline{V_l})_e = (\overline{X_2})_e\\
\endalign$$
where

\ \ (i) $(\overline{\phi})_e = (\overline{\psi_{l-1}})_e \circ (\overline{\psi_{l-2}})_e \circ
\cdot\cdot\cdot
\circ
(\overline{\psi_2})_e
\circ
(\overline{\psi_1})_e$,

\ \ (ii) $(\overline{\psi_i})_e$ are isomorphisms over $(\overline{U})_e$, and

\ \ (iii) either $(\overline{\psi_i})_e:(\overline{V_i})_e \dashrightarrow
(\overline{V_{i+1}})_e$ or
$(\overline{\psi_i})_e^{-1}:(\overline{V_{i+1}})_e
\dashrightarrow (\overline{V_i})_e$ is a morphism obtained by blowing up smooth $G_e$-invariant
center $(\overline{C_i})_e$ disjoint from $(\overline{U})_e$.

\vskip.1in

Taking the conjugate induced by the action of the element $g_j \in G$ in the Galois group, we
obtain the corresponding factorization for $(\overline{\psi_i})_{g_j}:(\overline{X_1})_{g_j}
\dashrightarrow (\overline{X_2})_{g_j}$

$$\align
(\overline{X_1})_{g_j} &= (\overline{V_1})_{g_j}
\overset{(\overline{\psi_1})_{g_j}}\to{\dashrightarrow} (\overline{V_2})_{g_j}
\overset{(\overline{\psi_2})_{g_j}}\to{\dashrightarrow} \\ &\cdot\cdot\cdot
\overset{(\overline{\psi_{i-1}})_{g_j}}\to{\rightarrow} (\overline{V_i})_{g_j}
\overset{(\overline{\psi_i})_{g_j}}\to{\rightarrow} (\overline{V_{i+1}})_{g_j}
\overset{(\overline{\psi_{i+1}})_{g_j}}\to{\rightarrow} \cdot\cdot\cdot \\
&\overset{(\overline{\psi_{l-2}})_{g_j}}\to{\dashrightarrow} (\overline{V_{l-1}})_{g_j}
\overset{(\overline{\psi_{l-1}})_{g_j}}\to{\rightarrow} (\overline{V_l})_{g_j} =
(\overline{X_2})_{g_j}.\\
\endalign$$

Collectively, we obtain the factorization over $\overline{K}$ for
$$\overline{\varphi}:\overline{X_1} \dashrightarrow \overline{X_2}$$
into blowups with centers $\coprod_{g_j \in C} (\overline{C_i})_{g_j}$, which are $G$-invariant
and smooth over $\overline{K}$.

Therefore, we obtain the factorization over $K$ for $\phi:X_1 \dashrightarrow X_2$

$$X_1 = V_1 \overset{\psi_1}\to{\dashrightarrow} V_2 \overset{\psi_2}\to{\dashrightarrow}
\cdot\cdot\cdot \overset{\psi_{i-1}}\to{\rightarrow} V_i \overset{\psi_i}\to{\rightarrow}
V_{i+1} \overset{\psi_{i+1}}\to{\rightarrow} \cdot\cdot\cdot
\overset{\psi_{l-2}}\to{\dashrightarrow} V_{l-1}
\overset{\psi_{l-1}}\to{\rightarrow} V_l = X_2$$ where

\ \ (i) $\phi = \psi_{l-1} \circ \psi_{l-2} \circ \cdot\cdot\cdot \circ \psi_2 \circ \psi_1$,

\ \ (ii) $\psi_i$ are isomorphisms over $U$, and

\ \ (iii') either $\psi_i:V_i \dashrightarrow V_{i+1}$ or $\psi_i^{-1}:V_{i+1} \dashrightarrow
V_i$ is a morphism obtained by blowing up a center $C_i$ smooth over $K$ and disjoint from $U$.

In order to guarantee that the centers of blowups to be irreducible and satisfy the condition
(iii), instead of blowing up $C_i$ which may be reducible, we blow up each irreducible component
one by one and replace the original sequence with the corresponding refinement.

The ``moreover" part concerning the projectivity is also obvious from the construction.

This completes the proof of the weak factorization theorem over a (possibly) NON-algebraically
closed field.

\vskip.1in

\S 5-4. Factorization in the logarithmic category

\vskip.1in

\proclaim{Theorem 5-4-1 (Weak Factorization Theorem in the Logarithmic Category)} Let
$(U_{X_1},X_1)$ and
$(U_{X_2},X_2)$ be complete nonsingular toroidal embeddings over an algebraically closed field $K$
of characteristic zero.  Let $\phi:(U_{X_1},X_1) \dashrightarrow (U_{X_2},X_2)$ be a birational
map which is an isomorphism over $U_{X_1} = U_{X_2}$.  Then
$\phi$ can be factored into a sequence of blowups and blowdowns with smooth ADMISSIBLE and
irreducible centers disjoint from
$U_{X_1} = U_{X_2}$.  That is to
say, there exists a sequence of birational maps between complete nonsingular
toroidal embeddings
$$\align
(U_{X_1},X_1) &= (U_{V_1},V_1) \overset{\psi_1}\to{\dashrightarrow} (U_{V_2},V_2)
\overset{\psi_2}\to{\dashrightarrow} \\
&\cdot\cdot\cdot \overset{\psi_{i-1}}\to{\rightarrow} (U_{V_i},V_i)
\overset{\psi_i}\to{\rightarrow} (U_{V_{i+1}},V_{i+1}) \overset{\psi_{i+1}}\to{\rightarrow}
\cdot\cdot\cdot \\
&\overset{\psi_{l-2}}\to{\dashrightarrow} (U_{V_{l-1}},V_{l-1})
\overset{\psi_{l-1}}\to{\rightarrow} (U_{V_l},V_l) = (U_{X_2},X_2) \\
\endalign$$ 
where

\ \ (i) $\phi = \psi_{l-1} \circ \psi_{l-2} \circ \cdot\cdot\cdot \circ \psi_2 \circ \psi_1$,

\ \ (ii) $\psi_i$ are isomorphisms over $U_{V_i}$, and

\ \ (iii) either $\psi_i:(V_i,U_{V_i}) \dashrightarrow (U_{V_{i+1}},V_{i+1})$ or
$\psi_i^{-1}:(U_{V_{i+1}},V_{i+1}) \dashrightarrow (U_{V_i},V_i)$ is a morphism obtained by
blowing up a smooth irreducible center disjoint from $U_{V_i} = U_{V_{i+1}}$ and transversal to
the boundary
$D_{V_i} = V_i - U_{V_i}$ or $D_{V_{i+1}} = V_{i+1} - U_{V_{i+1}}$, i.e., at each point $p \in
V_i$ or
$p
\in V_{i+1}$ there exists a regular coordinate system $\{x_1, \cdot\cdot\cdot, x_n\}$ in a
neighborhood
$p \in U_p$ such that $D_{V_i} \cap U_p (\text{or\ } D_{V_{i+1}} \cap U_p) = \{\prod_{j \in J} x_j
= 0\}$ and
$C_i
\cap U_p =
\{\prod_{j
\in J \cup J'}x_j = 0\}$ for some subsets $J, J' \subset \{1, \cdot\cdot\cdot, n\}$.

Moreover, if both $(U_{X_1},X_1)$ and $(U_{X_2},X_2)$ are projective, then we can choose a
factorization so that all the intermediate toroidal embeddings $(U_{V_i},V_i)$ are projective.
\endproclaim

\demo{Proof}\enddemo We specify the modifications we have to make at each step of the proof of
the strategy for the proof described in Chapter 0. Introduction.

\vskip.1in

Step 1. Elimination of points of indeterminacy

\vskip.1in

Lemma 1-4-1 holds in the above-mentioned logarithmic category where all the centers of the
blowups necessary for elimination of points of indeterminacy can be taken to be admissible as
well as smooth.  Thus we may assume that $\phi:(U_{X_1},X_1) \rightarrow (U_{X_2},X_2)$ is a
projective morphism which is the blowup of $X_2$ along an ideal sheaf $J \subset {\Cal O}_{X_2}$
with the support of ${\Cal O}_{X_2}/J$ being disjoint from $U_{X_2} = U_{X_1}$.

\vskip.1in

Step 2. Construction of a birational cobordism

\vskip.1in

In an identical manner to the proof of Theorem 2-2-2, we can construct a compactified birational
cobordism $\tau:(U_{\overline B},\overline{B}) \rightarrow (U_{X_2},X_2)$ for $\phi$ projective
over $X_2$ such that $\tau^{-1}(D_{X_2}) = D_{\overline{B}}$ is a divisor
with only normal crossings where $D_{X_2} = X_2 - U_{X_2}$ and $D_{\overline{B}} = \overline{B}
- U_{\overline{B}}$ and that
$$\align
(U_B,B)_+/K^* &= (U_{X_2},X_2) \\
(U_B,B)_-/K^* &= (U_{X_1},X_1) \\
\endalign$$
as toroidal embeddings.

When we construct a locally toric chart for a point $p \in B$ (Though $(U_B,B)$ is a nonsingular
toroidal embedding, the action of $K^*$ may not be toroidal.  Thus we still call the chart a
locally toric chart), we use the method (A) in Proposition 1-3-4 so that $U_p \cap D_B = V_p
\cap D_B$ is included in $\eta_p^{-1}(D_{X_p})$, where $D_B = B - U_b$ and $D_{X_p} = X_p -
U_{X_p}$.  This can be done, in the notation of the proof (A) of Proposition 1-3-4, by choosing a
regular system of parameters
$$f_1, \cdot\cdot\cdot, f_n \in m_q$$
around the limit point $q \in \overline{O(p)}$, consisting of eigenfunctions, so that each
irreducible component $D$ of $V_q \cap D_{\overline{B}}$ is defined by some
coordinate among $f_1, \cdot\cdot\cdot, f_n$ (shrinking $V_q$ if necessary).  In fact, since the
whole boundary $D_{\overline{B}}$ is invariant under the action of $K^*$, so is $D$
and hence its defining ideal $I_D \subset A(V_q)$ is also $K^*$-invariant.  Since $K^*$ is
reductive, the ideal $I_D$ splits into the eigenspaces according to the characters under the
action of $K^*$.  As
$D$ is nonsingular at
$q$, there exists an eigenfunction $f_D \in I_D$ and a coordinate , say $f_1$, such that
$$f_D\ \text{mod}\ m_q = f_1\ \text{mod}\ m_q.$$
By replacing $f_1$ with $f_D$ in the regular system of parameters, we achieve the goal.

\vskip.1in

Step 3. Torification

Step 4. Recovery from Singular to Nonsingular

\vskip.1in

In the torification of the quasi-elementary cobordism $B = B_{a_i}$, after blowing up the
torific ideal to obtain $\mu:B^{tor} \rightarrow B$, we only have to show that
$$(B^{tor} - \{D^{tor} \cup \mu^{-1}(D_B)\},B^{tor})$$
has a toroidal structure with respect to which the induced action of $K^*$ is toroidal.  But
this follows immediately, since already
$$(B^{tor} - D^{tor},B^{tor})$$
has a toroidal structure with respect to which the induced action of $K^*$ is toroidal and since
the additional boundary divisors coming from $D_B$ coorespond to the toric coordinate divisors in
the locally toric charts chosen as above.  (In the language of Torification Principle for Toric
Varieties, we do NOT have to ``remove" the divisors coming from $D_B$ from the boundary.) 
The argument for recovery from Singular to Nonsingular also goes without any change, considering
the toroidal structures with the divisors coming from $D_B$ added to the boundaries.

\vskip.1in

This completes the proof of the weak factorization theorem in the logarithmic category.

\vskip.1in

\proclaim{Corollary 5-4-2} The Hodge numbers are birational invariants for nonsingular minimal
models, i.e., if $X_1$ and $X_2$ are projective nonsingular varieties over
${\Bbb C}$ whose canonical divisors are nef and if they are birational to each other, then we
have
$$h^{p,q}(X_1) = h^{p,q}(X_2) \hskip.1in \forall p,q.$$
\endproclaim

\demo{Proof}\enddemo This application of the weak factorization theorem to a theorem of
Batyrev was communicated to us by J. Denef, F. Loeser and W. Veys.  It is easy to see that $X_1$
and
$X_2$, being birational and minimal, are isomorphic outside of some codimension 2 loci.  By some
sequences of blowups with admissible (in the usual sense) centers over the specified loci, we can
reach the situation of Theorem 5-4-1, where again the birational map is factored as a sequence of
blowups and blowdowns with smooth and admissible (in the sense of Theorem 5-4-1) centers over the
same specified loci.  Then it is straightforward to see the required invariance of the Hodge
numbers via the invariance of Batyrev's stringy
$E$-function through such blowups and blowdowns with admissible centers.  We refer the details
of the argument to the papers in the related field (cf. Batyrev [1] Denef-Loeser [1] Koll\'ar [2]
Wang [1]).

There is also a slightly different approach to the above theorem by considering the Grothendieck
group of the Hodge structures and looking for some elementary but explicit invariant form within
the group under blowups and blowdowns with admissible centers (cf. Arapura-Matsuki [1]).

\vskip.1in

\proclaim{Remark 5-4-3}\endproclaim

(i) We can also combine some of the generalizations above, as in Remark 5-2-2, and prove, e.g., 

\ \ \ $\circ$ an equivariant version of
the weak factorization theorem in the logarithmic category or in the usual category over a
(possibly) NON-algebraically closed field of characteric zero

\ \ \ $\circ$ the weak factorization theorem
of bimeromorphic maps in the logarithmic category with or without a group action in
consideration.

The modifications we have to make are simple combinations of those discussed in \S 5-1 through
\S 5-4 and we leave the details to the reader.

\vskip.1in

(ii) One can try to present the generalizations in \S 5-1 through \S 5-4 in a more unified
manner, including the case of algebraic spaces (cf. Section 5 in
Abramovich-Karu-Matsuki-W{\l}odarczyk [1]).  We chose, however, our rather repetitive manner in
order to emphasize minor but subtle differences.

\newpage

$$\bold{CHAPTER\ 6.\ PROBLEMS}$$

\vskip.2in

In this chapter, we discuss a couple of problems related to our proof of the weak factorization
theorem for birational maps.

\vskip.1in

\S 6-1. Effectiveness of the construction

\vskip.1in

As our proof factorizes a given birational map explicitly via the birational cobordism, it is
not merely existential but also constructive.  On the other hand, if one asks the following
effectiveness question, our construction falls short of giving an affirmative answer at several
places.

\proclaim{Question 6-1-1}  Suppose that nonsingular projective varieties $X_1 \subset {\Bbb
P}^{n_1}, X_2
\subset {\Bbb P}^{n_2}$ as well as a birational map between them $\phi:X_1 \dashrightarrow
X_2$ are given by a specific set of equations in terms of homogeneous coordinates of the
projective spaces.  Can we give a factorization of $\phi$ into blowups and blowdowns with smooth
centers in an effective way from the given set of equations ?
\endproclaim 

We list the places where our construction falls short of being effective (or rather to say, the
places where the AUTHOR does not know how to carry out the process in an effective way).

\vskip.1in 

\noindent Step 1: Reduction to the case where $\phi$ is a projective birational morphism

Hironaka's method of elimination of the points of indeterminacy in Lemma 1-4-1

\ \ \ $\circ$ the construction of the ideal $I$ in Remark 1-4-2 on $X_1$ whose blowup
factors through $X_2$ (after taking the graph of
$\phi$ and assuming
$\phi^{-1}$ is a morphism)  

\ \ \ $\circ$ canonical principalization of the ideal $I$ (Canonical resolution of
singularities as well as canonical principalization of ideals are constructive as presented in
Bierstone-Milman [1] Villamayor [1] Encinas-Villamayor [1], while the author does not know if
their algorithms can be made effective.)

\vskip.1in

\noindent Step 2: Factorization into locally toric birational maps

Construction of a birational cobordism in Theorem 2-2-2

\ \ \ $\circ$ the construction of the ideal $J \subset {\Cal O}_{X_2}$ whose blowup gives the
projective morphism $\phi:X_1 \rightarrow X_2$, more specifically finding the generators of $J$
in terms of the homogeneous coordinates for $X_2$

\ \ \ $\circ$ canonical resolution of singularities

\ \ \ $\circ$ the construction of a relatively ample line bundle ${\Cal L}$ on $\overline{B}$,
equipped with a $K^*$-action so that for sufficiently large $l \in {\Bbb N}$ the direct image
sheaf ${\Cal E} =
\tau_*({\Cal L}^{\otimes l})$ has the decomposition into the eigen spaces ${\Cal E} = \oplus_{b
\in {\Bbb Z}}{\Cal E}_b$ and that the compactified birational cobordism $\overline{B}$ has an
equivariant embedding $\overline{B} \hookrightarrow {\Bbb P}({\Cal E})$

\vskip.1in

Actually the effectiveness for the above could be an immediate consequence of the
effectiveness for canonical resolution.  In fact, in the process of canonical resolution using
any one of the algorithms given by Bierstone-Milman [1] Villamayor [1] Encinas-Villamayor [1]
$$\overline{B} = W_m \rightarrow \cdot\cdot\cdot \rightarrow W_{i+1}
\overset{\mu_i}\to{\rightarrow} W_i \rightarrow \cdot\cdot\cdot \rightarrow W_1 = X_2 \times
{\Bbb P}^1$$ 
if we can find inductively an equivariant embedding
$$W_i \subset X_2 \times {\Bbb P}^{k_i} \subset {\Bbb P}^{n_2} \times {\Bbb P}^{k_i}$$
where $K^*$ acts only on the second factor ${\Bbb P}^{k_i}$ with the semi-invariant homogeneous
coordinates, which induces a linearization on ${\Cal L}_i = p_2^*{\Cal O}_{{\Bbb P}^{k_i}}(1)$,
and if we can find effectively some generators consisting of bihomogeneous polynomials (in
terms of homogeneous coordinates for ${\Bbb P}^{n_2}$ and ${\Bbb P}^{k_i}$) of some fixed
bidegree for the ideal of the center of blowup for
$\mu_i$, which is $K^*$-invariant, then it is easy to realize the inductive situation for
$W_{i+1}$.  Then we can set ${\Cal L} = {\Cal L}_m$ and $l = 1$.

\vskip.1in

\ \ \ $\circ$ how to find Luna's locally toric charts as described in Lemma 2-4-4 satisfying
the condition $(\star)$

\vskip.1in

\noindent Step 3: Factorization into toroidal birational maps

Torification

Once Luna's toric charts satisfying the condition $(\star)$ are chosen, all the analysis is
reduced to the toric case, where the process of torification is effctive as well as
constructive.

\vskip.1in

\noindent Step 4: Recovery of nonsingularity

\vskip.1in

Again via Luna's locally toric (toroidal) charts all the analysis here is reduced to the toric
case, where the process, including the canonical resolution of singularities and canonical
principalization of (toric) ideals, is effective as well as constructive.

\vskip.1in

\noindent Step 5. Factorization of toroidal birational maps among nonsingular toroidal embeddings

Morelli's combinatorial algorithm

There are two places where the algorithm refers to some existential theorems in Morelli [1][2]
Abramovich-Matsuki-Rashid [1].

\ \ \ $\circ$ Toric Moishezon Theorem

Though we refer to Hironaka's elimination of points of indeterminacy in the two papers above,
the algorithm in Deconcini-Procesi [1] for toric birational maps is effective as well as
constructive.

\ \ \ $\circ$ Sumihiro's equivariant completion theorem

In the toric case, ``filling up" cones between the lower face and upper face of the
``vacuum" for the fan representing the birational cobordism can be made effective 
as well as effective as in Matsuki-Rashid [1]. 

\vskip.2in

\S 6-2. Toroidalization Conjecture and Strong Factorization Conjecture

\vskip.1in

\proclaim{Conjecture 6-2-1 (Toroidalization Algorithm Conjecture)} Let $f:(U_X,X)
\rightarrow (U_Y,Y)$ be a proper (not necessarily birational) morphism between nonsingular
toroidal embeddings such that $f^{-1}(Y - U_Y) = X - U_X$ and that
$f|_{U_X}:U_X \rightarrow U_Y$ is smooth.  Then there should exist an algorithm
which specifies the way to take sequences of blowups $\sigma_X:(U_{X'},X') \rightarrow
(U_X,X)$ and $\sigma_Y:(U_{Y'},Y') \rightarrow (U_Y,Y)$ with smooth (and possibly
admissible (!?)) centers in the boundaries so that the resulting morphism between nonsingular
toroidal embeddings
$\phi':(U_{X'},X') \rightarrow (U_{Y'},Y')$ is toroidal.
\endproclaim

It is easy to see that via the method of elimination of points of indeterminacy and (canonical)
resolution of singularities the strong factorization conjecture follows immediately from the
toroidalization conjecture.

As we briefly said in the introduction, a morphism between toroidal embeddings being TOROIDAL
(cf. Abramovich-Karu [1] Karu [1]) is equivalent to its being log-smooth (cf. Kato [1] Iitaka
[1]).  Thus the above conjecture asserts that in the logarithmic category of toroidal embeddings
we should have an algorithm which resolves the singularities of a morphism $\phi:X \rightarrow Y$,
much like the usual ones which resolves the singularities of varieties.  The algorithms for
resolution of singularities by Bierstone-Milman [1] Villamoyor [1] Encinas-Villamayor [1] and
others take advantages of their internal inductional structures together with some ingenious
introduction of invariants measuring how singular your variety is.  Our hope here is that,
going away from the original factorization problem which only applies to the birational maps,
via the flexiblity of the toroidalization conjecture dealing with morphisms between toroidal
embeddings of any dimensions we should start seeing the internal inductional structure with
some invariants measuring how singular your morphism is. 

Though in dimension 2 we can check the toroidalization conjecture as in Abkulut-King [1]
Abramovich-Karu [1] Abramovich-Karu-Matsuki-W{\l}odarczyk [1] Karu [1] and Cutkosky-Piltant [1],
none of the existing methods fully reveals this (conjectural) inductional structure and hence
yields a proof even in dimension 3.  Probably the approach of Cutkosky-Piltant [1] is the
closest to our hope in spirit and the recent solution by Cutkosky [4] to the case where $\dim X
= 3$ and $\dim Y = 2$ should provide us more information toward the general solution.

Instead of starting from a morphism between nonsingular toroidal emebeddings, we could restate
the toroidalization conjecture for a proper morphism $\phi:X \rightarrow Y$ between any
nonsingular varieties claiming after appropriate sequences of blowups $\sigma_X:X' \rightarrow
X$ and $\sigma_Y:Y' \rightarrow Y$ the resulting morphism $\phi':(U_{X'},X') \rightarrow
(U_{Y'},Y')$ is toroidal for some appropriate choice of open subsets $U_{X'}$ and $U_{Y'}$.  In
this restated form, the toroidalization conjecture in the case when $\dim Y = 1$ is
(almost) equivalent to resolution of singularities of hypersurface singularities, revealing
another close connection between the problem of toroidalization and that of resolution of
singularities.  (Note that any proper morphism between NONSINGULAR toroidal embeddings stated in
Conjecture 6-2-1 is automatically toroidal in the case when $\dim Y = 1$ and hence trivial.)

\vskip.1in

Seeing a solution to the (weak) factorization problem of birational maps is hardly an end to
this exciting field surrounding resolution of singularities, toroidalization, factorization,
semistable reduction, log category ... etc. and we feel that we just had a sneak preview
of what is still to come.

\newpage

$$\bold{REFERENCES}$$

\vskip.2in

Abramovich, D. and De Jong, A.J.

\ \ \ [1] \it Smoothness, semistability, and toroidal geometry, \rm J. Alg. Geom.
$\bold{6}$ (1997), 

\ \ \ 789-801

\vskip.1in

Abramovich, D. and Karu, K.

\ \ \ [1] \it Weak semistable reduction in characteristic 0, \rm preprint (1997), 1-24

\vskip.1in

Abramovich, D., Karu, K., Matsuki, K. and W{\l}odarczyk, J.

\ \ \ [1] \it Torification and factorization of birational maps, \rm preprint (1999)

\vskip.1in

Abramovich, D., Matsuki, K. and Rashid, S.

\ \ \ [1] \it A note on the factorization theorem of toric birational maps after
Morelli 

\ \ \ and its toroidal extension, \rm to appear from Tohoku Mathematical
Journal 

\ \ \ (1999), 1-56

\vskip.1in

Abramovich, D. and Wang, J.

\ \ \ [1] \it Equivariant resolution of singularities in characteristic 0, \rm Math.
Res. 

\ \ \ Letters (1997)

\vskip.1in

Abhyankar, S.
 
\ \ \ [1] \it Local uniformization on algebraic surfaces over ground fields of

\ \ \ characteristic $p \neq 0$, \rm Annals of Math. $\bold{63}$ (1956), 491-526

\vskip.1in

Abkulut, S. and King, H.

\ \ \ [1] \it Topology of algebraic sets, \rm MSRI publications $\bold{25}$

\vskip.1in

Arapura, D. and Matsuki, K.

\ \ \ [1] \it Invariance of the Hodge numbers for smooth birational minimal models via

\ \ \ factorization theorem, \rm in preparation
(1999)

\vskip.1in

Batyrev, V.

\ \ \ [1] \it On the Betti numbers of birationally isomorphic projective varieties with 

\ \ \ trivial
canonical bundles, \rm alg-geom/9710020

\vskip.1in

Bierstone, E. and Milman, D.

\ \ \ [1] \it Canonical desingularization in characteristic zero by blowing up 

\ \ \ the
maximum strata of a local invariant, \rm Invent. Math. $\bold{128}$ (1997), 207-302

\vskip.1in

Bogomolov, F.A. and Pantev, T.G.

\ \ \ [1] \it Weak Hironaka Theorem, \rm Math. Res. Lett. $\bold{3}$ (1996), 299-308

\vskip.1in

Brion, M. and Processi, C.

\ \ \ [1] \it Action d'un tore dans une vari\'et\'e projective, \rm in Operator
algebras, unitary 

\ \ \ representations, enveloping algebras and invariant theory
(Paris, 1989), Birkh\"auser, 

\ \ \ Progr. Math. $\bold{92}$, 509-539

\vskip.1in

Christensen, C.

\ \ \ [1] \it Strong domination/weak factorization of three dimensional regular
local 

\ \ \ rings, \rm Journal of the Indian Math. Soc. $\bold{45}$ (1981), 21-47

\vskip.1in

Corti, A.

\ \ \ [1] \it Factorizing birational maps of 3-folds after Sarkisov, \rm J. Alg. Geom.
$\bold{4}$ 

\ \ \ (1995), 23-254

\vskip.1in

Crauder, B.

\ \ \ [1] \it Birational morphisms of smooth algebraic threefolds collapsing three

\ \ \ surfaces to a point, \rm Duke Math. J. $\bold{48}$ (1981), 589-632

\vskip.1in

Cutkosky, S.D.

\ \ \ [1] \it Local factorization of birational maps, \rm Advances in math.
$\bold{132}$ (1997), 

\ \ \ 167-315

\ \ \ [2] \it Local factorization and monomialization of morphisms, \rm preprint
(1997), 

\ \ \ 1-133

\ \ \ [3] \it Local factorization and monomialization of morphisms,
\rm math.AG/9803078, 

\ \ \ March 17 (1998), 1-141

\ \ \ [4] \it Monomialization for morphisms from 3-folds to surfaces, \rm a talk at 

\ \ \ Oberwolfach (1999)

\vskip.1in

Cutkosky, S.D. and Piltant, O.

\vskip.1in

\ \ \ [1] \it Monomial resolution of morphisms of algebraic surfaces, \rm preprint (1999)

\vskip.1in

Danilov, V.I. 

\ \ \ [1] \it The geometry of toric varieties, \rm Russian Math.
Surveys $\bold{33:2}$ (1978), 

\ \ \ 97-154

\ \ \ [2] \it The birational geometry of toric 3-folds, \rm Math. USSR-Izv.
$\bold{21}$ (1983), 

\ \ \ 269-280

\vskip.1in

De Concini, C. and Procesi, C.
 
\ \ \ [1] \it Complete Symmetric Varieties II, \rm in Algebraic Groups
and Related Topics 

\ \ \ (R. Hotta, ed.), Adv. Studies in Pure Math. $\bold{6}$
(1985), 481-513 

\vskip.1in

De Jong, A.J.

\ \ \ [1] \it Smoothness, semistability, and alterations, \rm Publ. Math. I.H.E.S.
$\bold{83}$ (1996), 

\ \ \ 51-93

\vskip.1in

Denef, J. and Loeser, F.

\ \ \ [1] \it Germs on singular algebraic varieties and motivic integration, \rm math.AG/9803039

\ \ \ (1998)

\vskip.1in

Dolgachev, I. and Hu, Y.

\ \ \ [1] \it Variantion of geometric invariant theory quotient, \rm Inst. Hautes
\'Etudes 

\ \ \ Sci. Publ. Math. $\bold{87}$ (1998), 5-56

\vskip.1in

Encinas, S. and Villamayor, O.

\ \ \ [1] \it A course on constructive desingularization and equivariance, \rm preprint (1999)

\vskip.1in

Ewald, G.

\ \ \ [1] \it Blow-ups of smooth toric 3-varieties, \rm Abh. Math. Sem. Univ. Hamburg

\ \ \ $\bold{57}$ (1987), 193-201

\vskip.1in

Frankel, T.

\ \ \ [1] \it Fixed points and torsion on K\"ahler manifolds, \rm Annals of Math.
$\bold{70}$, 1-8 

\ \ \ (1959)

\vskip.1in

Fulton, W.

\ \ \ [1] \it Intersection Theory, \rm Ergebnisse der Mathematik und ihrer Grenzgebiete

\ \ \ 3. Folge Band 2 (1984), Springer-Verlag

\ \ \ [2] \it Introduction to Toric varieties, \rm Annals of Mathematics studies
$\bold{131}$ (1993), 

\ \ \ Princeton University press

\vskip.1in

Gabrielov, A.M.

\ \ \ [1] \it The formal relations between analytic functions, \rm Funct. Anal. and Appl. 

\ \ \ $\bold{5}$
(1971), 318-319

\ \ \ [2] \it Formal relations among analytic functions, \rm Math. USSR Izv. $\bold{37}$ (1973),

\ \ \ 1056-1088

\vskip.1in

Grothendieck, A.

\ \ \ [1] \it Technique de Construction in Geometrie Analytique VI, \rm Etudes local des 

\ \ \ Morphismes.
Sem. H. Cartan, 13e anne 1960/61, exp. $\bold{13}$

\vskip.1in

Hartshorne, R.

\ \ \ [1] \it Algebraic Geometry, \rm Graduate Texts in Mathematics
$\bold{52}$ (1997), 

\ \ \ Springer-Verlag

\vskip.1in

Hironaka, H.

\ \ \ [1] \it On the theory of birational blowing-up, \rm Thesis (Harvard) (1960)

\ \ \ [2] \it Resolution of singularities of an algebraic variety over a field of

\ \ \ characteristic zero, \rm Annals of Math. $\bold{79}$ (1964), 109-326

\ \ \ [3] \it Flattening theorem in complex analytic geometry, \rm Amer. J. of Math.
$\bold{97}$ 

\ \ \ (1975), 503-547

\vskip.1in

Hu, Y.

\ \ \ [1] \it The geometry and topology of quotient varieties of torus actions,
\rm Duke 

\ \ \ Math. J. $\bold{68}$ (1992), 151-184: Erratum, Duke Math. J. $\bold{68}$
(1992), 609

\ \ \ [2] \it Relative geometric invariant theory and universal moduli spaces,
\rm Internat. 

\ \ \ J. Math. $\bold{7}$ (1996), 151-181

\vskip.1in

H\"ubl, R.

\ \ \ [1] \it Completions of local morphisms and valuations, \rm preprint (1999)

\vskip.1in

Iitaka, S.

\ \ \ [1] \it Algebraic Geometry (An Introduction to Birational geometry of
Algebraic 

\ \ \ Varieties), \rm Graduate Texts in Math. $\bold{76}$ (1982),
Springer-Verlag

\vskip.1in

Karu. K.

\ \ \ [1] \it Semistable reduction in characteristic zero, \rm Thesis (Boston University) 

\ \ \ (1998)

\vskip.1in

Kato, K.

\ \ \ [1] \it Toric singularities, \rm Amer. J. Math. $\bold{116}$ (1994), 1073-1099

\vskip.1in

Kawamata, Y.

\ \ \ [1] \it On the finiteness of generators of the pluri-canonical ring for a
threefold of 

\ \ \ general type, \rm Amer. J. Math. $\bold{106}$ (1984), 1503-1512

\ \ \ [2] \it The cone of curves of algebraic varieties, \rm Annals of Math.
$\bold{119}$ (1984), 

\ \ \ 603-633

\ \ \ [3] \it Crepant blowing-ups of three dimensional canonical singularities and
its 

\ \ \ application to degenerations of surfaces, \rm Annals of Math. $\bold{127}$
(1988), 93-163

\vskip.1in

Kawamata, Y., Matsuda, K. and Matsuki, K.

\ \ \ [1] \it Introduction to the minimal model problem, \rm Advanced Studies in Pure

\ \ \ Mathematics $\bold{10}$ (1987), 283-360

\vskip.1in

Kempf, G., Knudsen, F., Mumford, D. and Saint-Donat, B.

\ \ \ [1] \it Toroidal embeddings I, Lecture Notes in Mathematics $\bold{339}$
(1973), 

\ \ \ \rm Springer-Verlag

\vskip.1in

King, H.

\ \ \ [1] \it Resolving Singularities of Maps, Real algebraic geometry and topology

\ \ \ (East Lansing, MI 1993) (1995), \rm Contemp. Math., Amer. Math. Soc.

\ \ \ [2] \it A private e-mail to Morelli (1996) \rm

\vskip.1in

Kirwan, F.

\ \ \ [1] \it Cohomology of quotients in symplectic and algebraic geometry, Notes
$\bold{31}$ 

\ \ \ (1984), \rm Princeton University Press

\vskip.1in

Kleiman, S.

\ \ \ [1] \it Toward a numerical theory of ampleness, \rm Annals of Math.
$\bold{84}$ (1966), 

\ \ \ 293-344

\vskip.1in

Koll\'ar, J.

\ \ \ [1] \it The cone theorem.  Note to a paper ``The cone of curves of algebraic

\ \ \ varieties" by Y. Kawamata, \rm Annals of Math. $\bold{120}$ (1984), 1-5

\ \ \ [2] \it Flops, \rm Nagoya Math. J. $\bold{113}$ (1989), 15-36

\vskip.1in

Kulikov, V.S.

\ \ \ [1] \it Decomposition of a birational map of three-dimensional varieties
outside 

\ \ \ codimension 2, \rm Math. USSR Izv. $\bold{21}$ (1983), 187-200

\vskip.1in

Matsuki, K.

\ \ \ [1] \it Introduction to the Mori program, \rm the manuscript of a textbook to
appear 

\ \ \ from Springer-Verlag (1998) 

\vskip.1in

Matsuki, K. and Rashid, S.

\ \ \ [1] \it An easy and inductive proof for Sumihiro's equivariant completion theorem 

\ \ \ for toric
varieties, \rm in preparation (2000)

\newpage

Matsuki, K. and Wentworth, R.

\ \ \ [1] \it Mumford-Thaddeus principle on the moduli space of vector bundles on
an 

\ \ \ algebraic surface, \rm Internat. J. Math. $\bold{8}$ (1997), 97-148

\vskip.1in

Matsumura, H.

\ \ \ [1] \it Commutative ring theory, \rm Cambridge Studies in Advanced Mathematics $\bold{8}$ 

\ \ \ (1930)

\vskip.1in

Milnor, J.

\ \ \ [1] \it Morse Theory, Annals of Math. Stud. $\bold{51}$ (1963), \rm Princeton
Univ. Press

\vskip.1in

Moishezon, B.

\ \ \ [1] \it On $n$-dimensional compact varieties with $n$ algebraically
independent 

\ \ \ meromorphic functions, \rm Amer. Math. Soc. Trans. $\bold{63}$
(1967), 51-177

\vskip.1in

Morelli, R.

\ \ \ [1] \it The birational geometry of toric varieties, \rm J. Alg. Geom.
$\bold{5}$ (1996), 751-782

\ \ \ [2] \it Correction to ``The birational geometry of toric varieties", \rm homepage at

\ \ \ the Univ. of Utah (1997), 767-770

\vskip.1in

Mori, S.

\ \ \ [1] \it Projective manifolds with ample tangent bundles, \rm Annals of Math.
$\bold{110}$ 

\ \ \ (1986), 593-606

\ \ \ [2] \it Threefolds whose canonical bundles are not numerically
effective, \rm Annals 

\ \ \ of Math. $\bold{116}$ (1982), 133-176

\ \ \ [3] \it Flip theorem and the existence of minimal models for 3-folds, \rm J. Amer.

\ \ \ Math. Soc. $\bold{1}$ (1988), 117-253

\vskip.1in

Mumford, D.

\ \ \ [1] \it Abelian Varieties, 2nd Edition (1974), \rm Oxford University Press

\vskip.1in

Mumford, D., Fogarty J. and Kirwan, F.

\ \ \ [1] \it Geometric Invarinat Theory (Third Enlarged Edition), \rm Ergebnisse der

\ \ \ Mathematik und ihrer Grenzgebiete $\bold{34}$ (1992), Springer-Verlag

\vskip.1in

Oda, T.
 
\ \ \ [1] \it Lectures on Torus Embeddings and Applications (Based on joint
work with 

\ \ \ K. Miyake), \rm Tata Inst. of Fund. Research $\bold{58}$ (1978),
Springer-Verlag

\ \ \ [2] \it Convex Bodies and Algebraic Geometry (An Introduction to the Theory
of 

\ \ \ Toric Varieties), \rm Ergebnisse der Mathematik und ihrer Grenzgebiete 3.
Folge 

\ \ \ Band 15 (1988), Springer-Verlag 

\vskip.1in

Pandripande, R.

\ \ \ [1] \it A compactification over $\overline{M_g}$ of the universal moduli space
of 

\ \ \ slope-semistable vector bundles, \rm J. Amer. Math. Soc. $\bold{9}$ (1996), 425-471

\vskip.1in

Paranjape, K.H.

\ \ \ [1] \it Bogomolov-Pantev resolution-an expository account, \rm in New Trends in 

\ \ \ Algebraic
Geometry, Warwick July 1996, Cambridge University Press (1998)

\newpage

Pinkham, H.

\ \ \ [1] \it Factorization of birational maps in dimension 3, \rm Proceedings of Symposia

\ \ \ in Pure Math. $\bold{40}$ (1983)

\vskip.1in

Reid, M.

\ \ \ [1] \it Canonical threeefolds, \rm G\'eom\'etrie Alg\'em\'etrie Alg\'ebrique Angers
1979, 

\ \ \ A. Beauville ed. (1980), Sijthoff and Nordhoff, 273-310

\ \ \ [2] \it Minimal models of canonical 3-folds, \rm Adv. Stud. in Pure Math.
$\bold{1}$ (1983), 

\ \ \ 131-180

\ \ \ [3] \it Decomposition of Toric Morphisms, \rm Arithmetic and Geometry, papers

\ \ \ dedicated to I. R. Shafarevich on the occasion of his 60th birthday, vol. II,

\ \ \ Progress in Math. (M. Artin and J. Tate eds.) $\bold{36}$ (1983), 395-418

\ \ \ [4] \it Birational geometry of 3-folds according to Sarkisov, \rm preprint (1991)

\vskip.1in

Sally, J.

\ \ \ [1] \it Regular overrings of regular local rings, \rm Trans. Amer. math. Soc.
$\bold{171}$ 

\ \ \ (1972), 291-300

\vskip.1in

Sarkisov, V.G.

\ \ \ [1] \it Birational maps of standard ${\Bbb Q}$-Fano fiberings, \rm I. V. Kurchatov
Institute 

\ \ \ Atomic Energy preprint (1989)

\vskip.1in

Schaps, M.

\ \ \ [1] \it Birational morphisms of smooth threefolds collapsing three surfaces to

\ \ \ a curve, \rm Duke Math. J. $\bold{48}$ (1981), 401-420

\vskip.1in

Shafarevich, I.R.

\ \ \ [1] \it Basic Algebraic Geometry (1977), \rm Springer-Verlag

\vskip.1in

Shannon, D.L.

\ \ \ [1] \it Monoidal transforms, \rm Amer. J. Math. $\bold{45}$ (1973), 284-320

\vskip.1in

Shokurov, V.

\ \ \ [1] \it The non-vanishing theorem, \rm Math. USSR Izv. $\bold{19}$ (1985), 591-604

\vskip.1in

Sumihiro, H.

\ \ \ [1] \it Equivariant Completion I, \rm J. Math. Kyoto Univ. $\bold{14}$ (1974), 1-28

\ \ \ [2] \it Equivariant Completion II, \rm J. Math. Kyoto Univ. $\bold{15}$ (1975),
573-605

\vskip.1in

Thaddeus, M.

\ \ \ [1] \it Stable pairs, linear systems and the Verlinde formula, \rm Invent. Math.
$\bold{117}$ 

\ \ \ (1994), 317-353

\ \ \ [2] \it Geometric invariant theory and flips, \rm J. Amer. Math. Soc. $\bold{9}$
(1996), 

\ \ \ 691-723

\vskip.1in

Villamayor, O.

\ \ \ [1] \it Constructiveness of Hironaka's resolution, \rm Ann. Sci. \'Ecole Norm. Sup. 

\ \ \ (4) $\bold{22}$
no.1 (1989), 1-32

\vskip.1in

Wang, C.-L.

\ \ \ [1] \it On the topology of birational minimal models, \rm math.AG/9804050

\newpage

W{\l}odarczyk, J.

\ \ \ [1] \it Decomposition of Birational Toric Maps in Blow-Ups and Blow-Downs.  

\ \ \ A
Proof of the Weak Oda Conjecture, \rm Transactions of the AMS $\bold{349}$
(1997), 

\ \ \ 373-411

\ \ \ [2] \it Birational cobordism and factorization of birational maps, \rm to appear
from 

\ \ \ J. of Alg. Geom. (1998)

\ \ \ [3] \it Combinatorial structures on toroidal varieties and a proof of the weak

\ \ \ factorization theorem, \rm preprint math.AG/9904076 (1999)

\vskip.1in

Zariski, O.

\ \ \ [1] \it Algebraic Surfaces, (1934), \rm Springer-Verlag

\enddocument